\documentclass[12pt, notitlepage]{article}
\usepackage{amsfonts, amsmath, amssymb, amsthm, color, hyperref, enumitem, fancyhdr, relsize, tikz, graphicx, mathrsfs}
\usepackage[utf8]{inputenc}
\usepackage{lipsum}

\title{\scshape 
On Hastings' approach to Lin's Theorem for Almost Commuting Matrices \sffamily }
\author{\scshape \sffamily David Herrera\footnote{Rutgers University. dh708@math.rutgers.edu}
}
\date{\today}

\newtheorem{prop}{Proposition}[section]
\newtheorem{corollary}[prop]{Corollary}
\newtheorem{thm}[prop]{Theorem}
\newtheorem{lemma}[prop]{Lemma}

\theoremstyle{definition}
\newtheorem{remark}[prop]{Remark}
\newtheorem{example}[prop]{Example}

\newtheorem{defn}[prop]{Definition}

\newcommand{\R}{\mathbb{R}}
\newcommand{\C}{\mathbb{C}}
\newcommand{\N}{\mathbb{N}}

\newcommand{\Z}{\mathbb{Z}}

\renewcommand{\Im}{\operatorname{Im}}
\renewcommand{\Re}{\operatorname{Re}}
\newcommand{\rank}{\operatorname{rank}}

\newcommand{\diam}{\operatorname{diam}}
\newcommand{\dist}{\operatorname{dist}}

\begin{document}


\maketitle

\abstract{Lin's theorem states that for all $\epsilon > 0$, there is a $\delta > 0$ such that for all $n \geq 1$ if self-adjoint contractions $A,B \in M_n(\C)$ satisfy $\|[A,B]\|< \delta$ then there are self-adjoint contractions $A',B' \in M_n(\C)$ with $[A',B']=0$ and $\|A-A'\|,\|B-B'\|<\epsilon$. We present fully explained and corrected details of the approach in \cite{Hastings}, which was the first version of Lin's theorem  to provide asymptotic estimates. 

We also apply this method to the case where $B$ is a normal matrix with spectrum lying in some nice 1-dimensional subset of $\C$.}

\section{Introduction}\label{Intro}
The following is known as Lin's theorem  
\begin{thm}
If $A, B \in M_n(\C)$ are self-adjoint contractions with \[\|[A,B]\| < \delta,\] then there are commuting self-adjoint $A', B'\in M_n(\C)$ such that \[\|A - A'\|, \|B-B'\| < \epsilon(\delta),\]
where
$\lim_{\delta\to0^+}\epsilon(\delta)=0$.
\end{thm}
It is important to note that $\epsilon(\delta)$ does not depend on $n$. We will discuss results that depend on $n$ later.
We formulated Lin's theorem similar to \cite{HastingsLoring}, \cite{KS}. One can also express this result in terms of $\delta$ depending on $\epsilon$ as in the abstract.

If $\|[A,B]\|$ is small, we refer to $A$ and $B$ as ``almost commuting'' and if there are nearby commuting matrices (with whatever specified properties), we refer to $A$ and $B$ as ``nearly commuting''. In these terms, Lin's theorem says that almost commuting Hermitian matrices (with bounded norm) are near Hermitian commuting matrices. Lin's theorem can also be formulated in terms of almost normal matrices being nearly normal, where a matrix $N$ is almost normal if its self-commutator $[N^\ast ,N]$ has small norm and $N$ is nearly normal if there is a normal matrix nearby $N$.

Lin's theorem is the positive answer to a question (for the operator norm) stated by Rosenthal in \cite{Rosenthal} in\footnote{The dates stated in this paper use published dates and received dates when possible.} 1969 and in 1976 was listed as one of various unsolved problems by Halmos in \cite{Halmosproblems}. 
By the time that Lin's proof appeared, there were various results that fell short of Lin's theorem in important ways. 
Positive results were obtained when $\epsilon$ is allowed to depend on $n$ and negative results were obtained when we are considering general self-adjoint operators acting on an infinite dimensional Hilbert space:
\cite{Lux} in 1969 (with the Euclidean norm) and \cite{Subnormal} in 1974 (with the operator norm, with proof attributed to Halmos) presented a result with dimensional dependence of $\delta$. 
Also, \cite{Subnormal} presented related results that fail in infinite dimensions and \cite{berg-olsen} in 1980 presented an example of an obstruction related to the Fredholm index in the infinite dimensional analogue to Lin's theorem. 
In 1979, Pearcy and Shields in \cite{P-S} gave the dimensional dependence of $ \epsilon=((n-1)\delta/2)^{1/2}$, but only used that one of the operators is assumed to be self-adjoint. The dependence of $\epsilon$ on $\delta$ has the optimal exponent for Lin's theorem.

Davidson in \cite{Davidson} in 1985 provided an example of normal matrices (gotten by an approximation result of Berg in \cite{Berg}) and a Hermitian matrices that almost commute but are not nearby commuting normal and Hermitian matrices using a spectral projection argument and a property of the shift on $\C^n$. This shows that the Pearcy and Shields result cannot be improved to remove its dimensional dependence, as we will discuss later. 
See Section \ref{A and N} for more about generalizations.

In 1990, Szarek in \cite{Szarek} proved that Lin's theorem holds for $\epsilon=c(n^{1/2}\delta)^{2/13}$, where $c$ is some constant. Szarek states that the exponent of $\epsilon$ can be reduced by using a more complicated argument but the aim was to use the approach discussed in \cite{Davidson} and the assumption that both operators are Hermitian to obtain stronger dimensional dependence than the counter-examples seen thus far had exhibited and hence ``the situation in the Hermitian setting is completely different''.

In 1997, Lin in \cite{Lin} provided a proof of Lin's Theorem and in 1996 Friis and R{\o}rdam in \cite{F-R} published a simplified version of Lin's result (referencing \cite{Lin} which was to appear in Operator Algebras and Their Applications).

Up until this point, the proof of Lin's theorem was abstract and neither gave the asymptotic dependence of $\epsilon$ solely on $\delta$ nor gave a way of constructing these matrices. 
In 2008, \cite{Hastings1sted} by Hastings presented an approach whose starting point is similar to one of \cite{Davidson}'s reformulations of Lin's theorem. \cite{Hastings1sted} claimed to provide a constructive method of finding the nearby commuting matrices and also presenting an asymptotic dependence of $\epsilon$ on $\delta$. 
The methods in \cite{Hastings1sted} for a diagonal and tridiagonal pair of Hermitian matrices (which was already implicitly solved without the ideal exponent in \cite{Szarek}) held. However, the arguments for the general case were amended to try to resolve errors and consequently the main result changed. 
Namely, the most recent version of this paper, being \cite{Hastings} pre-published on arxiv.org in 2010, states similar statements for asymptotic dependence of $\epsilon$ on $\delta$ but not that the result is constructive and the dependence only holds for $\delta$ small enough (but neither how small nor any constant in the dependence of $\epsilon$ on $\delta$ is given).

In 2015, Kachkovskiy and Safarov in \cite{KS} proved a result that not only gives the optimal homogeneous dependence of $\epsilon$ on $\delta$ of $\epsilon= Const.\delta^{1/2}$ and simultaneously addresses the infinite dimensional situation where the index obstruction of Berg and Olsen applies. The proof also appears to be essentially constructive if the nearby commuting matrices sought for belong to the von Neumann algebra generated by the two almost commuting matrices.

A list of applications of almost/nearly commuting matrices can be found in the introduction of \cite{LS}. The first example given there is related to a paper of von Neumann concerning macroscopic observables. Although at least three almost commuting Hermitian matrices are not necessarily nearly commuting, Ogata in \cite{Ogata} showed that a special case of this related to von Neumann's example is true. 

There have also been many ``spin-offs'' of this problem, involving different norms, algebraic objects that are not operators on a Hilbert space, and for bounded operators on an infinite dimensional Hilbert space. Notable mentions are that Lin's theorem holds for the normalized Schatten-$p$ norms (\cite{Said}) and that Lin's theorem for the normalized Hilbert-Schmidt norm holds for multiple Hermitian matrices (\cite{Glebsky}). 

\vspace{0.05in}

Given that this problem has effectively been solved in the general case, one might wonder why would one explore Hastings' earlier approach that does not provide the optimal exponent and does not provide a constructive proof unless relying on \cite{KS} as a black-box. One reason is that it demonstrates how to obtain an asymptotic estimate from even the version of Lin's theorem stated without any explicit estimate provided. That is, some of the arguments in Hastings' approach are interesting in and of themselves.

A main reason that the author still finds Hastings' approach interesting is that the unlike the approach of \cite{KS}, Hastings' approach involves following one of Davidson's reformulations of Lin's theorem in terms of constructing a certain almost invariant subspace of a block tridiagonal matrix. Szarek seems to be the first to effectively use this reformulation and Hastings provided a nice proof of this lemma in the case of tridiagonal matrices in \cite{Hastings}. Although Lin's theorem has been proved with the optimal asymptotic exponent, it is still an open problem whether the exponent in this reformulation of Lin's theorem can be improved. See Remark \ref{Linexpbootstrap} below. So, we still have an interesting matrix theory question to explore. 

Also, even in the case that one uses \cite{KS}'s result, the structure gotten by Davidson's reformulation of Lin's theorem is nice enough that one might find reason to make a compromise and apply \cite{KS}'s result to Davidson's reformulation. For instance, we explore some of this in Section \ref{A and N}.

As a final reason, this reformulation of Lin's theorem also recently led to a constructive proof of Ogata's theorem for $d=2$ in the soon to be announced \cite{OgataC}.
Because finding an explicit estimate for the problem of two almost commuting unitary matrices is still open, one also might hope to find inspiration for a solution in this reformulation of Lin's theorem.

\vspace{.1in}

\noindent\underline{Purpose of Paper}: The goal of this paper is to present a full and correct account of Hastings' approach to Lin's theorem in \cite{Hastings}, which is the second revision since the result was published (\cite{Hastings1sted}).
In this aim, we factor out certain claims made in \cite{Hastings}, providing references or proofs for clarity. The arguments and calculations in the proof of Hastings' result are intended to be in full (and perhaps too much) detail.
For various modifications to \cite{Hastings}, the author is indebted to Hastings for clarifications and resolutions of a number of the errors found by the author. This is cited as \cite{PC}.

\vspace{0.1in}

\noindent\underline{Changes in the second edition of this paper}: The title of the paper was changed from ``Constructive Approaches to Lin's Theorem for Almost Commuting Matrices'' to better align with the above stated purpose of the paper.
Various comments and discussions that were not central to the paper were removed. The introductory sections have been condensed and  rewritten. 

Typos were corrected and various explanations in the proof of Hastings' result were slightly rewritten or expanded for clarity. 
Comments in the abstract and introduction concerning the constructiveness of \cite{KS} were updated. Inequalities in  Remark \ref{linindep} and Proposition \ref{Gap} were corrected. 
No substantial changes to the proof of Hastings' result were made.

\vspace{0.1in}

\noindent \underline{Notation Conventions}: Hilbert spaces are always finite dimensional and complex with an inner product $(-,-)$ that is conjugate linear in the first argument. $\|-\|$ always denotes the operator norm of a linear transformation on a Hilbert space (i.e. a matrix), $|-|$ always denotes the norm of a vector (or absolute value of a number), and $[C,D] = CD - DC$ is the commutator of two matrices.  It is common in the literature on this problem to see self-adjoint matrices referred to as Hermitian. Often, this will be done in this paper when discussing the statements of Lin's theorem or results from other papers. However, in the bulk of the paper, we will simply refer to them as ``self-adjoint''.

\section{Outline of Hastings' proof}\label{Outline}
One of the equivalent forms of Lin's theorem that Davidson uses is the following (we use notation from \cite{Hastings}):

($Q'$): For every $\epsilon > 0$, there is a $L_0 > 0$ with the following property. For any self-adjoint contraction $J$ that is block tridiagonal with $L \geq L_0$ blocks, there is a projection $P$ that contains the first block of the basis, is orthogonal to the last block of the basis, and satisfies $\|[J,P]\| < \epsilon.$

Moreover, a constructive proof of Lin's theorem gives a constructive proof of ($Q'$) and vice-versa. See Section \ref{reformulations} for more about this.

The key result of \cite{Hastings} that we will prove is of the form:
\begin{thm}\label{maintheorem}
For $\delta > 0$ small enough there exists an $\epsilon = \epsilon(\delta) > 0$ such that for any self-adjoint contractions $A, B \in M_n(\C)$ such that $\|[A,B]\| \leq \delta$ there exist self-adjoint $A', B'$ such that $[A',B'] = 0$ and $\|A-A'\|, \|B-B'\| \leq \epsilon.$ We can choose $\epsilon(\delta) = E(1/\delta)\delta^\beta$, where $\beta = 1/6$ and $E$ is a function growing slower than any polynomial.

\end{thm}
The value of $\beta$ gotten in \cite{Hastings} depends on the type of matrix considered. Our exposition will prove this result with $\beta = 1/6$ without condition on the structure of $A$ and $B$. This is the exponent of the original published article \cite{Hastings1sted}, but smaller than the proposed update in \cite{Hastings} given some complications from its Lemma 4. Hastings in \cite{Hastings1sted}  proves this result with the optimal exponent $\beta = 1/2$ and a simple construction when $B$ is diagonal and $A$ is tridiagonal.

We now proceed to the discussion of the proof of Theorem  \ref{maintheorem}. We start with two almost commuting self-adjoint contractions $A, B \in M_n(\C)$. We apply Lemma \ref{finite range} to replace $A$ with $H$. $H$ and $B$ are almost commuting self-adjoint contractions and $H$ has the additional property that it is ``finite range'' (we can choose how far the ``range'' $\Delta$ is) with respect to the eigenspaces of $B$ by which we mean that it maps eigenspaces of $B$ into the span nearby eigenspaces. The price paid is that the distance between $A$ and $H$ will depend on the commutator $[A,B]$ of $A$ and $B$ and on the range $\Delta$.

Next, break the spectrum of $B$ up into many small disjoint intervals $I_i$. We choose the ``finite range'' to be much smaller than the length of these intervals, so that each of these intervals breaks up into many subintervals of length at least $\Delta$ and at most $2\Delta$. 
We form orthogonal subspaces by grouping the eigenvectors of $B$ associated with each subinterval. 
We project $H$ onto the eigenvectors of $B$ associated with the interval $I_i$ (the range of $E_{I_i}(B)$). Then on the range of each $E_{I_i}(B)$ the projection of $H$ will be block tridiagonal with respect to the many subspaces associated with the subintervals.

We apply ($Q'$), stated as Lemma \ref{modified lemma}, here to then obtain for each interval $I_i$ a subspace $\mathcal W_i$ that is almost invariant under the restriction of $H$, contains the eigenvectors associated with the first subinterval, and is orthogonal to the eigenvectors associated with the last subinterval. 
Each subspace has an orthogonal complement as a subspace of $R(E_{I_i}(B))$ which we call $\mathcal W_i^\perp$. 
Now, we use these to form $A'$ by projecting $H$ onto the subspaces $\mathcal W_i^\perp\oplus \mathcal W_{i+1}$. 
We define $B'$ to be a multiple of the identity on $\mathcal W_i^\perp\oplus \mathcal W_{i+1}$. The details can be found in Section \ref{reformulations}. 
The estimates for these various steps are collected in Section \ref{estimates}.

At this point, we would focus all of our attention on proving Lemma \ref{main lemma}, which is a formulation of $(Q')$. After a simple projection construction used in \cite{Szarek}, this implies Lemma \ref{modified lemma}. The set-up that we have is a self-adjoint contraction which is (block) tridiagonal with respect to some orthogonal subspaces $\mathcal V_1, \dots, \mathcal V_{L}$. In Section \ref{General Approach} there is a discussion of the general approach from both \cite{Hastings} and \cite{Szarek}.

A starting point is that the first two parts of Lemma \ref{main lemma} can be satisfied using a simple construction: break the spectrum of $J$ up into intervals $I_i$ (which are unrelated to the intervals discussed above), then constructions similar to $\tilde{\mathcal W}=\bigoplus\tilde{\mathcal X_i}$ with $\tilde{\mathcal X_i}= \chi_{I_i}(J)(\mathcal V_1$) are used to construct the subspace $\mathcal W$. The idea is that this subspace contains $\mathcal V_1$ as a subspace and because each subspace is almost invariant under $J$ and they are orthogonal, the entire subspace is almost invariant under $J$. The key obstacle to be avoided is that we need $\tilde{\mathcal W}$ to be almost orthogonal to $\mathcal V_L$. The details and further motivation are in Section \ref{General Approach}. 

Section \ref{Dim section} is a break from the ideas behind Hastings' proof to discuss the results with dimensional dependence.
Section \ref{smooth partitions} lists some of the properties of the smooth cut-off functions ${\mathcal F }_{\omega_0}^{r,w}(\omega)$ used and provides a definition of the key subspaces $\mathcal X_i$. 

Section \ref{use of Lin's theorem} singles out Lemma \ref{application of Lin's theorem}, the ``nonconstructive bottleneck'' of Hastings' proof. We only use this lemma for one value of $\epsilon < 1$ in the construction of the $\mathcal N_i$ to ensure that $\|N_{i-1}N_iN_{i+1}\| \leq 1-\chi < 1$ for some $\chi > 0$. This involves applying Lin's theorem for only one value of $\epsilon$ less than $1/22$.  As described in \cite{Hastings}, Hastings' result bootstraps Lin's theorem for one value of $\epsilon$ to get a result for an asymptotic dependence of $\delta$ on $\epsilon$. 

Section \ref{Many subspaces} discusses many of the subspaces derived from the $\mathcal X_i$ and the first three paragraphs contains some of the general motivation in the construction of $\mathcal W$. An issue is that if $x = \sum x_i$ for $x_i \in \mathcal X_i$ then we do not clearly (or even necessarily) have control of $\sum |x_i|^2$ in terms of $|x|$ due to having no control of the orthogonality of the $\mathcal X_i$. To approach this, we explore the subspace $\mathcal X = \sum_i \mathcal X_i$ and representation space $\mathcal R$ which is the exterior direct sum of the $\mathcal X_i$. We form the linear map $A:\mathcal R \to \mathcal X$ that is the identity on each $\mathcal X_i$ (which we call $\mathcal R_i$ when it is a subspace of $\mathcal R$) and we are interested in having some control of $\sum |x_i|^2=\sum_i|r_i|^2$ in terms of $|x|^2=|A\sum_i r_i|^2$ by restricting $A$ to a subspace. 

In Section \ref{N_i} we form the subspace $\mathcal U^\perp$ from subspaces $\mathcal N_i$ ($i$ even) and $\mathcal N_i'$ ($i$ odd) which satisfy certain properties and set $\mathcal W = A(\mathcal U)$. Roughly speaking, the localized subspaces $\mathcal N_i$ belong to the span of the eigenvectors of $\rho = {A^\ast}A$ with ``small eigenvalues''  and can be put together to approximately recover in a reasonable way the eigenvectors of $\rho$ with small eigenvalues. This uses the tridiagonal nature of $\rho$, Lemma \ref{operator Lieb-Robinson modified}, and Lemma \ref{application of Lin's theorem}. For $i$ odd, we cut out the part of $\mathcal N_i$  that is not orthogonal enough to the neighboring (even) $\mathcal N_{i-1}$ and $\mathcal N_{i+1}$ to obtain $\mathcal N_i'$. This semi-orthogonality is used throughout the remainder of the proof.

In Section \ref{U and W} some various inequalities concerning representation of vectors in $\mathcal W$ are proved along with inequalities concerning the projection onto $\mathcal W$. Finally, in Section \ref{Verifying} we conclude the proof with verifying the validity of the desired properties of $\mathcal W$.

Section \ref{LA} (Linear Algebra Preliminaries), Section \ref{Spectral Projections} (Lemmas on Spectral Projections and Commutators), and Section \ref{Lieb-Robinson section} (Relevant Lieb-Robinson Bounds) provide various important results that seem outside the ``main story'' and have been separated from the rest of the proof for clarity. 

Section \ref{A and N} contains some consequences of Hastings' result that rely upon the fact that it uses Davidson's reformulation of Lin's theorem.

\section{Linear Algebra Preliminaries}\label{LA}
Note that $|-|$ will always denote the (Hilbert space) norm of a vector and $\|-\|$ will always denote the operator norm of a matrix. If $S$ is a set and $N$ is normal, $E_S(N)$ is the spectral projection of $N$ on $S$. If $W$ is a subspace, $P_W$ is a projection onto $W$. Note that $S_1, S_2\subset \R$ will always denote disjoint sets where the relevant constant is the distance between them and $S'', S'\subset \R$ denote nested sets and the relevant constant is the distance between $S''$ and $\R \setminus S'$. We will use the convention for the Fourier transform: $\hat{f}(k) = \frac{1}{2\pi}\int f(x)e^{-ikx}dx$.

\begin{defn}
If we have a sequence of vectors $u_k$ (resp.  subspaces $\mathcal U_k$) such that if $i, j$ with $|i-j|\geq 2$ then $u_i$ and $u_j$ (resp. $\mathcal U_i$ and $\mathcal U_j$) are orthogonal, we call this sequence nonconsecutively orthogonal.
\end{defn}

One type of estimate used often in the proof of \cite{Hastings}'s Lemma 2 is that if $x_1, \dots, x_n$ are nonconsecutively orthogonal vectors with sum $v$ then we have
\begin{align} \label{constant 2}
|v|^2 = \left|\sum_{i} x_i\right|^2 = \left|\sum_{i\,  odd} x_i + \sum_{i\, even} x_i\right|^2 \leq 2 \left|\sum_{i\, odd} x_i\right|^2 + 2\left|\sum_{i\, even} x_i\right|^2 = 2\sum_{i}|x_i|^2.
\end{align}

If one has uniform control of the inner products of vectors $x_i$ then one has a reverse inequality. Specifically, if there is also $C < 1/2$ such that 
\begin{align}\label{degree of orthogonality}
|(x_i, x_j)| \leq C|x_i||x_j|
\end{align}
for $j=i\pm 1$, we have 
\begin{align} \label{semi-pythagorean}
\nonumber |v|^2 &= |\sum_{i} x_i|^2 \geq \sum_{i} |x_i|^2 - 2\sum_{i<j} |(x_i, x_j)| \geq \sum_{i} |x_i|^2 - 2C\sum_{1\leq i< n} |x_i||x_{i+1}| \\
&\geq \sum_{i} |x_i|^2 - 2C\sum_{1\leq i< n} \frac{|x_i|^2+|x_{i+1}|^2}{2} \geq (1-2C)\sum_{i} |x_i|^2.
\end{align}
This condition of some degree of minimal orthogonality is what makes Hastings' tridiagonal result possible, because given any nonconsecutively orthogonal subspaces $\mathcal X_k$ which are each one dimensional, either we have control on the inner products as in Equation (\ref{degree of orthogonality}) for each vector $x_i\in\mathcal X_i, x_{i+1}\in\mathcal X_{i+1}$ or we have the opposite inequality. This is the key element to the construction in Lemma 6 of \cite{Hastings}.

We consider Proposition 2.2 from \cite{Exp} which after some crude estimates, gives the next result. Recall that a matrix $A$ is called $m$-banded if $A_{i,j} = 0$ for $|i-j| > m/2.$ Being $2$-banded is equivalent to being tridiagonal.
\begin{prop}\label{exponential decay of inverse}
Let $A$ be a tridiagonal, strictly positive definite matrix with $a = \min \sigma(A) > 0, b = \max \sigma(A).$ Then 
\[|(A^{-1})_{i,j}| \leq C\alpha^{|i-j|},\]
where $C = C(a, b)>0, \alpha = \alpha(a,b) \in (0,1)$.
\end{prop}
\begin{remark}\label{Benzi_Remark}
Examples of results like this for analytic functions, instead of the function $f(x) = 1/x$, of banded self-adjoint matrices can be found in \cite{Benzi}.
\end{remark}

The following lemma is a part of the proof of Lemma 5 of \cite{Hastings} which contains a sketch of this result as a claim.
\begin{lemma} \label{positive} Suppose that $c_i, d_i$, $i = 1, \dots, n$ are non-negative constants. If $M$ is a self-adjoint $n \times n$ tridiagonal matrix such that $M_{i,i} \geq c_i^2 + d_i^2$ for $i \leq n$ and $|M_{i,i+1}| \leq d_ic_{i+1}$ for $1\leq i < n$, then $M$ is positive.
\end{lemma}
\begin{proof}
We will compare it to the Hermitian tridiagonal matrix $D$ defined by $D_{i,i} = a_i^2 + b_i^2$ and $D_{i,i+1} = b_{i}\overline{a_{i+1}}$. Here is the case $n = 4$:
\[D = \begin{pmatrix} 
|a_1|^2 + |b_1|^2 & \overline{b_1}a_2 &  &    \\
b_1\overline{a_2} & |a_2|^2 + |b_2|^2 & \overline{b_2}a_3 &    \\
& b_2\overline{a_3} & |a_3|^2 + |b_3|^2 & \overline{b_3}a_4  \\
&  & b_3\overline{a_4} & |a_4|^2 + |b_4|^2
\end{pmatrix}.\]
This matrix is positive because if we consider the matrix $G$ defined with columns \\
$(a_1, b_1,0, \dots, 0)^T, (0, a_2, b_2, \dots,0)^T, \dots, (0, \dots, 0, a_n)^T$ then we have that $G^\ast G + b_n^2e^{n,n} = D$, where $(e^{n,n})_{i,j} = \delta^n_i\delta^n_j$. So, $D$ is positive.

When $a_i = c_i, b_i = d_i$, we see that the matrix $M$ is compared to $D$ by having its diagonal entries larger than those of $D$ and its off-diagonal entries smaller in absolute value than $D$. We now pick $a_i, b_i$ so that it is clear that $M$ is positive.
Let $a_i=c_i$ for $i=1,\dots,n$. If $d_ic_{i+1}=0$ then let $b_i = 0$. If $d_ic_{i+1}>0$, let $b_i = M_{i,i+1}/\overline{c_{i+1}}$ for $i=1,\dots,n-1$. This is so that $M_{i,i+1}=b_i\overline{a_{i+1}}$. Let $b_n=d_n$.

Now, $|b_i\overline{a_{i+1}}|\leq d_ic_{i+1}$ so $a_i=c_i$ and $|b_i| \leq d_i$. 
$|a_i|^2+|b_i|^2 \leq c_i^2+d_i^2 \leq M_{i,i}$.
This tells us that $D$ has the same off-diagonal terms as $M$ and its diagonal terms are less than the diagonal terms of $M$. This tells us that $M-D$ is positive, so $M$ is positive.

\end{proof}

We use Lemma 2 of \cite{P-S} concerning the Schur product
\begin{lemma}\label{P-S}
Let $(T_{i,j})$ be a matrix and let $a_i, b_j$ be real numbers with $a_i - b_j \geq d$. Then
\[\left\|\left(\frac{1}{a_i-b_j}T_{i,j}\right)\right\| \leq \frac{1}{d}\|(T_{i,j})\|.\]
\end{lemma}

The following (constructive) lemma which is Lemma 2.2 from \cite{Davidson} will serve many uses, including simplifying the statement of a main lemma in \cite{Hastings} and the lemma where Lin's theorem is applied. 
\begin{lemma}\label{projection lemma}
Let $E\leq G$ be projections on a Hilbert space and $\epsilon  > 0$. If $F'$ is a projection with $\|E{F'}^{\perp}\| < \epsilon$ and $ \|F'G^{\perp}\| < \epsilon$, then there is a projection $F$ such that $E \leq F \leq G$ with $\|F - F'\| \leq 5\epsilon$.
\end{lemma}

The following result concerning projections has been called ``Jordan's Lemma''. We restate the proof because it is simple and because we want to emphasize that the decomposition is orthogonal, because the results cited below that state the result in these terms do not clearly mention this property.

\begin{prop}\label{Jordan-block}
Let $P, Q$ be two projections on finite dimensional Hilbert space $\mathcal H$. They induce an orthogonal decomposition of $\mathcal H$ into one and two dimensional spaces $\mathcal H_i$ that are invariant under both $P$ and $Q$ and irreducible in the sense that $P|_{\mathcal H_i}=Q|_{\mathcal H_i}$ only if $\dim \mathcal H_i = 1$.
\end{prop}
\begin{remark}
Jordan's lemma shows that orthogonal projections $P, Q$ onto subspaces in any dimension is just the direct sum of simple cases that we already understand: one dimensional projections in at most two dimensions. With this in mind, one sees parallels between Euler's decomposition of rotation matrices in $\R^n$.
\end{remark}

\begin{proof}
Consider the reflections $R = 2P - 1, S = 2Q-1$ and the unitary operator $U = RS$. We only need to prove the result for $R$ and $S$, instead of $P$ and $Q$.

If $v$ is an eigenvector of $U$ with eigenvalue $\lambda$, then we claim that $H_1 = \operatorname{span}\{v, Rv\}$ is invariant under $R$ and $S$. Clearly, it is invariant under $R = R^{-1}$ so we check invariance under $S$.
\[Sv = R^{-1}Uv = \lambda Rv.\]
Also, $U^{-1} = S^{-1}R^{-1} = SR$ so
\[S(Rv)=U^\ast v = \overline{\lambda}v.\]
So, $H_1$ is invariant under $R$ and $S$. This is a subspace of at most two dimensions and because $R$ and $S$ are self-adjoint we obtain that $H_1$ reduces $R$ and $S$. If $P$ and $Q$ agree on $H_1$, we can break down $H_1$ into one dimensional subspaces on which $P$ and $Q$ agree. Thus we can restrict $R$ and $S$ to the orthogonal complement of $H_1$ and the result then follows by infinite descent.  
\end{proof}
A proof of this result and a discussion about this from the perspective of research in Quantum Computation can be found in Section 3.3.1 of \cite{Quantum Algorithms}. Note that sometimes, as in Section 2.1 of \cite{QMA} where the above proof is primarily taken, this result is stated in terms of unitary matrices with spectrum in $\{-1, 1\}$, which are just the reflection across the ranges of $P,Q$ given by $2P-1, 2Q-1$, as in \cite{Pironio}. Note that more ``functional analytic'' perspectives for results related to this can be found in  \cite{Jordan}, \cite{Dixmier}, and \cite{Halmos}.

Note that there is no such generalization of Jordan's lemma to more than two projections. In fact, Davis proved in \cite{three projections} that the Banach algebra of all bounded linear operators on a separable Hilbert space is generated by only three projections (and the identity).

The way that we use Jordan's lemma is the form from \cite{Hastings}:
\begin{prop}\label{Jordan}
Let $P, Q$ be two projections on Hilbert space $H$. Then there is a basis $\{p_i\}$ of the range of $P$ such that $(p_i, Qp_j)=0$ for $i \neq j$.
\end{prop}
\begin{proof}
Let ${\mathcal H} = \sum_i {\mathcal H}_i$ as in Proposition \ref{Jordan-block} just above. 
If $\mathcal H_i$ is one dimensional, then it is an eigenspace for both $P$ and $Q$, so if it is a $1$-eigenspace of $P$ let $p_i$ be a unit vector spanning $\mathcal H_i$, otherwise we do nothing.
If $\mathcal H_i$ is two dimensional, then because $\mathcal H_i$ is invariant under $P$, it has an eigenvector there. Because $P$ is not a multiple of the identity when restricted to $\mathcal H_i$, we obtain a $1$-eigenvector $p_i$ for $P$ which spans the image of $P$ restricted to $\mathcal H_i$.

We obtain that the span of the $p_i$ is the range of $P$ and because the $\mathcal H_i$ are orthogonal, the $Qp_i\in\mathcal H_i$ are as well. So, the $Qp_i$ are orthogonal and so $(p_i, Qp_j)=0$ for $i \neq j$.
\end{proof}
\begin{remark}
This result has the following geometric interpretation. 

If $P \leq Q$, then we can pick any basis $\{p_i\}$ of the range of $P$. If $Q \leq P$, then we can form $\{p_i\}$ as an orthonormal basis of the range of $Q$ and extend it to an orthonormal basis of the range of $P$.

If $P\not\leq Q$, then there is an annoying fact that it may be true that two vectors $v, w$ in the range of $P$ may be orthogonal, but $Qv$ and $Qw$ may not be. For example, let $\mathcal H = \C^3$ and let $P$ project onto the subspace spanned by the first two standard basis vectors $e_1, e_2$. If the range of $Q$ is the span of $e_3$ and $e_1 + e_2$, then $Qe_1 = Qe_2$. In other words, $e_1, e_2$ are orthogonal but by applying $Q$ we have eliminated the components of $e_1$ and $e_2$ that contribute to their orthogonality. A way to avoid this phenomenon is to pick basis vectors $v = \frac{1}{\sqrt{2}}(e_1 + e_2)$ and $w = \frac{1}{\sqrt{2}}(e_1 - e_2)$ for the range of $P$ so that $Qv = v$ and $Qw = 0$, so that $\{Qv, Qw\}$ is an orthogonal set of vectors. 

This construction is more complicated in the general case when $P$ and $Q$ do not intersect orthogonally or when there are multiple two dimensional subspaces in the decomposition. In particular, even though $\{Qp_i\}$ are orthogonal, we are not guaranteed a lower bound for the norm of these vectors. This is easily seen in the case that $P$ and $Q$ project onto arbitrary lines in $\C^2$.
\end{remark}

\section{Lemmas on Spectral Projections and Commutators}\label{Spectral Projections}
\begin{prop}\label{D-K theorem} (Davis-Kahan $\sin\theta$ Theorem)
There exists a constant $c > 0$ such that for self-adjoint $A, B \in M_n(\C)$ and $S_1, S_2 \subset \R$ we have
\[\|E_{S_1}(A)E_{S_2}(B)\| \leq \frac{c}{\operatorname{dist}(S_1, S_2)}\|A - B\|.\]
If there is a $\delta > 0, \alpha, \beta \in \R$ with  $S_1 \subset [\alpha, \beta]$ and $S_2 \subset (-\infty, \alpha - \delta] \cup [\beta + \delta, \infty)$, we have
\[\|E_{S_1}(A)E_{S_2}(B)\| \leq \frac{1}{\delta}\|A - B\|.\]
\end{prop}
\begin{proof}
For a proof see Sections 10 and 11 of \cite{Bhatia}.
\end{proof}

\begin{remark}\label{Davis-Kahan interpretation}
If $P, Q$ are projections, the quantity $\|PQ\|$ can be thought of as the ``minimal $|\cos \theta|$'' between any two lines in the range of $P$ and $Q$, respectively, because 
\[\|PQ\| = \max_{|v|, |w| = 1} |(v, PQw)| = \max_{|v|=|w|=1,  Pv = v, Qw = w} |(v, w)|.\]

With Jordan's lemma if $\dim \mathcal H_i = 2$ for all $i$, this gives $\|PQ\| = \max_i\cos\theta_i$, where $\theta_i$ are the angles between the rank one projections of $P$ and $Q$ restricted to $\mathcal H_i$.

If there are one dimensional $\mathcal H_i$, then $\|PQ\| = \max(\max_i\cos\theta_i,  \max_j\|PQ|_{\mathcal H_j}\|)$, where $i$ ranges over the two dimensional subspaces and $j$ ranges over the one dimensional subspaces. $\|PQ|_{\mathcal H_j}\|$ equals zero if $P|_{\mathcal H_j}$ and $Q|_{\mathcal H_j}$ are not identical and equals one otherwise.

In our notation, \cite{Bhatia} states that ``the name `$\sin\theta$ theorem' comes from the interpretation of $\|PQ\|$ as the sine of the angle between $\operatorname{Ran}(P)$ and $\operatorname{Ran}(Q)^\perp$.''
\end{remark}

\begin{example}
Consider $A = \begin{pmatrix} 1 & 0 \\ 0 & 0 \end{pmatrix}, A_1 = \begin{pmatrix} 1+\epsilon_1 & 0 \\ 0 & 0 \end{pmatrix}, A_2 = \begin{pmatrix} 1 & \epsilon_2 \\ 0 & 0 \end{pmatrix}$. 
Perturbing $A$ to get $A_1$ causes the eigenvalues (but not eigenvectors) to drift with $\|A - A_1\| = \epsilon_1$ and $E_{\{1\}}(A) = E_{\{1+\epsilon_1\}}(A_1)$.
For $A_2$, we get that the eigenvectors rotate but the eigenvalues remain unchanged.\\
(See \cite{Davis} for more about this behavior in general.)
\end{example}

Here are two results regarding spectral projections and commutators. One result has the commutator small and the other has that the operators have a small difference. This result is part of an argument used in \cite{P-S}.
\begin{prop}\label{comm-proj}
Suppose that $C, D\in M_n(\C)$ with $D$ self-adjoint. Then for sets $S_1 \subset (-\infty, \alpha], S_2 \subset [\alpha+\delta, \infty)$, we have
\[\|E_{S_1}(D)CE_{S_2}(D)\| \leq \frac{\|[C,D]\|}{\delta}.\]
\end{prop}
\begin{proof}
Fix a vector $v$ in the range of $E_{S_1}(D)$ and a vector $w$ in the range of $E_{S_2}(D)$. Because $v, w$ are arbitrary, we wish to show that 
\[|(v,Cw)|\leq \frac{1}{\delta}\|[C,D]\|.\]

Let the notation: $v_\lambda$ (which may be zero) represent a vector such that $D v_{\lambda} = \lambda v_{\lambda}$ and $w_\mu$ likewise. That is, $v_\lambda, w_\mu$ are eigenvectors or zero. We write the orthogonal eigenspace decompositions $v = \sum_{\lambda}v_{\lambda}, w = \sum_{\mu}w_{\mu}$, where for the rest of the proof $\lambda$ will be an element of $S_1$ and $\mu$ an element of $S_2$.

Because $\lambda, \mu \in \R$ and $D^\ast = D$, we see that
\[(\lambda - \mu)(v_{\lambda}, Cw_{\mu}) = (Dv_{\lambda},  Cw_{\mu}) - ( v_{\lambda}, CD w_{\mu}) = -(v_{\lambda}, [C,D]w_{\mu}).\]

List the eigenvectors $v_\lambda\neq 0$ as $u_1, \dots, u_r$ and the corresponding $\lambda$'s as $a_1, \dots, a_r$. Also list the eigenvectors $w_\mu \neq 0$ as $u_{s}, \dots, u_n$ and the corresponding $\mu$'s as $b_{s}, \dots, b_n$. Let $u_{r+1}, \dots, u_{s-1}$ be some unit vectors so that $\mathscr B=\{\frac{1}{|u_i|}u_i\}$ forms an orthonormal basis for $\C^n$. Define otherwise $a_i = 0, b_j = \delta$. We define a matrix $T \in M_n(\C)$ in the basis $\mathscr B$ so that if $1\leq i \leq r$ and $s \leq j \leq n$ then 
$T_{i,j} = (\frac{u_i}{|u_i|}, [C,D]\frac{u_j}{|u_j|})$ and $T_{i,j} = 0$ otherwise. This is so that our auxiliary operator $T$ satisfies $(u_i, Tu_j) = (u_i, [C,D]u_j)$ when $1\leq i \leq r$ and $s \leq j \leq n$ and $T =  E_{S_1}(D)[C,D]E_{S_2}(D)$.

Applying Lemma \ref{P-S} we obtain
\begin{align*}
|(v,  Cw)| &= \left|\sum_{\lambda,\mu} (v_{\lambda}, Cw_\mu)\right| = \left|\sum_{\lambda,\mu} \frac{1}{\lambda - \mu}(v_{\lambda}, [C,D]w_\mu)\right| = \left|\sum_{\lambda,\mu} \frac{1}{\lambda - \mu}(v_{\lambda}, Tw_\mu)\right| \\
&= \left|\sum_{i,j} \frac{1}{a_i - b_j}(u_i, Tu_j)\right| = \left|(v, \left(\frac{1}{a_i - b_j}T_{i,j}\right)w)\right|  \leq \frac{1}{\delta}\|[C,D]\|.
\end{align*}
\end{proof}

\begin{example}\label{strict comm-proj}
This result is nicely illustrated by taking $D = \begin{pmatrix} a & 0 \\ 0 & b \end{pmatrix}$ and $C = \begin{pmatrix} 0 & \epsilon \\ \epsilon & 0 \end{pmatrix}$.
Then $[C,D] = \begin{pmatrix} 0 & -(b-a)\epsilon \\ (b-a)\epsilon & 0 \end{pmatrix}$ and $\|E_{\{a\}}(D)CE_{\{b\}}(D)\| = \epsilon$.

This example shows that the result is sharp. We also see the behavior that if $C$ and $D$ almost commute and if $b-a$ is large then $\epsilon$ is small. Alternatively, if $b-a$ is small then $\epsilon$ is not required to be small (but it cannot be large). This can be interpreted as saying that $C$ approximately does not send vectors in one eigenspace of $D$ into a ``far away'' eigenspace  of $D$ but can for nearby eigenspaces.
\end{example}

We use the convention for the Fourier transform from \cite{Hastings}. With the definition
\[\hat{f}(k) = \frac{1}{2\pi}\int_\R f(x)e^{-ikx}dx,\]
we have
\[f(x) = \int_\R \hat{f}(k)e^{ikx}dk\]and
\[\widehat {f \ast g}(k) = 2\pi \hat{f}(k)\hat{g}(k),\]
The following is a simplification of Theorem 3.2.32 from \cite{B&R} that has been modified so that it agrees with this convention. 
\begin{prop}\label{B-R proposition} Let $f \in C^0(\R)\cap L^1(\R)$ with \begin{align}\label{Const.}
C_f = \int_\R |k||\hat{f}(k)|dk < \infty.
\end{align}
If $A, B\in M_n(\C)$ with $A$ self-adjoint, then
\[f(A) = \int_\R \hat{f}(k)e^{ikA}dk,\] 
\[[f(A), B]= i\int_\R  k\hat{f}(k) \int_0^1 e^{itkA}[A,B]e^{i(1-t)kA}dtdk\]
and so
\[\|[f(A), B]\| \leq C_f\|[A, B]\|.\]
\end{prop}

We apply Proposition \ref{B-R proposition} to get the following:

\begin{lemma}\label{spectral gap}
If $A, B \in M_n(\C)$ with $A$ self-adjoint and having no spectrum in $(a,b)$, then 
\[\|\,[E_{(-\infty, a]}(A), B]\,\| \leq \frac{c_2}{b-a}\|[A,B]\|.\]
Here $c_2 = 4\inf\|\hat\rho\|_{L^1(\R)}$, where the infimum is taken over all $\rho \in C^0(\R)\cap L^1(\R)$ supported in $[-1,1]$ with $\int_\R \rho = 1$. 
\end{lemma}
\begin{proof}
Write $b = a+2\epsilon_0$. We restrict to $0<\epsilon < \epsilon_0$. Let $f_\epsilon(x) = \chi_{[R-\epsilon, a + \epsilon]}\ast \rho_\epsilon$, where $\rho$ is as above, $\rho_\epsilon(x) = \frac{1}{\epsilon}\rho(x/\epsilon)$, and $R = \min \sigma(A)$. Then $f_\epsilon(x)$ equals $1$ on $[R, a]$ and equals zero outside $(-R - 2\epsilon, b)$. Because $A$ has no spectrum in $(-\infty, R) \cup (a,b)$, we see that $f_{\epsilon}(A) = E_{(-\infty, a]}(A)$.
Recalling that
\[\hat{\chi}_{[c,d]}(k) = -e^{-ik\frac{c+d}{2}}\frac{\sin(\frac{d-c}{2}k)}{\pi k},\] we obtain
\begin{align*} C_{f_\epsilon} &= 
2\pi\int_{\R}\left|k\hat{\chi}_{[R-\epsilon, a + \epsilon]}(k)\hat{\rho}(\epsilon k)\right|dk \leq \frac{2}{\epsilon}\int_{\R}\left|\hat{\rho}(\epsilon k)\right|\epsilon dk  = \frac{2}{\epsilon}\|\hat{\rho}\|_{L^1(\R)}\\
\end{align*}
and the result follows.
\end{proof}
\begin{remark}
There are other ways to pick the interpolating function $f_\epsilon$ in the proof, but ultimately we know that this result, up to the constant, is sharp and the best constant is at least $1$.  This is because 
\begin{align*}
\|E_{(-\infty,a]}&(A)B - BE_{(-\infty,a]}(A)\| \geq \|(E_{(-\infty,a]}(A)B - BE_{(-\infty,a]}(A))E_{(-\infty,a]}(A)\| \\
&= \|(1-E_{(-\infty,a]}(A))BE_{(-\infty,a]}(A)\| = \|E_{[b,\infty)}(A)BE_{(-\infty,a]}(A)\|
\end{align*}
and we know that we have equality in $\|E_{[b,\infty)}(A)BE_{(-\infty,a]}(A)\| \leq \frac{c}{b-a}\|[A,B]\|$ from Example \ref{strict comm-proj}.
\end{remark}
\begin{remark}
If we choose $\rho = \chi_{[-1/2,1/2]}\ast  \chi_{[-1/2,1/2]}$, then $\rho$ satisfies the required properties with
\[\|\hat{\rho}\|_{L^1} = 2\pi\int_{\R}\left(\frac{\sin(k/2)}{\pi k}\right)^2 dk = \frac1\pi\int_{\R}\left(\frac{\sin(k)}{ k}\right)^2 dk = 1.\]
So, we obtain that the result holds with $c_2 =  4$.
\end{remark}

\section{Relevant Lieb-Robinson Bounds}\label{Lieb-Robinson section}
The following result appeared in \cite{Davidson} in the discussion following its Lemma 3.1, while the statement and proof appearing here is modified from \cite{Hastings}. 

Note that  the results in Lemma \ref{finite range}, Lemma \ref{modified finite range}, and Corollary \ref{finite range normal} still hold if $A$ is not self-adjoint, though in that case $H$ is not necessarily self-adjoint. In our applications, $A$ will always be self-adjoint.
\begin{lemma}\label{finite range}
There exist constants $c_0, c_1 > 0$ such that given $\Delta > 0$ and self-adjoint $B \in M_n(\C)$, there exists $H\in M_n(\C)$ such that $\|H\|\leq c_1\|A\|$ with
\[\|A-H\|\leq \frac{c_0}{\Delta}\|[A,B]\|,\]
\[\|[H,B]\| \leq c_1\|[A,B]\|,\]
and
$E_{S_1}(B)HE_{S_2}(B) = 0$ for any $S_1, S_2 \subset\R$ with $\operatorname{dist}(S_1, S_2) \geq \Delta$. If $A$ is self-adjoint then $H$ can be chosen to be self-adjoint.
\end{lemma}

\begin{proof}
Let $f \in C^0(\R) \cap L^1(\R)$ be supported in $[-1,1]$ with $f(0) = 1$ such that the constants $c_0, c_1$ defined below are finite. Write 
\[H = \int_{\mathbb{R}}e^{i\frac{k}{\Delta}B}Ae^{-i\frac{k}{\Delta}B}\hat{f}(k)dk.\] 

To show that $E_{S_1}(B)HE_{S_2}(B) = 0$ if $\operatorname{dist}(S_1, S_2) \geq \Delta$, we pick $v_\lambda, v_\mu$ two eigenvectors of $B$ with $\lambda \in S_1, \mu \in S_2$. Then 
\begin{align*}
(v_{\lambda}, Hv_{\mu}) &= \int_{\R}\left(e^{-i\frac{B}{\Delta}k}v_\lambda, Ae^{-i\frac{B}{\Delta}k}v_\mu\right) \hat{f}(k) dk = (v_\lambda, Av_\mu)\int_{\R}e^{i\frac{\lambda-\mu}{\Delta}k}\hat{f}(k) dk \\
&= (v_\lambda, Av_\mu)f\left(\frac{\lambda - \mu}{\Delta}\right) = 0.
\end{align*}

A key fact (see \cite{B&R} Lemma 3.2.31) is that for $k \in \mathbb{R}$, $\|[A, e^{ikB}]\| \leq |k|\|[A,B]\|$. We then see that because $f(0) = 1$,
\begin{align*}
\|A - H\| &=
\left\|A\int_{\mathbb{R}}e^{i\frac{k}{\Delta}B}e^{-i\frac{k}{\Delta}B}\hat{f}(k)dk -\int_{\mathbb{R}}e^{i\frac{k}{\Delta}B}Ae^{-i\frac{k}{\Delta}B}\hat{f}(k)dk\right\|
\\
&= \left\|\int_{\mathbb{R}}[A,e^{i\frac{k}{\Delta}B}]e^{-i\frac{k}{\Delta}B}\hat{f}(k)dk\right\| \leq \frac{\|[A,B]\|}{\Delta}\int_{\mathbb{R}}|k\hat{f}(k)|dk = \frac{c_0}{\Delta}\|[A,B]\|,
\end{align*}
where $c_0 = \int_{\mathbb{R}}|k\hat{f}(k)|dk$.
Also,
\[\|[H,B]\| = \left\|\int_{\mathbb{R}}e^{i\frac{k}{\Delta}B}[A,B]e^{-i\frac{k}{\Delta}B}\hat{f}(k)dk\right\| \leq \|[A,B]\|\int_{\mathbb{R}}|\hat{f}(k)|dk = c_1 \|[A,B]\|,\]
where $c_1 = \int_{\mathbb{R}}|\hat f(k)|dk$.

\end{proof}
\begin{remark}
The ``best choice'' of the constants $c_0, c_1$ depends on the function $f$ that we chose. A similar remark concerning the constant in Proposition \ref{D-K theorem} can be made where the geometry of $S_1$ and $S_2$ are more general. See \cite{Bhatia} and \cite{Extremal}.

It is stated in \cite{Hastings} that $c_1$ can be chosen to equal $1$ with the provided function $f(x)=(1-x^2)^3\chi_{[-1,1]}(x)$, but this is not so because\footnote{As seen by a simple mathematical software calculation.} $\hat{f}(0) > 0, \hat{f}(10) < 0$ so $\|\hat{f}\|_{L^1} > f(0) = 1$.
In \cite{Davidson}, Davidson uses a very similar but less direct proof for this result (saying that it is a modification of Theorem 4.1 in \cite{Bhatia Perturbation}) and by citing some literature obtains $c_0 = 8, c_1 = 4$.
\end{remark}
We also have the following which generalizes the above result to localize $A$ with respect to commuting self-adjoint matrices $B_1, \dots, B_m$.
\begin{lemma}\label{modified finite range}
Let $c_0, c_1$ be the constants from Lemma \ref{finite range}. For $\Delta > 0$, $A\in M_n(\C)$, and commuting self-adjoint $B_1, \dots, B_m \in M_n(\C)$, there exists $H\in M_n(\C)$ such that $\|H\|\leq c_1\|A\|$ with
\[\|A-H\|\leq \frac{c_0c_1^{m-1}}{\Delta}\sum_{j=1}^m\|[A,B_j]\|,\]
\[\|[H ,B_j]\| \leq c_1^m\|[A,B_j]\|,\]
and
$E_{S_1}(B_j)HE_{S_2}(B_j) = 0$ for any $S_1, S_2 \subset\R$ with $\operatorname{dist}(S_1, S_2) \geq \Delta$. If $A$ is self-adjoint, $H$ can be chosen to be self-adjoint.
\end{lemma}
\begin{proof}
We essentially iterate the above construction because the $B_j$ commute. Let $f$, $c_0$, and $c_1$ be as in the proof of Lemma \ref{finite range} and set
\begin{align*}
H &= \int_{\mathbb{R}^m}e^{i\sum_{j=1}^m\frac{k_j}{\Delta}B_j}Ae^{-i\sum_{j=1}^m\frac{k_j}{\Delta}B_j}\hat{f}(k_1)\cdots\hat{f}(k_m)dk_1\dots dk_m.
\end{align*}
Following the calculation in the argument in the proof of Lemma \ref{finite range}, that the $B_j$ commute, that $[A,-]$ is a derivation, and $f(0) = 1$,
\begin{align*}
\|A-H\| &= \left\|\int_{\mathbb{R}^m}\left[A,\prod_{j=1}^me^{i\frac{k_j}{\Delta}B_j}\right]e^{-i\sum_{j=1}^m\frac{k_j}{\Delta}B_j}\hat{f}(k_1)\cdots\hat{f}(k_m)dk_1\dots dk_m\right\|\\
&\leq \sum_{j=1}^m\|[A,B_j]\|\int_{\mathbb{R}^m}\frac{|k_j|}{\Delta}|\hat{f}(k_1)|\cdots|\hat{f}(k_m)|dk_1\dots dk_m\\
&= \frac{c_0c_1^{m-1}}{\Delta}\sum_{j=1}^m\|[A,B_j]\|.
\end{align*}
Also, for any index $j_0$ between $1$ and $m$,
\begin{align*}
\|[H,B_{j_0}]\| &= \left\|\int_{\mathbb{R}^m}e^{i\sum_{j=1}^m\frac{k_j}{\Delta}B_j}[A,B_{j_0}]e^{-i\sum_{j=1}^m\frac{k_j}{\Delta}B_j}\hat{f}(k_1)\cdots\hat{f}(k_m)dk_1\dots dk_m\right\| \\
&\leq c_1^m\|[A,B_{j_0}]\|.
\end{align*}
Let $j'$ index the integers in $[1,m]\setminus \{j_0\}$. Because the $B_j$ commute, we have that
\[H = \int_{\R}e^{i\frac{k_{j_0}}{\Delta}B_{j_0}}\left(\int_{\mathbb{R}^{m-1}}e^{i\sum_{j'}\frac{k_{j'}}{\Delta}B_{j'}}Ae^{-i\sum_{j'}\frac{k_{j'}}{\Delta}B_{j'}}\prod_{j'}\hat{f}(k_{j'})dk_{j'}\right)e^{-i\frac{k_{j_0}}{\Delta}B_{j_0}}\hat{f}(k_{j_0})dk_{j_0}.\] 
So, $H$ has finite range of at least $\Delta$ with respect to $B_{j_0}$ by the argument in the proof of Lemma \ref{finite range}.
\end{proof}
For $N$ normal, we can write it as the sum of the commuting $\frac{N+N^\ast}{2},i\frac{N-N^\ast}{2i}$ and so we have the following consequence. 
\begin{corollary}\label{finite range normal}
Let $c_0, c_1$ be the constants from Lemma \ref{finite range}. For $\Delta > 0$,  $N$ normal and $A$ in $M_n(\C)$, there exists $H \in M_n(\C)$ such that $\|H\|\leq c_1\|A\|$ with
\[\|A-H\|\leq \frac{2c_0c_1}{\Delta}\|[A,N]\|,\]
\[\|[H,N]\|\leq 2c_1^2\|[A,N]\|,\]
and $E_{S_1}(N)HE_{S_2}(N)=0$ for any $S_1, S_2 \subset \C$ with $\operatorname{dist}(S_1,S_2)\geq \sqrt{2}\Delta$. If $A$ is self-adjoint, $H$ can be chosen to be self-adjoint.
\end{corollary}
\begin{proof}
Apply Lemma \ref{modified finite range} with $B_1=\Re N = \frac{N+N^\ast}{2}, B_2=\Im N = \frac{N-N^\ast}{2i}$ using
\[\|[A,B_j]\|\leq \|[A,N]\|,\]
\[\|[H,N]\|\leq \|[H,\Re N]\|+\|[H,\Im N]\|.\]
These two inequalities give the two inequalities of the lemma.

To obtain the third result, note that the distance between points $\lambda \in S_1, \mu \in S_2$ is bounded above by
\[\sqrt{2}\max(| \lambda_1- \mu_1|,| \lambda_2- \mu_2|),\]
where $\Re \lambda=\lambda_1, \Im \lambda = \lambda_2, \Re\mu = \mu_1, \Im\mu=\mu_2$.
So, if $v$ is a $\lambda$-eigenvector for $N$ and $w$ is a $\mu$-eigenvector for $N$ with $\lambda \in S_1, \mu \in S_2$ then $v$ is a $\lambda_j$-eigenvector for $B_j$ and $w$ is a $\mu_j$-eigenvector for $B_j$ and
\[\sqrt{2}\Delta \leq |\lambda-\mu|\leq \sqrt{2}\max(| \lambda_1- \mu_1|,| \lambda_2- \mu_2|).\]
Using the last result of Lemma \ref{modified finite range} we obtain $(v,Hw)=0$.

\end{proof}

The following Lieb-Robinson type result is from \cite{Hastings} where its statement and proof originate. Note that we are now requiring that $H$ be a contraction.
\begin{thm} \label{Lieb-Robinson}
Let $H, B$ self-adjoint be such that $\|H\| \leq 1$ and $E_{S_1}(B)HE_{S_2}(B) = 0$ for any $S_1, S_2 \subset\R$ with $\operatorname{dist}(S_1, S_2) \geq \Delta$. Let $v_{LR} = e^2\Delta$. Then  for $|t| \leq \operatorname{dist}(S_1,S_2)/v_{LR}$,
\begin{align}
\|E_{S_1}(B)e^{itH}E_{S_2}(B)\| \leq e^{-\operatorname{dist}(S_1,S_2)/\Delta}.
\end{align}
\end{thm}
\begin{proof}
By iteration we get that the range of $H^nE_{S_2}(B)$ lies in the range of spectral projection of $B$ on the open $n\Delta$ neighborhood of $S_2$. So, $E_{S_1}(B)H^nE_{S_2}(B) = 0$ if $\operatorname{dist}(S_1, S_2) \geq n\Delta$. Then expressing $e^{itH}$ in $E_{S_1}(B)e^{itH}E_{S_2}(B)$ as the standard exponential power series, we get for $m = \lceil \frac{\operatorname{dist}(S_1, S_2)}{\Delta} \rceil$
\[\|E_{S_1}(B)e^{itH}E_{S_2}(B)\| \leq \sum_{n \geq m}\frac{|t|^n}{n!} \leq \frac{1}{e}\sum_{n \geq m}\left(\frac{e|t|}{n}\right)^n \leq \frac{1}{e}\left(\frac{(e|t|/m)^m}{1-e|t|/m}\right),\]
where the second inequality follows from the following reductions. $n! \geq e(n/e)^n$ follows from the inequality \[\frac{\log(n!)}{n} \geq \log\left(\frac{n(n+1)}{2n}\right) \geq \log\left(\frac{n}{e^{1-1/n}}\right).\] The first equality is a convexity inequality and the second inequality follows by removing $\log$'s giving
$n(1/2-1/e^{1-1/n}) +1/2 \geq 0$. A computation verifies the cases $n = 1, 2, 3$ and when $n > 1/(1-\log(2))\approx 3.26$,   we have $1/2-1/e^{1-1/n} > 0$.

Then because $e|t|/m \leq \frac{\operatorname{dist}(S_1,S_2)}{em\Delta}\leq e^{-1}$ we have \[\|E_{S_1}(B)e^{itH}E_{S_2}(B)\| \leq e^{-1}e^{-m}/(1-e^{-1}) = e^{-m}/(e-1) \leq e^{-\frac{\operatorname{dist}(S_1, S_2)}{\Delta} }.\]
\end{proof}

\begin{remark}
Note that this does not actually use the fact that $H$ is self-adjoint or that $t$ is real, but the power series expression of $f(x) = e^{itx}$ applied to the matrix $H$ with norm at most $1$. Likewise, one might expect there to be similar estimates for other analytic functions $f$. See the reference in Remark \ref{Benzi_Remark}.

This is in line with the interpretation of Theorem \ref{Lieb-Robinson} in terms of the coefficients of a matrix as follows. Let $\beta = (v_1, \dots, v_n)$ be some basis of $\C^n$ and $Bv_j = jv_j$, be a ``position'' operator which scales each basis vector by its index. Then the condition $E_{S_1}(B)HE_{S_2}(B) = 0$ for $\operatorname{dist}(S_1,S_2) \geq \Delta$ tells us that $H$ is $2\lceil\Delta\rceil$-banded. Then, as mentioned in Remark \ref{Benzi_Remark}, one expects exponential decay of the entries away from the diagonal for analytic $f$ (with some conditions). This is similar to what we have in the next result, which is key to our application of the Lieb-Robinson result. More generally, we will not look at $m$-banded matrices, but block tridiagonal matrices.
\end{remark}

\begin{corollary}\label{operator Lieb-Robinson}
Let $H, B$ self-adjoint be such that $\|H\| \leq 1$ and $E_{S_1}(B)HE_{S_2}(B) = 0$ for any $S_1, S_2 \subset\R$ with $\operatorname{dist}(S_1, S_2) \geq \Delta$. Then for $f \in C^0(\R) \cap L^1(\R)$, 
\begin{align}\label{superpolynomial decay}
\|E_{S_1}(B)f(H)E_{S_2}(B)\| \leq \int_{|k| > \frac{\operatorname{dist}(S_1, S_2)}{e^2\Delta}}|\hat{f}(k)|dk +  \|\hat{f}\|_{L^1(\R)}e^{-\operatorname{dist}(S_1,S_2)/\Delta}.
\end{align}
\end{corollary}
\begin{proof}
By the representation result in Proposition \ref{B-R proposition} and by Theorem \ref{Lieb-Robinson} we see that
\begin{align*}
\|E_{S_1}(B)f(H)E_{S_2}(B)\| &\leq  \int_\R \|E_{S_1}(B)e^{ikH}E_{S_2}(B)\| |\hat{f}(k)|dk \\
&\leq \int_{|k| > \frac{\operatorname{dist}(S_1, S_2)}{e^2\Delta}}|\hat{f}(k)|dk + \int_{|k| \leq \frac{\operatorname{dist}(S_1, S_2)}{e^2\Delta}} e^{-\frac{\operatorname{dist}(S_1,S_2)}{\Delta}}|\hat{f}(k)|dk \\
&\leq \int_{|k| > \frac{\operatorname{dist}(S_1, S_2)}{e^2\Delta}}|\hat{f}(k)|dk +  \|\hat{f}\|_{L^1(\R)}e^{-\frac{\operatorname{dist}(S_1,S_2)}{\Delta}}.
\end{align*}
\end{proof}
\begin{remark}
Following \cite{Hastings}, we will use the smoothness of the function $f$, so that the tail estimate for $\hat{f}$, $\Phi(t) = \int_{|k| \geq t}|\hat{f}(k)|dk$ decreases faster than any polynomial and $\|\hat{f}\|_{L^1}$ is at most a constant, to obtain fast decay of (\ref{superpolynomial decay}). Note that if we do not care for ``faster than any polynomial'', then we could just assume some smoothness for $f$. See the end of Section \ref{smooth partitions} for some estimates.

\vspace{0.1in}

The above two statements are more-or-less explicit in \cite{Hastings}. The following are implicitly used and, along with the former results, go under the umbrella of     ``Lieb-Robinson bounds'':
\begin{thm}
Let $H, B$ self-adjoint be such that $\|H\| \leq 1$ and $E_{S'}(B)HE_{S''}(B) = 0$ for any $S'' \subset S' \subset \R$ with $\operatorname{dist}(S'', \R\setminus S') \geq \Delta$. Let $v_{LR} = e^2\Delta$. Then for $|t| \leq \operatorname{dist}(S'', \R\setminus S')/v_{LR}$,
\begin{align}
\|\left[e^{itH} - e^{itE_{S'}(B)HE_{S'}(B)}\right]E_{S''}(B)\| \leq 3e^{-\operatorname{dist}(S'', \R\setminus S')/\Delta}.
\end{align}
\end{thm}
\begin{proof}The proof proceeds essentially as that of the ``original'' Lieb-Robinson result, Theorem \ref{Lieb-Robinson}, taking $S_1= \R \setminus S'$ and $S_2=S''$. For $n > 0$, let $S_n''$ be the open $n\Delta$ neighborhood of $S''$ and $S_0''=S''$.

Just as in the proof of Theorem \ref{Lieb-Robinson}, $H$ maps the range of $E_{S_n''}(B)$ into the range of $E_{S_{n+1}''}(B)$ and hence $H^n$ maps the range of $E_{S''}(B)$ into the range of $E_{S_n''}(B)$.
We prove by induction that 
\[[E_{S'}(B)HE_{S'}(B)]^{n}E_{S''}(B)=E_{S'}(B)H^nE_{S''}(B)\] for $n \leq m = \lceil \frac{\operatorname{dist}(S_1, S_2)}{\Delta} \rceil$ by noting it is clearly true when $n=0,1$ and that when $n \leq m-1$, $S_{n}'' \subset S'$ and hence
\[[E_{S'}(B)HE_{S'}(B)]^{n+1}E_{S''}(B)=[E_{S'}(B)HE_{S'}(B)][E_{S'}(B)H^nE_{S'}(B)]E_{S''}(B).\] Because $H^n$ maps the range of $E_{S''}(B)$ into the range of $E_{S_{n}''}(B)\leq E_{S'}(B)$ so
\begin{align*}[E_{S'}(B)&HE_{S'}(B)][E_{S'}(B)H^nE_{S'}(B)]E_{S''}(B)=[E_{S'}(B)HE_{S'}(B)]E_{S_{n}''}(B)H^nE_{S''}(B)\\
&=E_{S'}(B)HE_{S_{n}''}(B)H^nE_{S''}(B)=E_{S'}(B)H^{n+1}E_{S''}(B).\end{align*} This gives us the desired result.

Expressing the exponentials in $E_{S'}(B)\left[e^{itH} - e^{itE_{S'}(B)HE_{S'}(B)}\right]E_{S''}(B)$ as the standard power series, we obtain
\begin{align*}
\|E_{S'}(B)&\left[e^{itH} - e^{itE_{S'}(B)HE_{S'}(B)}\right]E_{S''}(B)\|\\
&=\left\|E_{S'}(B)\left[\sum_{n=0}^\infty\frac{(it)^nH^n}{n!} -\sum_{n=0}^\infty\frac{(it)^n\left(E_{S'}(B)HE_{S'}(B)\right)^n}{n!}\right]E_{S''}(B)\right\|\\
&\leq 2\sum_{n \geq m}\frac{|t|^n}{n!} \leq 2e^{-\operatorname{dist}(S'', \R\setminus S')/\Delta},
\end{align*}
where the last inequality follows by the argument in the proof of Theorem \ref{Lieb-Robinson}. Also, by Theorem \ref{Lieb-Robinson}, 
\[\|E_{S'}(B)e^{itH}E_{S''}(B) - e^{itH}E_{S''}(B)\| = \|E_{\R \setminus S'}(B)e^{itH}E_{S''}(B)\| \leq e^{-\operatorname{dist}(S'', \R\setminus S')/\Delta}.\]
The result then follows because, expanding as a power series, we see that 
\[E_{S'}(B)e^{itE_{S'}(B)HE_{S'}(B)}E_{S''}(B) =  e^{itE_{S'}(B)HE_{S'}(B)}E_{S''}(B).\]
\end{proof}
\begin{corollary}\label{operator Lieb-Robinson modified}
Let $H, B \in M_n(\C)$ be self-adjoint such that $\|H\| \leq 1$ and $E_{S'}(B)HE_{S''}(B) = 0$ for any $S'' \subset S' \subset \R$ with $\operatorname{dist}(S'', \R\setminus S') \geq \Delta$. Then for $f \in C^0(\R) \cap L^1(\R)$, let $H' = E_{S'}(B)HE_{S'}(B)$ so
\begin{align}\label{superpolynomial decay}\nonumber
\|[f(H) - f(H')]E_{S''}(B)\| \leq 2\int_{|t| > \frac{\operatorname{dist}(S'', \R\setminus S')}{e^2\Delta}}|\hat{f}(k)|dk +  3\|\hat{f}\|_{L^1(\R)}e^{-\operatorname{dist}(S'', \R\setminus S')/\Delta}.
\end{align}
\end{corollary}
\begin{proof}
By the representation result in Proposition \ref{B-R proposition} and using the previous theorem, we have
\begin{align*}
\|[f(H) - f(E_{S'}(B)HE_{S'}(B))]&E_{S''}(B)\| \\
&\leq \int_{\R}\|[e^{itH} - e^{itE_{S'}(B)HE_{S'}(B)}]E_{S''}(B)\||\hat{f}(k)|dk  \\
&\leq 2\int_{|t| > \frac{\operatorname{dist}(S'', \R\setminus S')}{e^2\Delta}}|\hat{f}(k)|dk +  3\|\hat{f}\|_{L^1(\R)}e^{-\frac{\operatorname{dist}(S_1,S_2)}{\Delta}}.
\end{align*}
\end{proof}

How we will use this is is to form nonconsecutively orthogonal subspaces $\mathcal Y_i$ where neighboring subspaces have significant overlap. Then for a vector $v$ in the span of these spaces, we can break it up into an orthogonal sum $\sum v_i$ where $v_i \in \mathcal Y_i'$ is well nested inside $\mathcal Y_i$. Then we can guarantee that $f(H)v_i$ is approximately equal to $f(H_i)v_i$, where $H_i$ is $H$ restricted to $\mathcal Y_i$. See the properties of the spaces $\mathcal N_i$ in Section \ref{N_i}.
\end{remark}

\section{Davidson's Reformulations of Lin's Theorem}\label{reformulations}

Davidson's three equivalent reformulations of Lin's theorem are:

($Q$): For every $\epsilon > 0$ there is a $\delta > 0$ such that if $A, B$ are self-adjoint contractions with $\|[A,B]\| \leq \delta$, then there are commuting contractions $A' , B'$ such that $\|A-A'\|, \|B-B'\| \leq \epsilon$.

\vspace{0.05in}

($Q'$): For every $\epsilon > 0$, there is an $L_0$ with the following property. If $J$ is a self-adjoint contraction that is  block tridiagonal with respect to the following orthogonal subspaces $\mathcal V_1, \dots, \mathcal V_{L_0}$, then there is a subspace $\mathcal W$ such that $\mathcal V_1 \subset \mathcal W \subset \bigoplus_{i=1}^{L_0-1} \mathcal V_i$ and $\|[J,P_{\mathcal W}]\| \leq \epsilon.$

\vspace{0.05in}

($Q''$): For every $\epsilon > 0$, there is an $L_0$ with the following property. If $\mathcal K$ is a finite dimensional subspace of $L^2([0,1])$ and $M_x$ is the multiplication operator by $x$ on $L^2([0,1])$, then there is a subspace $\mathcal W$ such that $\mathcal K \subset \mathcal W \subset \operatorname{span}_{0\leq i\leq L_0-1}M_x^i(\mathcal K)$ and $\|[J,P_{\mathcal W}]\| \leq \epsilon.$
 
\vspace{0.1in}

We now discuss the reduction of ($Q$) to ($Q'$). Although Section 3 of \cite{Davidson} does much of what we discuss in this section, we follow the notation and argument of \cite{Hastings}, where Davidson's reformulations are not explicitly mentioned. We start with $A, B$ self-adjoint with $\|A\|, \|B\| \leq 1$ and $\|[A,B]\| \leq \delta$. 
We will pick $\Delta = \delta^{\gamma_0}<<\delta$, but leave it as is for now (for simplicity and also for intuition).

We now proceed into Section III of \cite{Hastings}, where we construct ``the new basis''.
Because $\sigma(A) \subset I = [-1,1]$, we will cut up the interval $I$ into $n_{cut}$ (chosen later to equal $\lceil1/\Delta^{\gamma_1}\rceil$ for $0<\gamma_1<1$) many disjoint intervals
$I_i$ of the form $I_{i} = [-1 + i\frac{2}{n_{cut}} , -1 + (i+1)\frac{2}{n_{cut}})$ for $0 \leq i \leq n_{cut}-2$ and $I_{n_{cut}-1} = [1-\frac{2}{n_{cut}} , 1]$. 
Then we ``pinch'' $H$ by the projections\footnote{If $A$ is a matrix and we have orthogonal projections $P_1, \dots, P_k$ such that $\sum_i P_i = I$ then $\sum_i P_iAP_i$ is the ``pinching'' of $A$ by $\{P_i\}$ by the terminology in \cite{Davis}.} $E_{I_i}(B)$ getting matrices $J_i$ acting on a spaces ${\mathcal B}_{i+1} = \operatorname{Ran}(E_{I_i}(B))$, $i=0,\dots,n_{cut}-1$.

Pictorially, we now focus only on $I_i$ and on that it has length $\frac{2}{n_{cut}}$. 
We will pick $n_{cut}$ later so that $\Delta = o(1/n_{cut})$ and hence we can partition $I_i$ into at least $\lfloor\frac{2/n_{cut}}{\Delta}-1\rfloor =: L$ many intervals $I_i^j$ of length at least $\Delta$ and at most $2\Delta$ (only the first and last subinterval may have length greater than $\Delta$). If $H$ has finite range $\Delta$, we obtain that $J_i=E_{I_i}(B)HE_{I_i}(B)$ is block tridiagonal with respect to the subspaces $\mathcal V_i^j$ that the $E_{I_i^j}(B)$ project onto.
 
Naturally, given our choices above, we will get $L$ increasing like a negative power of $\Delta$, so either $J_i$ has many blocks or, because at least one block is empty because $\mathcal V_i^j = 0$, we will obtain a nontrivial reducing subspace for $J_i$ such that the following lemma trivially holds.

Now, we state the main lemma (Lemma 2) for the argument in \cite{Hastings}.

\begin{lemma}\label{main lemma}
Let $J$ be self-adjoint with $\|J\| \leq 1$ acting on ${\mathcal B}$ with $L$ orthogonal subspaces ${\mathcal V}_i$ with respect to which $J$ is block tridiagonal. Then there is a subspace ${\mathcal W}$ of  ${\mathcal B}$ satisfying \begin{enumerate}
\item For any $v \in {\mathcal V}_1, |P_{{\mathcal W}^\perp}(v)| \leq \epsilon_3|v|$.
\item For any $w \in {\mathcal W}, |P_{{\mathcal W}^\perp}(Jw)| \leq \epsilon_4|w|$.
\item For any $w \in {\mathcal W}, |P_{{\mathcal V}_L}(w)| \leq \epsilon_5|w|$,
\end{enumerate}
where for $i=3,4$, $\epsilon_i =\frac1{L^{\gamma_i}}E_i(L)$ where $E_i(t)$ grows slower than any (positive) power (of $t$) and $\epsilon_5$ decays faster than any power of $L$. We prove $\gamma_i \geq 1/4$.
\end{lemma}
The ``construction'' in the proof only works if $L$ is large, which is masked by the undefined nature of the $E_i$ in the lemma above, because we can imagine defining $E_i(1/L)$ to be large for all $L$ small. 
How large $L$ needs to be is undetermined because of the non-constructive step in the proof discussed in Section \ref{application of Lin's theorem}.

\begin{remark} Intuitively, Item 1 above means that ${\mathcal V}_1$ is almost contained in ${\mathcal W}$, Item 2 means that $J$ is almost invariant under ${\mathcal W}$, and Item 3 means that ${\mathcal W}$ is almost orthogonal to ${\mathcal V}_L$.
\end{remark}

\begin{remark}\label{alternate}
Note that Item 3 above is formulated differently in \cite{Hastings}, but they are equivalent because they both say that for all $w \in {\mathcal W}, v \in {\mathcal V}_L, |(w,v)|\leq \epsilon_5|w||v|$. The form stated above is what is proved in \cite{Hastings}, whose proof we follow.

Similarly we have a dual statement for Lemma \ref{main lemma}: 
\begin{enumerate}\label{duality}
\item For any $w \in \mathcal W^\perp, |P_{\mathcal V_1}(w)| \leq \epsilon_3|w|$.
\item For any $w \in \mathcal W^\perp, |P_{\mathcal W}(Jw)| \leq \epsilon_4|w|$.
\item For any $v \in \mathcal V_L, |P_{\mathcal W}(v)| \leq \epsilon_5|v|.$
\end{enumerate}
where the second statement uses that $J$ is self-adjoint.
\end{remark}

Although the following result is what we will use in the proof of Lin's theorem (because it simplifies the discussion of the $\mathcal W_i$ later in the section), Lemma \ref{main lemma} is called the main lemma, because most of this paper is dedicated to proving it. This formulation more closely follows Davidson's and Szarek's treatments.
\begin{lemma}\label{modified lemma}
Let $J$ be self-adjoint with $\|J\| \leq 1$ acting on ${\mathcal B}$ with $L$ orthogonal subspaces ${\mathcal V}_i$ with respect to which $J$ is block tridiagonal. Then there is a subspace ${\mathcal W}$ of  ${\mathcal B}$ satisfying  $\mathcal V_1 \leq \mathcal W \perp \mathcal V_L$ and
\[\|P_{\mathcal W}^{\perp}JP_{\mathcal W}\| \leq \epsilon_4 + 10\max(\epsilon_3,\epsilon_5) =: \epsilon_2,\]
where $\epsilon_3,\epsilon_4,\epsilon_5$ are as in Lemma \ref{main lemma}.
\end{lemma}
\begin{proof}
This is a direct application of Lemma \ref{projection lemma} to Lemma \ref{main lemma} which gives us a projection $P'=P_{\mathcal W}$ such that $\| P_{\mathcal V_1}{P'}^\perp\| \leq \epsilon_3$, $\|{P'}^\perp J P'\| \leq \epsilon_4,$ and $\|P'P_{\mathcal V_L}\| \leq \epsilon_5$.

Then we get a projection $P$ such that $P_{\mathcal V_1} \leq P \perp P_{\mathcal V_L}$ and $\|P - P'\| \leq 5\max(\epsilon_3,\epsilon_5)$. Then we have
\begin{align*}
\|(1-P)JP\| &\leq \|(1-P')JP'\| + \|(1-P')JP'-(1-P)JP\| \leq \epsilon_4 + 10\max(\epsilon_3,\epsilon_5).
\end{align*}
\end{proof}
\begin{remark}\label{Linexpbootstrap}
The proof of Theorem 3.2 of \cite{Davidson} (the equivalence of this lemma and Lin's theorem) shows that we can pick $\epsilon_2 = 1/L^{1/2}$ by simply applying Lin's theorem from \cite{KS}. As we discuss in Section \ref{estimates}, obtaining $\epsilon_2 = 1/L$ provides the optimal exponent for Lin's theorem.

Hastings showed that for $H$ a tridiagonal matrix then one has $\epsilon_2 = E(L)/L$. It appears that it is not known whether this holds in general. 
\end{remark}

Now, using Lemma \ref{modified lemma} we obtain Lin's theorem as detailed in \cite{Hastings}. Consider ``the new basis'' of subspaces for $0 \leq j \leq n_{cut},$
\[\tilde {\mathcal B}_j: {\mathcal W}_1, {\mathcal W}_1^\perp \oplus  {\mathcal W}_2, \dots,  {\mathcal W}_i^\perp \oplus  {\mathcal W}_{i+1}, \dots, {\mathcal W}_{n_{cut}-1}^\perp \oplus  {\mathcal W}_{n_{cut}}, {\mathcal W}_{n_{cut}}^\perp.\]
See Figure 1 below.
\begin{figure}[htp]
    \centering
    \includegraphics[width=12cm]{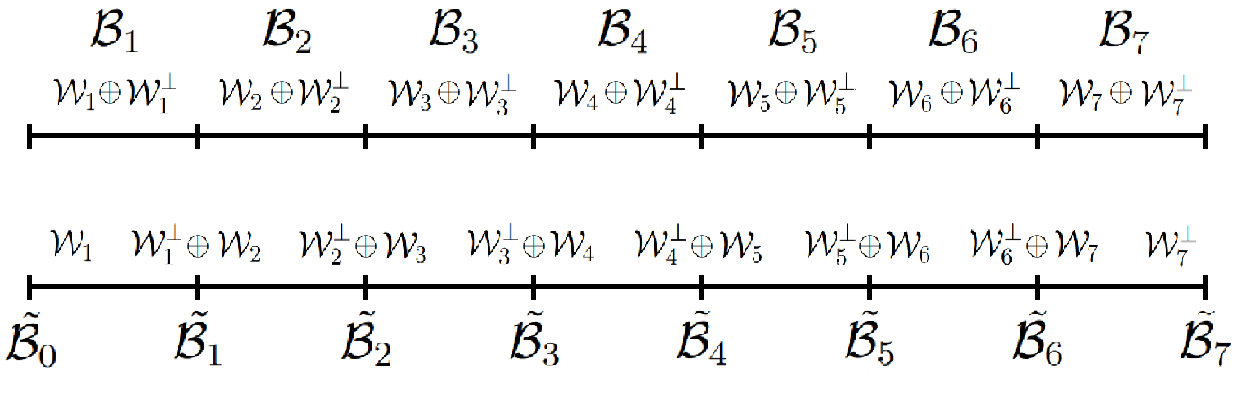}
    \caption{The formation of the new basis from the old. Note that the interval illustrates the spectrum of $B$. Note $n_{cut}=7$.}
\end{figure}
Note the abuse of notation that we will use for the rest of the \\section: ${\mathcal W}_i^\perp =  {\mathcal B}_i \ominus {\mathcal W}_i$ and ${\mathcal B}_i^\perp = {\mathcal B} \ominus {\mathcal B}_i$.

For simplicity, set ${\mathcal W}_0={\mathcal W}_0^\perp = {\mathcal B}_0 = 0$ and ${\mathcal W}_{n_{cut}+1}={\mathcal W}_{n_{cut}+1}^\perp = {\mathcal B}_{n_{cut}+1} = 0$ and $I_{n_{cut}} = \{1\}$ so that for $0 \leq i \leq n_{cut}$, $\tilde {\mathcal B}_i = {\mathcal W}_i^\perp \oplus  {\mathcal W}_{i+1} \subset {\mathcal B}_i \oplus {\mathcal B}_{i+1}$.
Now, let $B'$ be the block identity operator which equals the right endpoint of $I_i$ multiplied by the identity on $\tilde {\mathcal B}_i \subset {\mathcal B}_i \oplus {\mathcal B}_{i+1}$, $0 \leq i \leq n_{cut}-1$. Because $B$ has eigenvalues in $I_{I}$ on ${\mathcal B}_{I+1}$, we see that $\|B-B'\| \leq \frac{2}{n_{cut}}$.

We now show that the $\tilde {\mathcal B}_i$ are almost invariant subspaces for $H$, where we let $H'$ be the pinching of $H$ along $\tilde {\mathcal B}_i$ so then $\|H-H'\|\leq 2\epsilon_2$.
This amounts to showing that applying $H$ to any vector in ${\mathcal W}_i$ remains in ${\mathcal W}_{i-1}^\perp\oplus {\mathcal W}_i$ except that very small amounts are  permitted to ``leak out''. A similar statement would hold for ${\mathcal W}_{i-1}^\perp$. The idea is that applying $H$ to ${\mathcal W}_{i} \subset {\mathcal B}_{i}$ can only make ${\mathcal W}_{i}$ leak out into ${\mathcal V}_1^{i+1}$ and ${\mathcal V}_{L}^{i-1}$ through ${\mathcal V}_{L}^i$ and ${\mathcal V}_1^{i}$, respectively. In details, let ${\mathcal B}_i$ be written as the orthogonal direct sum of ${\mathcal V}_1^i, \dots, {\mathcal V}_L^i$ as in the statement of Lemma 2.

Recall that by Item 2, $\|P_{{\mathcal W}_i^\perp}HP_{{\mathcal W}_i}\| \leq \epsilon_2.$ 
By the block tridiagonality of $H$, $P_{{\mathcal B}_i^\perp}HP_{{\mathcal W}_i} = P_{{\mathcal B}_i^\perp} H(P_{{\mathcal V}_1^i} + P_{{\mathcal V}_L^i})P_{{\mathcal W}_i}$.
Now, $P_{{\mathcal V}_L^i}P_{{\mathcal W}_i} =0$ and $P_{{\mathcal B}_i^\perp} HP_{{\mathcal V}_1^i}P_{{\mathcal W}_i} = P_{{\mathcal B}_i^\perp} HP_{{\mathcal V}_1^i}$ maps into  ${\mathcal V}_L^{i-1} \subset \mathcal W_{i-1}^\perp$ so $P_{{\mathcal B}_i^\perp}HP_{{\mathcal W}_i} = P_{{\mathcal V}_L^{i-1}} HP_{{\mathcal W}_i} = P_{{\mathcal W}_{i-1}^\perp} HP_{{\mathcal W}_i}$. Hence,
\begin{align*}
\|(1-P_{\tilde {\mathcal B}_{i-1}})&HP_{\tilde {\mathcal B}_{i-1}}P_{\mathcal W_i}\|=\|HP_{{\mathcal W}_i} - P_{{\mathcal W}_i}HP_{{\mathcal W}_i} - P_{{\mathcal W}_{i-1}^\perp} HP_{{\mathcal W}_i}\|\\
&\leq \|P_{{\mathcal B}_i}HP_{{\mathcal W}_i} - P_{{\mathcal W}_i}HP_{{\mathcal W}_i}\| + \|P_{{\mathcal B}_i^\perp}HP_{{\mathcal W}_i} - P_{{\mathcal W}_{i-1}^\perp} HP_{{\mathcal W}_i}\| \leq \epsilon_2 + 0.
\end{align*}

Likewise, by taking adjoints, $\|P_{\mathcal W_i}HP_{{\mathcal W}_i^\perp}\| \leq \epsilon_2$. 
By the block triangularity of $H$, $P_{{\mathcal B}_i^\perp}HP_{{\mathcal W}_i^\perp} = P_{{\mathcal B}_i^\perp}H(P_{{\mathcal V}_1^i} + P_{{\mathcal V}_L^i})P_{{\mathcal W}_i^\perp}$. 
Now, $P_{{\mathcal V}_1^i}P_{{\mathcal W}_i^\perp} = 0$ and $P_{{\mathcal B}_i^\perp}H P_{{\mathcal V}_L^i}P_{{\mathcal W}_i^\perp} = P_{{\mathcal B}_i^\perp}H P_{{\mathcal V}_L^i}$ 
maps into $\mathcal V_1^{i+1}\subset \mathcal W_{i+1}$, so $P_{{\mathcal B}_i^\perp}HP_{{\mathcal W}_i^\perp} = P_{{\mathcal V}_1^{i+1}}H P_{{\mathcal W}_i^\perp} = P_{{\mathcal W}_{i+1}}H P_{{\mathcal W}_i^\perp}$. Hence,
\begin{align*}
\|(1-P_{\tilde {\mathcal B}_{i}})&HP_{\tilde {\mathcal B}_{i}}P_{{\mathcal W}_i^\perp}\|=\|HP_{{\mathcal W}_i^\perp} - P_{{\mathcal W}_{i+1}}HP_{{\mathcal W}_i^\perp} - P_{{\mathcal W}_{i}^\perp}HP_{{\mathcal W}_i^\perp}\| \\
&\leq \|P_{{\mathcal B}_{i}^\perp}HP_{{\mathcal W}_i^\perp} - P_{{\mathcal W}_{i+1}}HP_{{\mathcal W}_i^\perp}\| + \|P_{{\mathcal B}_{i}}HP_{{\mathcal W}_i^\perp} - P_{{\mathcal W}_{i}^\perp}HP_{{\mathcal W}_i^\perp}\| \leq 0 + \epsilon_2.
\end{align*}
So, it follows that
\[\|P_{\tilde {\mathcal B}_j}^\perp HP_{\tilde {\mathcal B}_j}\| \leq 2\epsilon_2.\] 
Set $H' = \sum_j P_{\tilde {\mathcal B}_j} HP_{\tilde {\mathcal B}_j}$. 
Because the spaces $\tilde {\mathcal B}_i$ are orthogonal, we see that $\|H - H'\| \leq 2\epsilon_2$. By construction $B'$ and $H'$ commute and we conclude.

\section{Estimates, Putting it all Together}\label{estimates}
In this section we explore how the various constants discussed previously can be chosen to get the best possible decay in Lin's theorem. There are cases where one might not want the best possible estimates, because perhaps picking a different rate allows one to not use Lemma \ref{main lemma} but instead use a different construction such as that of \cite{Szarek}. 

We now summarize the steps taken thus far in the construction, assuming a solution of Lemma \ref{main lemma} with exponent $\gamma_2 = \min(\gamma_3, \gamma_4)$.
From $A, B$ self-adjoint with $\|A\|, \|B\| \leq 1$, $\|[A,B]\| \leq \delta$ we get (by applying Lemma \ref{finite range}) $H$ such that $\|H\| \leq 1, \|[H,B]\| \leq Const.\delta$, and $\|A-H\| \leq Const. \delta/\Delta$. We get $[H',B'] = 0$ with $\|B - B'\| \leq 2/n_{cut} \sim 2\Delta^{\gamma_1}$ for $0<\gamma_1<1$ and $\|H-H'\|\leq 2\epsilon_2 \leq 2L^{-\gamma_2 }E(\frac{1}{L})$, where $L \sim \frac{2}{n_{cut}\Delta} \sim  2\Delta^{\gamma_1-1}$.
So, writing $\Delta = \delta^{\gamma_0}$, $0<  \gamma_0 < 1$, we get $2L^{-1}\sim \Delta^{1-\gamma_1} = \delta^{\gamma_0(1-\gamma_1)}$ hence
\begin{align}\label{AH-est}
\|A - H'\| \leq Const. \delta/\Delta + 2\epsilon_2 \leq  Const.\delta^{1-\gamma_0} + \delta^{\gamma_0(1-\gamma_1)\gamma_2}E(\delta^{-1})
\end{align}
and
\begin{align}\label{B-est}
\|B - B'\| \leq 2/n_{cut} \leq  Const.\delta^{\gamma_0\gamma_1}.
\end{align}

We want to minimize these, noting that $\gamma_2$, which comes from Lemma 2, is the only constant that we cannot choose. Setting the three exponents equal we get $\gamma_0 = \frac{1}{1+\gamma_1}$ and $\gamma_1 = \frac{\gamma_2}{1+\gamma_2}$, so the common value of the exponents will be chosen to be $\gamma = \frac{\gamma_2}{1+2\gamma_2}$.
Because the value of $\gamma_2$ that we obtain is $\gamma_2 = 1/4$, we obtain $\gamma = 1/6$.

Note that this method cannot give the optimal exponent of $\gamma = \frac12$ (only potentially $\frac12+\epsilon$ if we are allowed to take $\gamma_2$ arbitrarily large). This is partly because of the averaging that we do to make $H$ ``finite range'' as expressed in the title of Section II.A of \cite{Hastings}, as a consequence of the $\delta/\Delta$ factor. 

If we were able to remove that factor (by improving the result or starting in a special case of $A$) then we would then have to compare $\gamma_0(1-\gamma_1)\gamma_2$ with $\gamma_0\gamma_1$ which gives $\gamma_1 = \frac{\gamma_2}{1+\gamma_2}$. Because we need to have $A$ finite range of distance $\Delta$ with respect to $B$, we see that $\|[A,B]\| \leq Const. \Delta$. This gives $\delta \leq Const.\delta^{\gamma_0}$ so $0 < \gamma_0 \leq 1$. Because we do not need to worry about bounding $A - H$ (because we are starting with $A = H$) we get
\begin{align}
\|A - H'\| \leq \delta^{\gamma_0(1-\gamma_1)\gamma_2}E_0(1/\delta)
\end{align}
and
\begin{align}
\|B - B'\| \leq 2\delta^{\gamma_0\gamma_1}.
\end{align}

So, we see that the rate we get is $\gamma = \frac{\gamma_2}{1+\gamma_2}$, which is an improvement. The value of $\gamma_2$ is $1$ to then get the optimal exponent $\gamma = 1/2$. 

Note that Davidson in the proof of the equivalence of ($Q$) and ($Q'$) showed that Lin's theorem implies that we can pick $\gamma_2 = 1/2$, which is an estimate that gives Lin's theorem with $\gamma = 1/4$.

\section{General Approach to the construction of $\mathcal W$.}\label{General Approach}
In this section we define subspaces $\tilde{\mathcal X}_i$ and $\tilde{\mathcal X}_i'$ to describe the motivating ideas used in \cite{Szarek}. We then define smooth cut off functions $\mathcal F^{r,w}_{\omega_0}(t)$ and use them to define the spaces $\mathcal X_i$, which will be put together in some sense to define $\mathcal W$ following \cite{Hastings}.

Because we need to essentially recover $\mathcal V_1$ as a subspace of $\mathcal W$ by Item 1, we might consider finding subspaces $\tilde{\mathcal X}_i$ such that $\mathcal V_1$ is approximately a subset of the sum of the $\tilde{\mathcal X_i}$. For example, consider breaking up $[-1,1]$ into $n_{win}$ many disjoint intervals $I_i$ similarly to how we did before then we could consider $\tilde {\mathcal X}_i = [\chi_{I_i}(J)](\mathcal V_1)$.
This space recovers $\mathcal V_1$ because any element $v$ in $\mathcal V_1$ can be written $v = \sum_i x_i$, where $x_i = \chi_{I_i}(J) v$.

These spaces are almost invariant: Let $\omega(i)$ be the midpoint of $I_i$ and $x \in \tilde {\mathcal X}_i$. We have that 
\begin{align}\label{243}
|(J - \omega(i))x| = |(J - \omega(i))\chi_{I_i}(J)x| \leq \|(J - \omega(i))\chi_{I_i}(J)\||x| \leq \frac{|I_i|}{2}|x|.
\end{align}
This shows that $Jx$ is almost in $\operatorname{span}(x) \subset \tilde {\mathcal X}_i$, where the error depends on the norm of $x$.
Also, these spaces $\tilde{\mathcal X}_i$ are also orthogonal, so we obtain that if $v = \sum_i x_i$ then $|v|^2 = \sum_i |x_i|^2$.

Putting all of this together, if we set $\tilde {\mathcal W} = \bigoplus_i \tilde{\mathcal X}_i$, we have Item 1 because it contains $\mathcal V_1$. It satisfies Item 2 because if $w \in \tilde {\mathcal W}$ then we can write it as $w = \sum_i x_i$, $x_i \in \tilde {\mathcal X}_i$. 
Because the $I_i$ are disjoint,
\begin{align}\label{264} \nonumber|Jw - \sum_i\omega(i)x_i|^2 &= |\sum_{i}(J - \omega(i))x_i|^2 = \sum_{i}|(J - \omega(i))x_i|^2 \\
&\leq \sum_i \left(\frac{|I_i|}{2}\right)^2|x_i|^2 = \left(\frac{|I_i|}{2}\right)^2|w|^2.
\end{align}
So because $x_i \in \tilde {\mathcal W}$,
\begin{align}\label{267}
|P_{\tilde{\mathcal W}^\perp}Jw| = |P_{\tilde{\mathcal W}^\perp}\left(Jw - \sum_i\omega(i)x_i\right)| \leq \frac{|I_i|}{2}|w|.
\end{align}
This general set-up is common ground for Szarek's and Hastings' arguments. They differ in how to modify this core argument to get a result concerning Item 3. We proceed with a modification used in both constructions. 

We want to show that Item 3 holds as well for $\tilde{\mathcal W}$, but it might not unless we address an issue.  Recall that $\tilde {\mathcal X}_i = \{\chi_{I_i}(J) v: v \in \mathcal V_1\}$. Let $S_1$ be a matrix whose columns are a fixed orthonormal basis of $\mathcal V_1$, where $d = \dim\mathcal V_1$. Then $S_1$ is an isometric isomorphism between $\mathbb{C}^d$ and $\mathcal V_1$ and $R(\chi_{I_i}(J)S_1)=\mathcal {\tilde X}_i$. 

We might hope that we have exponential decay of the columns of $S_1$ (which would provide what we want) along the lines of Remark \ref{Benzi_Remark}. This might not happen. To see this, \cite{Hastings} gives an example of an $n \times n$ matrix of the form:
\[\begin{pmatrix} 0 & 1/4 & 0 &  & \\
1/4 & 0 & 1/4 & 0 &  \\
0 & 1/4 & 0 & 1/4 & 0 \\
 & 0 & 1/4 & 0 & 1/4 & \ddots\\
 &  & 0 & 1/4 & \ddots & \ddots & 0 \\
 &  &  & \ddots & \ddots & 0 & 1/4 \\
  &  & & & 0 & 1/4 & 1/2\\
\end{pmatrix}\]
Here, we can take $\mathcal V_1$ to be the first basis vector. For $n$ large, its spectrum is distributed finely in $[-1/2, 1/2]$ with a single (with multiplicity one) eigenvalue at around $5/8$.

For example, MATLAB calculations give that if $S_1$ is just the first standard basis vector then $\chi_{[5/8-1/100, 5/8+1/100]}(J)S_1$ is, for $n = 10$,
\[10^{-3}\times 
    (0.0016,
   0.0040,
   0.0084,
   0.0171,
   0.0343,
   0.0686,
   0.1373,
   0.2747,
   0.5493,
1.0987)^T,\]
and for $n = 50$
\begin{align*}10^{-15} \times
(&0.0000, \dots, 0.0000, 0.0001, 0.0002, 0.0005,  0.001,  0.002,   0.004, 0.008, 0.016,\\
& 0.031, 0.062, 0.125,
   0.251, 0.502, 1.004)^T.
\end{align*}

The pattern is that the norm of this vector is very small, but does not have small projection onto $\mathcal V_L$ relative to its norm. One way of addressing this phenomenon is to disallow vectors that have too small norm, by removing columns of $\chi_{I_i}(J)S$ that are too small.
However, that does not exclude the case that all of its columns are large, but a linear combination of them are of this problematic type. For example, we might have chosen a different basis for $\mathcal V_1$. 
If one modifies the operators $\chi_{I_i}(J)S$ so that we still maintain Item 1 and Item 2 and removing the issue discussed above, we might be closer to showing Item 3. 

The approach to this problem is to write $\tilde \tau_i = \chi_{I_i}(J)S_1$ and remove its singular values that are too small by letting $\tilde Z_i$ project onto the eigenspace of ${\tilde \tau_i}^\ast \tilde \tau_i$ for eigenvalues less than some $\lambda_{min}$. 
Then define new spaces $\tilde {\mathcal X}_i'$ to be the range of $\tilde\tau_i(1-\tilde Z_i)$ and let $\tilde {\mathcal W}' = \bigoplus_i \tilde {\mathcal X}_i'$, an orthogonal direct sum.

Now, this definition of $\tilde {\mathcal X}_i'$ gives that for any $x \in \tilde {\mathcal X}_i'$, we have some $\textbf{x} \in \C^d$ such that $\tilde\tau_i(1-\tilde Z_i)\textbf{x} = x$ and we can choose $\textbf{x}$ to be in the range of $1-\tilde Z_i$ so that $\tilde\tau_i \textbf{x} = x$. We call $\textbf{x}$ \emph{the} representative of $x$ (by $\tilde\tau_i$). Writing $\textbf{x} = \sum_a \textbf{x}^a$ as the orthogonal sum of eigenvectors $\textbf{x}^a$ of ${\tilde \tau_i}^\ast \tilde \tau_i$ with eigenvalues $\lambda_a \geq \lambda_{min}$ we obtain
\[|\tilde\tau_i \textbf{x}|^2 = \sum_{a,b}(\tilde \tau_i\textbf{x}^a, \tilde {\tau_i}\textbf{x}^b) = \sum_{a}(\tilde {\tau_i}^\ast \tilde \tau_i\textbf{x}^a, \textbf{x}^a) \geq \sum_{a}\lambda_{min}|\textbf{x}^a|^2,\]
so
\[|\tilde\tau_i \textbf{x}| \geq \lambda_{min}^{1/2}|\textbf{x}|.\]

This guarantees that we do not have any exponentially small ``problematic'' vectors if we make $\lambda_{min}$ go to zero only like a power of $L$. This adjustment seemingly does not offer a way to prove Item 3, however, if it is now true.

Item 1 holds for $\tilde {\mathcal W}'$ because for any $v \in \mathcal V_1$, we can write it as the orthogonal direct sum $\sum_i x_i$, where $x_i = \chi_{I_i}(J)v = \tilde\tau_i \textbf{x}$ and $v = S\textbf{x}$. Then we get that \[|\tilde\tau_i \textbf{x} - \tilde \tau_i(1-\tilde Z_i)\textbf{x}|^2 = |\tilde\tau_i\tilde Z_i\textbf{x}|^2 \leq \lambda_{min}|\textbf{x}|^2 = \lambda_{min}|v|^2.\]
Because the ranges of the $\tilde\tau_i$ are orthogonal, we get that
\begin{align*}|P_{\tilde {\mathcal W}'^\perp}v| &\leq |\sum_i\tilde \tau_i\textbf{x} - \sum_i\tilde \tau_i(1-\tilde Z_i)\textbf{x}| = 
|\sum_i\tilde \tau_i\tilde Z_i\textbf{x}| = \left(\sum_{i=0}^{n_{win}-1}|\tilde \tau_i\tilde Z_i\textbf{x}|^2\right)^{1/2}
\\
&\leq (n_{win}\lambda_{min})^{1/2}|v|.\end{align*}

We see that Item 2 holds as follows. Let $\omega(i)$ be the midpoint of the interval $I_i$.
For $w = \sum_i x_i$ with $x_i \in \tilde{\mathcal X}_i'$ so $\chi_{I_i}(J)x_i = x_i$, we have \[|(J - \omega(i))x_i|  \leq \frac{|I_i|}{2}|x_i|.\] Recall that $x_i \in \tilde{\mathcal X}'_i \subset \tilde{\mathcal X}_i$ and $Jx_i \in\tilde{\mathcal X}_i$. So, because the $\tilde{\mathcal X}_i$ are orthogonal we obtain
\begin{align*}
|P_{\tilde W'^\perp}Jw|^2 &= |P_{\tilde {\mathcal W}'^\perp}\sum_i(J-\omega(i))x_i|^2 \leq |\sum_i(J-\omega(i))x_i|^2 = \sum_i|(J-\omega(i))x_i|^2 \\
&\leq \sum_i \left(\frac{|I_i|}{2}|x_i|\right)^2=\left(\frac{|I_i|}{2}|w|\right)^2 
\end{align*} so $|P_{\tilde{\mathcal W}'^\perp} Jw| \leq \frac{|I_i|}{2}|w|$. 

An important property that we used in approximating $Jw$ was that $Jw = \sum_iJx_i$ and not only is $\sum x_i = w$ an orthogonal decomposition but the $(J-\omega(i))x_i$ are also orthogonal because the $\tau_i$ are. Later when we define $\mathcal X_i$ to be a subspace of the range of ${\mathcal F}^{0,\kappa}_{\omega(i)}(J)$ (where ${\mathcal F}^{0,\kappa}_{\omega(i)}$ are smooth overlapping cut-off functions), we see that $J(\mathcal X_{i})$ are nonconsecutively orthogonal. One could see that the above argument would work, except that we seemingly do not have \begin{align}\label{recover}
\sum_i|x_i|^2 \leq Const. |w|^2.\end{align} 

These are some of the issues that Hastings deals with using the ${\mathcal F}^{0,\kappa}_{\omega(i)}$. The subspaces that we obtain will no longer be orthogonal (which causes its own issues), but in Lemma \ref{estimates for W}.1 we obtain the estimate
\[\sum_i|x_i|^2 \leq Const.l_b|w|^2\]
where $l_b$ grows like a power of $L$ and along with the decay given in Corollary \ref{operator Lieb-Robinson} we can obtain the third property. The smoothness of ${\mathcal F}^{0,\kappa}_{\omega(i)}$ ensures the sufficient decay  of its Fourier transform and hence that of Corollary \ref{operator Lieb-Robinson}.  As stated in \cite{Hastings}, ``The smoothness will be essential to ensure that the vectors ... have most of their
amplitude in the first blocks rather than the last blocks''.

\section{Dimensional-Dependent Results}\label{Dim section}

In this section we describe the dimensional dependent results prior to \cite{Hastings} which indicate the necessity of both matrices being self-adjoint and why it is necessary that the proof of Lin's theorem involve the structure of both almost commuting matrices. This section can be skipped if reading for the proof of Hastings' result.

One comment is that Lin's theorem is essentially the $m$ dimensional case if we assume that $A$ always has at most $m$ distinct eigenvalues. 
More precisely, Pearcy and Shields' argument in \cite{P-S} indicates that if $A$ is self-adjoint and has at most $m$ distinct eigenvalues, then there are commuting matrices $A', B'$ where $A'$ is self-adjoint, $B'$ is self-adjoint if $B$ is, and
\[\|A-A'\|,\|B-B'\|\leq \left(\frac{m-1}{2}\|[A,B]\|\right)^{1/2}.\]
This allows us to choose $\delta = \frac{2}{m-1}\epsilon^2$. The ideas in Pearcy and Shields' argument are illustrated in the following cheap result.  
\begin{prop}\label{Cheap}
Let $A, B \in M_n(\C)$, where $A$ is self-adjoint with at most $m$ distinct eigenvalues. Then there are commuting $A', B'$ such that $A'$ is self-adjoint and
\[\|A - A'\|, \|B-B'\|\leq \frac{m}{\sqrt{2}}\|[A,B]\|^{1/2}.\]
\end{prop}
\begin{proof}
We partition the spectrum of $A$ into intervals $I_i$ of length at most $c_m\delta^{1/2}$, where $\delta:= \|[A,B]\|$. We choose $c_m$ later so as to make our bound for $\|A-A'\|$ equal to that of $\|B-B'\|$. 
We write $A$ as block multiples of the identity $(a_{i}I_i)$, with the $a_i$ ordered increasingly, and $B$ as a block matrix $(B_{ij})$, so then $[A,B]_{rs} = (a_r-a_s)B_{rs}$.  

We cannot avoid that some off-diagonal terms of $(B_{ij})$ might not be small where the eigenvalues of $A$  are close because of behavior depicted in Example \ref{strict comm-proj}. A way of dealing with this is to use the trick in Hastings' proof of the tridiagonal case and merge together intervals that are too close. This will be doable, because our estimate will depend on the number of eigenvalues. 

Doing this, take our partition of half-open intervals $I_i$ of length $c_m\delta^{1/2}$ and discard the intervals that do not intersect the spectrum of $A$ to obtain half-open intervals $I_j$. Merge together neighboring intervals to obtain merged intervals $\tilde I_k$. Because there were at most $m$ eigenvalues, these merged intervals have length at most $mc_m\delta^{1/2}$ and there is a gap of at least $c_m\delta^{1/2}$ (the length of one $I_i$). We let $A'$ be the block matrix formed by merging the diagonal entries of $A$ that correspond to the merged intervals and replacing the entries with the midpoint of the merged interval. Then $\|A - A'\| \leq \frac{1}{2}mc_m\delta^{1/2}$. We let $B'$ be the matrix gotten by discarding the off-diagonal blocks $B_{rs}$ where $a_r$ and $a_s$ are from different merged blocks. Because $|a_r - a_s|\|B_{rs}\| \leq \|[A,B]\| = \delta$ and for $a_r, a_s$ from different merged blocks, $|a_r - a_s| \geq c_m\delta^{1/2}$ so we obtain for these $r,s$, $\|B_{rs}\| \leq \frac{1}{c_m}\delta^{1/2}$. 

Because matrices in $M_n(\C)$ are being expressed as square block matrices with at most $m$ rows, $\|B-B'\| \leq \frac{m}{c_m}\delta^{1/2}$. Setting $c_m = \sqrt{2}$, we obtain 
\[\|A-A'\|, \|B-B'\| \leq \frac{m}{\sqrt2}\delta^{1/2}.\] 
\end{proof}
\begin{remark}\label{P-S-D}
Davidson's counter-example in \cite{Davidson} of $(n^2+1)\times (n^2+1)$ self-adjoint $A_n$ and nearly normal $B_n$ that satisfy $\|[A_n,B_n]\| = 1/n^2$ with bounded norms without being nearby commuting matrices $A_n', B_n'$ where $A_n'$ is self-adjoint shows that \cite{P-S}'s result is asymptotically optimal.
Choi in \cite{Choi} also has a similar counter-example showing that there are no nearby commuting matrices at all with the same asymptotic estimates.
\end{remark}

The choice of the groupings of eigenvalues of $A$ in Pearcy and Shileds' result does not involve $B$ in any way, besides the use of $\|[A,B]\|$. A difficulty in solving Lin's theorem involves the fact that we need to use the structure of $B$ when constructing $A'$ (and vice-versa) unless $\|[A,B]\|$ is very small. We explore these two ideas in the next two results.
\begin{prop}
Suppose that $A_k, B_k \in M_{n_k}(\C)$ are self-adjoint contractions. If\\ $\|[A_k,B_k]\|\to 0$, by Lin's theorem there are self-adjoint, commuting $A_k', B_k'$ so that\\ $\|A_k-A_k'\|, \|B_k-B_k'\|\to 0$. 
Let $m_k$ be the number of  eigenvalues of $A_k'$. 

If the choice of $A_k'$ can be made independent of $B_k$, then $\sup_k m_k < \infty$. 
\end{prop}
\begin{proof}
We prove this by contradiction. Choose a subsequence and relabel so that $\lim_k m_k = \infty$. 

By the Pigeonhole principle, there are two distinct eigenvalues $\lambda_k^1, \lambda_k^2$ of $A_k'$ with eigenvectors $v_k^1, v_k^2$, respectively, so that $|\lambda_k^1-\lambda_k^2|\leq 2/m_k$. If we define $B_{k}$ to be the linear operator satisfying $B_{k}v_k^1 =v_k^2, B_{k}v_k^2 =v_k^1$, and is identically zero on the orthogonal complement of $\operatorname{span}(v_k^1, v_k^2)$. Then $\|[A_k',B_k]\| \leq 2/m_k$ so
\[\|[A_k, B_k]\| \leq \|A_k-A_k'\|\|B_k\|+\|[A_k',B_k]\| \to 0.\]

If $B_k'$ is any matrix that commutes with $A_k'$ then the eigenspaces of $A_k'$ are invariant under $B_k'$. This means that when writing $B_k'$ as a block diagonal matrix with respect to the eigenspaces of $A_k'$, it must be diagonal. This shows that $\|B_k - B_k'\|\geq 1$, which contradicts our assumption.
\end{proof}

The following result illustrates what seems to be a rather strict requirement on the size of the commutators because it destroys the invariants for three almost commuting Hermitians. Compare to Remark \ref{P-S-D}.
\begin{corollary}\label{3HermsCommute}
Let $A_k, B_k, C_k \in M_{n_k}(\C)$ be self-adjoint contractions with $A_k$ having $m_k$ (distinct) eigenvalues. If $\|[A_k,B_k]\|, \|[A_k, C_k]\| = o(1/m_k)$ and $\|[B_k,C_k]\| = o(1)$ then there are self-adjoint commuting contractions $A'_k, B'_k, C'_k$ such that \[\|A_k-A'_k\|, \|B_k-B'_k\|, \|C_k-C'_k\| =  o(1).\] 
\end{corollary}
\begin{proof}
Note that $o(1)$ only depends on $k$. Without loss of generality, suppose that $\|[A_k, C_k]\| \leq \|[A_k, B_k]\|$. Because the construction of $A_k'$ in \cite{P-S} only uses $B_k$ through a bound for $\|[A_k, B_k]\|$ and the construction of $B_k'$ is as below, we see that it works equally as well for $C_k$.

Applying \cite{P-S}'s construction to $A_k, B_k$, we get disjoint intervals $I_1^k, \dots, I_{r^k}^k$ that cover $\sigma(A_k)$. For the projections $P_i^k = E_{I_i^k}(A_k)$, we can modify the spectrum of $A_k^i := P_i^kA_kP_i^k$ to get multiples of identity matrices $({A'}_k)^i$ which form the diagonal entries of the block diagonal matrix $A'_k$. This construction satisfies the following properties: \[\|A_k - A_k'\| \leq  \left(\frac{m_k-1}2\|[A_k,B_k]\|\right)^{1/2} =: \epsilon_k =  o(1).\]
and for $B_k'' = \sum_i P_i^k B_k P_i^k,$ $C_k'' = \sum_i P_i^k C_k P_i^k$ we have $\|B_k - B_k''\|, \|C_k - C_k''\| \leq \epsilon_k = o(1)$.
Then \[\|[B_k'', C_k'']\| \leq \|[B_k,C_k]\| + 2\|B_k-B_k''\| + 2\|C_k-C_k''\| = o(1).\] By definition, \[\|[P_i^k B_k P_i^k, P_i^k C_k P_i^k]\| = \|P_i^k [B_k'',C_k'']P_i^k \|=o(1)\] and so we can apply Lin's theorem to these block matrices, so there are commuting self-adjoint contractions $B_k^i, C_k^i$ of the same block size as $P_i^k B_k P_i^k, P_i^k C_k P_i^k$ such that $\|B_k^i - P_i^k B_k P_i^k\|, \|C_k^i - P_i^k C_k P_i^k\| = o(1)$. 
Since $P_i^k A_k' P_i^k$ is a multiple of the identity on $R(P_i^k)$, we see that forming $B_k' = \bigoplus_i B_k^i, C_k' = \bigoplus_i C_k^i$ gives us commuting self-adjoint matrices $A_k', B_k', C_k'$ satisfying the statement of the lemma.
\end{proof}

Before continuing with a discussion of Hastings' approach to Lemma \ref{modified lemma}, we pause to discuss Szarek's approach to this lemma. 
Szarek proves a dimensional-dependent version of this lemma in the following result from \cite{Szarek}. Note that this result is constructive and although is expressed in terms of Davidson's third formulation of Lin's theorem, it is readily transformed into the form that we have discussed in this paper.
\begin{prop}\label{Szarek proposition}
Let $J$ be self-adjoint with $\|J\| \leq 1$ acting on ${\mathcal B}$ with $L$ orthogonal subspaces ${\mathcal V}_i$ with respect to which $J$ is block tridiagonal. If for some $i < L$ we set $m = \dim  \mathcal V_i$, then there is a subspace $\mathcal W$ such that $\mathcal V_1 \subset \mathcal W \perp \mathcal V_L$ and
\[\|P_{{\mathcal W}^\perp}JP_{\mathcal W}\| \leq \epsilon_1(L),\]
where $\displaystyle\epsilon_1(L) = Const.\left(\frac{m}{L}\right)^{1/9}$.
\end{prop}
\begin{remark}
Szarek provides an argument that reduces the proof of the previous proposition to the case that all the $\mathcal V_i$, in particular, $\mathcal V_1$ has dimension $m$. This argument actually can be extended to say that not only can all the blocks of $J$ be large, but we only need one off-diagonal block to have small rank and so can set
$m = \rank P_{\mathcal V_{i+1}}JP_{\mathcal V_i}$.

The point of this remark is that not does this result apply for a matrix $J$ that contains blocks like these:
\[\begin{pmatrix} 
\ast & 1 &  &  \\
1 & \ast & \ast & \ast \\
 & \ast & \ast & \ast \\
 & \ast & \ast & \ast
\end{pmatrix}\]
because there is a diagonal block (corresponding to $\dim \mathcal V_i$) of size one, but this result also applies if there are blocks like these (potentially under a change of basis for the subspaces $\mathcal V_i$): \[\begin{pmatrix} 
\ast & \ast & \ast & \ast &   &   &   \\
\ast & \ast & \ast & \ast &   &   &   \\
\ast & \ast & \ast & \ast &   &   &   \\
\ast & \ast & \ast & \ast & 1 &   &   \\
  &   &   & 1 & \ast & \ast & \ast \\
  &   &   &   & \ast & \ast & \ast \\
  &   &   &   & \ast & \ast & \ast
\end{pmatrix}.\]
\end{remark}

The argument in \cite{Szarek} involves using the general approach discussed in Section \ref{General Approach}, however some modifications are made to obtain the last inequality. The construction involves several ideas. One of these ideas is that to obtain the approximate orthogonality of $\mathcal W$ and $\mathcal V_L$, Szarek forms $\mathcal W$ exactly as $\tilde{\mathcal W}'$, however the intervals $I_i$ will not cover $\sigma(J)$ but will be constructed by the application of the following lemma.
\begin{lemma}\label{measure}
Let $\mu$ be a finite positive measure on $[0,1]$ and let $\kappa > 8\eta > 0$. Then there are disjoint intervals $I_i$ for $i = 1, \dots, n_{win}$ such that
\begin{enumerate}
\item $n_{win} \leq \displaystyle\frac{2}{\kappa}$
\item $|I_j| \leq \kappa$
\item $\operatorname{dist}(I_i, I_j) \geq \eta, \, i \neq j$
\item $\mu\left([0,1]\setminus\bigcup_{j} I_j\right) \leq \displaystyle\frac{4\eta}{\kappa}\mu([0,1])$
\end{enumerate}
\end{lemma}
Based on what we did in the previous section, what remains is to show is that $\mathcal V_1$ still belongs approximately to $\mathcal W$ and that $\mathcal W$ is approximately orthogonal to $\mathcal V_L$. The first statement follows from the construction of the intervals through the application of the last inequality in the interval lemma above. Because the construction of $\mu$ is related to projection onto $\mathcal V_1$, the dimension $m$ of $\mathcal V_1$ appears in the estimates as it is the trace of this projection.

The approximate orthogonality follows from the key observation that by using polynomial interpolation, there are polynomials $p_i$ such that $0\leq p_i \leq 1$, $p_i\approx \chi_{I_i}$ on $I_i$ and $p_i \approx 0$ on the other intervals $I_j$. How well an approximation obtained depends on the degree of the $p_i$ and the separation $\eta$ from the above lemma.  

We require 
$p_i(J)$ have degree at most $L-1$. So, using a trick similar to that of Equation (\ref{267}), one obtains the estimate 
\[|P_{\mathcal V_L^\perp}w| = |P_{\mathcal V_L^\perp}\left(w- \sum_i p_i(J)v_i\right)| \leq |w- \sum_i p_i(J)v_i|\] for $w \in \mathcal W$ and any $v_i \in \mathcal V_1$. This gives the desired result for the approximate orthogonality.

Having summarized the methods used by Szarek, we state the result gotten:
\begin{thm}\label{Szarek theorem}
If self-adjoint contractions $A, B \in M_n(\C)$ have rank at most $m$, then there are self-adjoint commuting $A', B' \in M_n(\C)$ such that \[\|A-A'\|, \|B-B'\| \leq Const.(m^{1/2}\|[A,B]\|)^{2/13}.\]
\end{thm}
As Szarek mentions, a key consequence of this was improving Pearcy and Shield's result so that one gets a factor of $m^{1/2}\|[A,B]\|$ instead of $m\|[A,B]\|$ by using the assumption that both $A$ and $B$ are  self-adjoint.

\section{Smooth Partitions and Hastings' Approach}\label{smooth partitions}

The method of proof of Lemma 2 in \cite{Hastings} in getting Item 3 is to utilize smooth step functions that overlap along the lines of Section \ref{General Approach}. So, we define the general profile function that we use. Let ${\mathcal F}$ be some smooth strictly decreasing function on $[0,1]$ with all derivatives at $0$ and $1$ equal to zero, $\overline{\mathcal F}(0)= 1, \overline{\mathcal F}(1) = 0$, and $\overline{\mathcal F}(x) + \overline{\mathcal F}(1-x) = 1$ for all $x \in [0,1]$. Define ${\mathcal F}_{\omega_0}^{r,w}(\omega)$ for $\omega_0 \in \R, w > 0, r \geq 0$ by translation
\[{\mathcal F }_{\omega_0}^{r,w}(\omega) =  {\mathcal F }_{0}^{r,w}(\omega-\omega_0)\]
where ${\mathcal F }_{0}^{r,w}$ is even, identically $1$ on $[0,r]$, and equal to $\overline{\mathcal F}((x-r)/w)$ on $[r,r+w]$. So ${\mathcal F }_{0}^{r,w}$ is identically equal to $1$ on an interval centered at $\omega_0$ of radius $r$ and smoothly decreases to zero within the width $w$.
See Figure 2 below. The symmetry property of $\overline{\mathcal F}$ implies that if the $\omega(i)$ are spaced by $\kappa$ and $I_i = (\omega(i) - \kappa, \omega(i)+\kappa)$ are nonconsecutively disjoint intervals whose union contains $[-1,1]$ then $\{{\mathcal F}_{\omega(i)}^{0,\kappa}\}$ forms a partition of unity covering $[-1,1]$. 

\begin{figure}[htp]
    \centering
    \includegraphics[width=10cm]{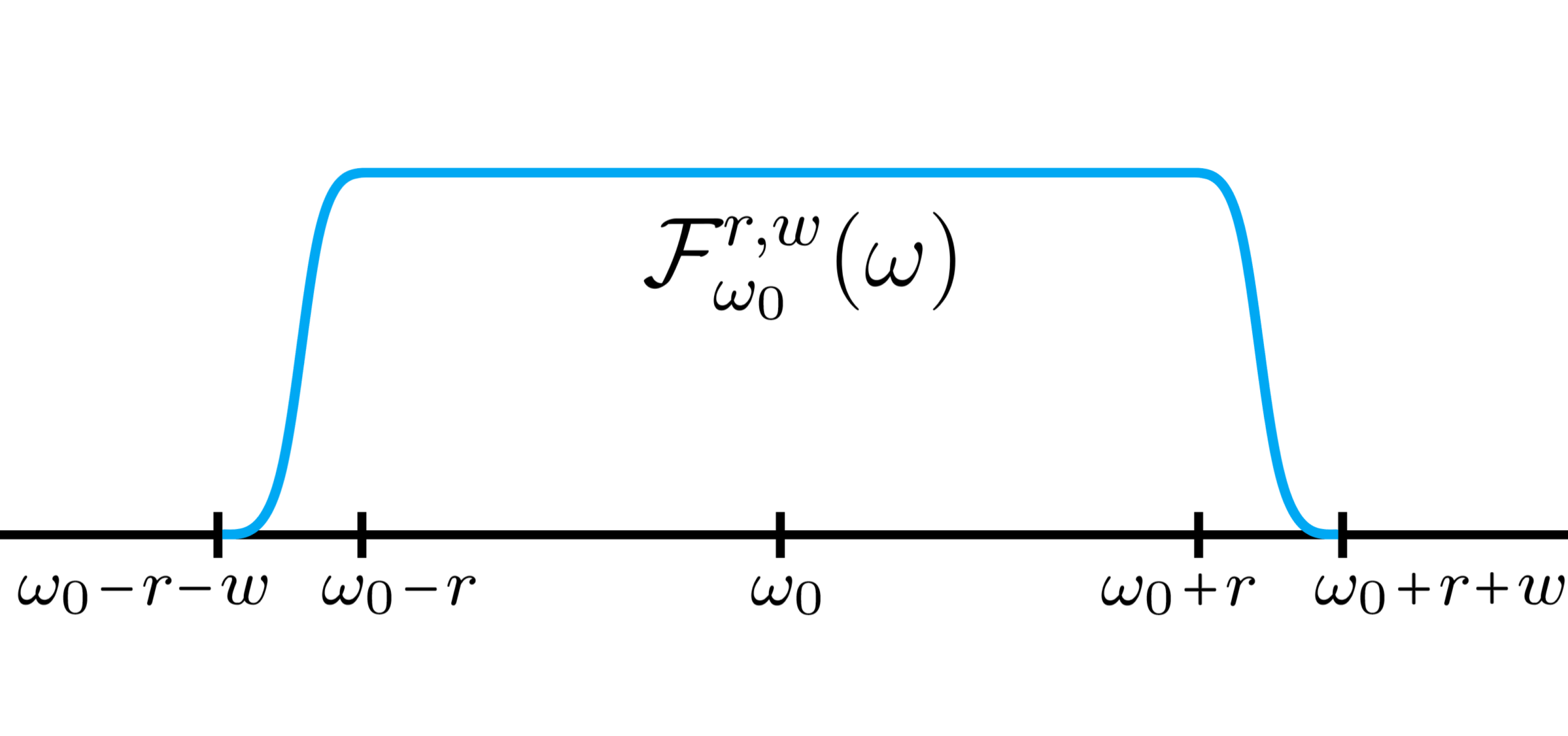}
    \caption{The graph of ${\mathcal F}_{\omega_0}^{r,w}(\omega)$.}
\end{figure}

Note that although we choose $\mathcal F^{r,w}_{\omega_0}$ to be smooth, we could pick it to be \emph{smooth enough}. In this case, terms like ``faster than any polynomial'' will become polynomial growth of a certain order. The modification\footnote{Compare this to how \cite{KS} characterizes Hastings' result.} changes $\delta^\gamma E(1/\delta)$ to $\delta^{\gamma - \epsilon}$ and how smooth that $\overline{\mathcal F}$ needs to be can be seen at the very end of the proof of the proposition where we put together the estimates for $T(l)$ and $G(l)$. This then eliminates the very slowly growing functions $E(1/\delta), F, T, G$ with much more explicit power functions.

We now return to the construction. Below, we choose $F(L)$ to grow slower than any power of $L$ and so that our Lieb-Robinson estimates work out. For $\beta_1 \in (0,1]$, let $n_{win} = \lceil L^{\beta_1}/F(L) \rceil$ writing $\omega
(i) = -1 + \kappa i$, where $\kappa = 2/n_{win}$, for $i = 0, \dots, n_{win}$. We then let $I_i$ be the interval of radius $\kappa$ centered at $\omega(i)$. So, the graphs of the ${\mathcal F}_{\omega(i)}^{0,\kappa}$,  are centered at $\omega(i)$ with radius $\kappa$,  have no ``flat part'', and form a partition of unity covering the interval $[-1,1]$.

\begin{defn} Let $S_1$ be a matrix whose columns $v_1, \dots, v_d$ are an orthonormal basis of $\mathcal V_1$. Let $\tau_i= [{\mathcal F}^{0,\kappa}_{\omega(i)}(J)]S_1$ and $Z_i = E_{[0, \lambda_{min}]}(\tau_i^\ast \tau_i)$. 
We define $\lambda_{min}$ later. Then let $\mathcal X_i$ be the range of $\tau_i(1-Z_i)$.

\end{defn}

We now state some estimates that will be useful later.
For $j = 0, 1$, the functions ${\mathcal F}_0^{jw, w}(t)$ are step functions that are equal to $1$ on $[-jw, jw]$ and zero outside $(-(j+1)w, (j+1)w)$, where we imagine $w$ being large. Now, note that
\[\int_{|k| \geq c} |\hat{\mathcal F}_{\omega_0}^{jw, w}(k)|dk = \int_{|wk| \geq cw} |w\hat{\mathcal F}_0^{j, 1}(wk)|dk = \int_{|k|\geq cw} |\hat{\mathcal F}_0^{j, 1}(k)|dk.\]
This type of argument gives us several identities:

\begin{align}\label{invariant under scaling}
\|\hat{\mathcal F}_{\omega_0}^{jw, w}\|_{L^1(\R)} = \|\hat{\mathcal F}_0^{j, 1}\|_{L^1(\R)},\end{align}
\begin{align}\label{scaling under scaling}
\|k\hat{\mathcal F}_{\omega_0}^{jw, w}(k)\|_{L^1(\R)} = \frac{1}{w}\|k\hat{\mathcal F}_0^{j, 1}(k)\|_{L^1(\R)},\end{align}
\begin{align}\label{defining G(l)}
\int_{|k| \geq \frac{l}{5e^2}} |\hat{\mathcal F}_0^{\frac{G(l)}{2l}, \frac{G(l)}{2l}}(k)|dk = \int_{|k|\geq \frac{G(l)}{10e^2}} |\hat{\mathcal F}_0^{1, 1}(k)|dk,
\end{align}
and
\begin{align}\label{defining F(L)}
\int_{|k| \geq \frac{L-1}{2e^2}} |\hat{\mathcal F}_0^{0, \kappa}(k)|dk = \int_{|k|\geq (L-1)/(e^2n_{win})} |\hat{\mathcal F}_0^{0, 1}(k)|dk.
\end{align}
Now, $(L-1)/(e^2n_{win}) = (L-1)/(e^2\lceil L^{\beta_1}/F(L) \rceil) \sim L^{1-\beta_1}F(L)/e^2$. We will later pick $\beta_1 = 1$ and $F$ to be slowly increasing in a way that equation (\ref{defining F(L)}) decreases faster than any power of $L$. We likewise pick $G(l)$ to be increasing slower than any power of $l$ in a way that equation (\ref{defining G(l)}) is decreasing faster than any polynomial of $l$. We also assume that $G\geq 2$.

Now, set
\begin{align}\label{defining S(l)}
S(L) = \int_{|k|\geq \frac{L-1}{e^2n_{win}}} |\hat{\mathcal F}_0^{0, 1}(k)|dk + \|\hat{\mathcal F}_0^{0, 1}\|_{L^1(\R)}e^{-\frac{L-1}{2}}
\end{align}
and
\begin{align}\label{defining T(l)}
T(l) = 2\int_{|k|\geq \frac{G(l)}{10e^2}} |\hat{\mathcal F}_0^{1, 1}(k)|dk + 3\|\hat{\mathcal F}_0^{1, 1}\|_{L^1(\R)}e^{-\frac{l}{5}}.
\end{align}
Note the implicit dependence of $n_{win}$ on $L$ in the definition of $S(L)$.
The intent of these definitions and the following estimates is to take advantage of the Lieb-Robinson estimates gotten in Section \ref{Lieb-Robinson section}.

\section{The use of Lin's theorem in Lemma 6.1}\label{use of Lin's theorem}
The following result is Lemma 3 from \cite{Hastings}. This is where Lin's theorem is used.

\begin{lemma} \label{application of Lin's theorem}
For any $\epsilon \in (0,1)$, there is a $\delta > 0$ such that if $A, B \in M_n(\C)$ are self-adjoint  with $\|A\|, \|B\| \leq 1$ and $\|[A,B]\| \leq \delta$ then there is a projection $P$ such that $\|[P,B]\| \leq \epsilon$ and $E_{[-1, -1/2]}(A)\leq P\leq 1-E_{[1/2,1]}(A)$.
\end{lemma}
\begin{remark}
Although we will only use this result for only one $\epsilon \in (0,1)$, the proof uses the ``For all $\epsilon > 0$ there exists a $\delta > 0$'' aspect of Lin's theorem for only one $\epsilon \in (0,1/22)$.
\end{remark}

This proof is a simplification of the proof in \cite{Hastings} using arguments from the proof of Theorem 3.2 in \cite{Davidson}.
\begin{proof} 
Let $A', B' \in M_n(\C)$ be any commuting self-adjoint matrices. Define $E=$ $E_{[-1,-1/2]}(A)$, $G = E_{[-1,1/2)}(A)$ and $P' = E_{(-\infty, 0)}(A')$. We have 
\[\|[P',B]\| = \|[P', B-B']\| \leq 2\|B-B'\|.\]
Also, by the Davis-Kahan theorem,
\[\|E{P'}^\perp\| = \|E_{[-1,-1/2]}(A)E_{[0,\infty)}(A')\| \leq 2\|A-A'\| \]
and
\[\|P'G^\perp\| = \|E_{(-\infty,0)}(A')E_{[1/2,1]}(A)\| \leq 2\|A-A'\|.\]

Then we apply Lemma \ref{projection lemma} to get that there is a projection $P$ such that $E \leq P \leq G$ with $\|P-P'\| \leq 10\|A-A'\|$ and
\[\|[P,B]\| \leq \|[P-P',B]\| + \|[P',B]\| \leq 20\|A-A'\| + 2\|B-B'\|.\]
Now by Lin's theorem, there is a $\delta > 0$ such that if $\|[A,B]\| \leq \delta$ we can pick $A', B'$ such that $\|A-A'\|, \|B-B'\| \leq \epsilon/22$ so that the projection $P$ satisfies the conditions of the lemma.
\end{proof}

\section{Many subspaces}\label{Many subspaces} 
Recall the spaces $\mathcal X_i$ from above. Spaces that are indexed with nonconsecutive indices are orthogonal, but otherwise we have no control over the orthogonality of these spaces. As mentioned at the end of Section \ref{General Approach}, what we want to do is be able to write $x \in \mathcal X_0 + \cdots + \mathcal X_{n_{win}-1} = \mathcal X$ as $\sum_{i=0}^{n_{win}}x_i,$ $x_i \in \mathcal X_i$ with $|x|^2$ comparable to $\sum_i |x_i|^2$. We cannot guarantee that this is possible, so we attempt to work around it. 

Make spaces $\mathcal R_i$ of the same dimension of $\mathcal X_i$ and define  $A:\bigcup_i\mathcal R_i \to \bigcup_i\mathcal X_i$ to be the natural identification. 
We then define $ \mathcal R$ as the (exterior) direct sum $ \bigoplus_i \mathcal R_i$ and extend $A: \mathcal R \to \mathcal X$.  

Note the key property that $A|_{\mathcal R_i}$ is an isometry. However, $A$ is not necessarily an isometry because the  subspaces are $\mathcal X_i$ are nonconsecutively orthogonal. We now estimate $\|A\|$. For any $r \in \mathcal R$ there are orthogonal $r_i\in \mathcal R_i$ such that $r = \sum_i r_i$. Then by Equation (\ref{constant 2}),
$|Ar|^2 \leq 2 \sum_i|Ar_i|^2 =2|r|^2$,
so $\|A\| \leq \sqrt2$.

Let $x \in \mathcal X$ with $r_i \in \mathcal R_i$ such that $x_i = Ar_i$. What we are looking for is a way to  make $|x|^2 = |\sum_i A r_i|^2 = |A(\sum_i r_i)|^2$ comparable to $\sum_i |x_i|^2 = \sum_i |r_i|^2 =  |\sum_i r_i|^2$. In the previous paragraph we showed that by nonconsecutive orthogonality, $|x|^2 \leq 2\sum_i|x_i|^2$. However, we need some control of the nonorthogonality between consecutive spaces in order to obtain some sort of reverse inequality. Indeed, if the $\mathcal X_i$ are linearly dependent this goal is hopeless.

\vspace{0.05in}

Just as before with $\mathcal F_{\omega(i)}^{0, \kappa}(J)\mathcal V_1$, in order to obtain a lower bound for $|A(\sum_i r_i)|^2$ we involve ourselves with finding a way to approximately ``cut'' the subspace $E_{\lambda \, small}(A^\ast A)$ out of our space $\mathcal R$. 
We will form a subspace $\mathcal U$ approximately out of spaces $\mathcal N_i$ so that $ E_{\lambda\, small}(A^\ast A)$ is approximately encompassed by the span of the $\mathcal N_i$ and then choose $\mathcal W = A\mathcal U^\perp$, which will be a subspace of $\mathcal X$.

Because we would like to use the fact that the $\mathcal X_i$ are approximate eigenspaces, we might hope that $E_{\lambda \, small}(A^\ast A)$ almost breaks down into an orthogonal sum of subspaces $\mathcal N_i$ of $\mathcal R_i$. 
It is not clear that this is possible, but what we can show is that by merging long chains of the $\mathcal R_i$ together into subspaces $\mathcal Y_i$ that have considerable consecutive overlap we can get $E_{\lambda \, small}(A^\ast A)$ to be approximately the span of $\mathcal N_i\subset \mathcal Y_i$. Because we do not expect there to be orthogonality (or even linear independence) between the $\mathcal N_i$, we need to concern ourselves with representing each vector projected onto by $E_{\lambda\, small}(A^\ast A)$ in a manageable way. Compare this motivation to that of the remark of Section C of \cite{Hastings}. We now proceed to the details of the constructions.

Let $\rho = A^\ast A.$ Then $\rho$ is positive and is block tridiagonal: if $r_k \in \mathcal R_k$  for all $k$, $Ar_k \in \mathcal X_k$ so $(r_i, \rho r_j) = (Ar_i, Ar_j) = 0$, if $|i-j| > 1$. Also, $\|\rho\| = \|A^\ast A\| \leq 2$.
Note that because $A|_{\mathcal R_i}$ is an isometry, the blocks of $\rho$ on the diagonal are identity matrices. 

Now, we merge some of the spaces $\mathcal R_i$ to attempt to take advantage of the block tridiagonal nature of $\rho$. As a notational convenience, if $I$ is an interval, then let $\mathcal R_{I} = \bigoplus_{i\in I}\mathcal R_i$. Define $\mathcal X_{I}$ likewise.  
Let $n_{b} \sim L^{\beta_0}$ be the number of ``superblocks'' that we form most of which have length $l_b$, which we will choose to be a multiple of $4$. That is, we define $k_b=\lfloor (n_{win}+1)/l_{b}\rfloor-1$ and then:
\[n_b = \left\{\begin{array}{ll}k_b & k_b \mbox{ is odd}.\\
k_b-1 & k_b \mbox{ is even}\end{array}\right..\] 
This is so that $n_b$ is always odd.

For $1\leq i \leq n_{b}-1$, let $\mathcal Y_i$ be the $i$th ``superblock'' defined by $\mathcal Y_i = \mathcal R_{[(i-1)l_b, (i+1)l_b)}$. See Figure \ref{subspacesImage} below.
\begin{figure}[htp]  \label{F1}
    \centering
    \includegraphics[width=16.5cm]{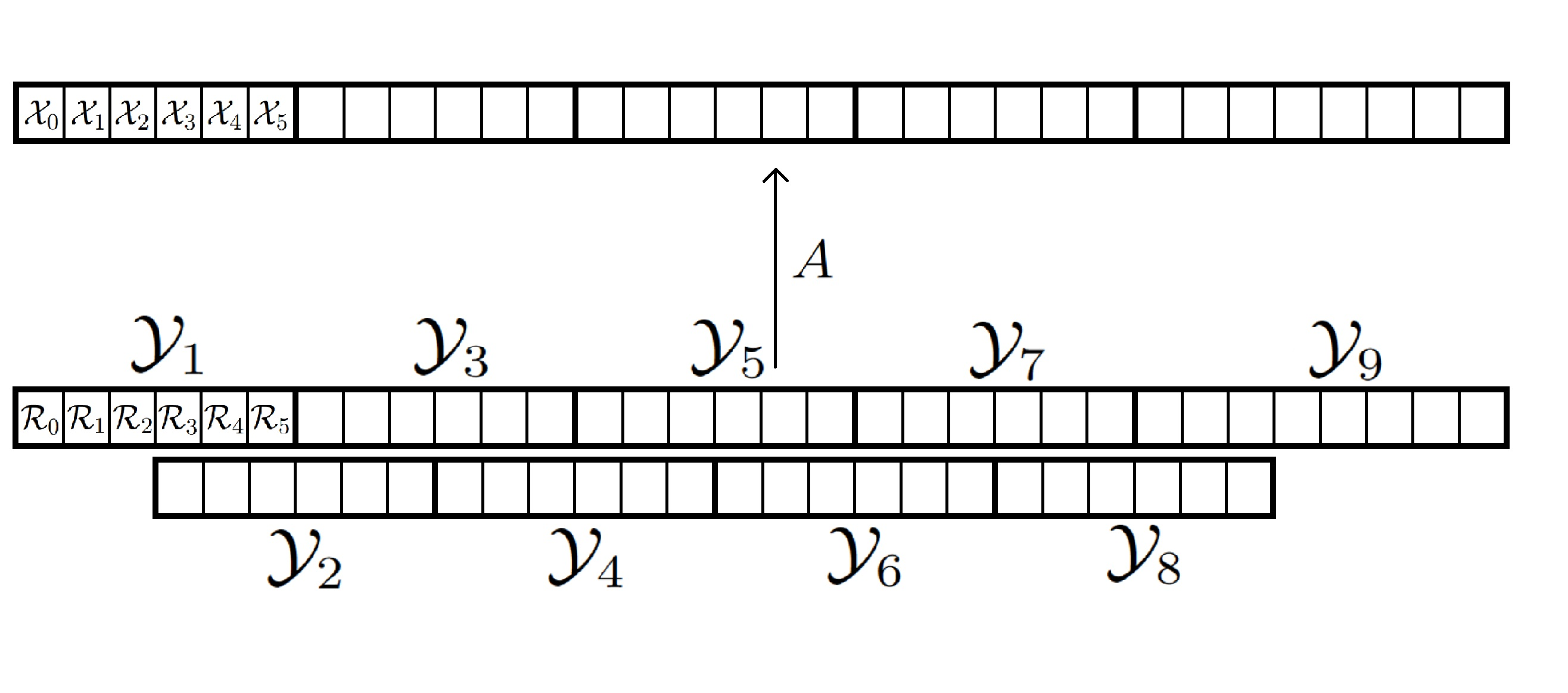}
    \caption{\label{subspacesImage}The subspaces $\mathcal X_i, \mathcal R_i, \mathcal Y_i$ are illustrated above with $32$ blocks, $l_b = 3$, $n_{win} = 31$, $n_b = 9$. (Note that in the proof, $l_b$ is assumed to be a multiple of $4$ but not in this illustration.) In this case we could define a subspace $\mathcal Y_{10}$ to be the last five blocks, however our construction choosing to make $\mathcal Y_{9}$ longer so that every subspace contains at least $6$ of the subspaces $\mathcal R_i$.
    Note the distinction between the first and last of the $\mathcal Y_i$. Note that the subspaces $\mathcal X_i$ are nonconsecutively orthogonal (and not even necessarily linearly independent), while the $\mathcal R_i$ are orthogonal.
    Also, the subspaces $\mathcal Y_i$ are nonconsecutively orthogonal and are displayed to show this.}
\end{figure}
Let $R_j$ be a projection in $\mathcal R$ onto $\mathcal R_j$, and for $1 < i \leq n_b-1$ let $Y_i$ project onto $\mathcal Y_i$, $Y_i'$ project onto $\mathcal R_{[(i-3/4)l_b, (i+3/4)l_b)}$, and $Y_i''$ project onto $\mathcal R_{[(i-1/2)l_b, (i+1/2)l_b)}$. See Figure \ref{YsubspacesImage} below.
Note that $Y_{i-1}'Y_i' = Y_i'Y_{i-1}'$ projects onto (a ``left'' subspace of the image of $Y_i'$) $\mathcal R_{[(i-3/4)l_b, (i-1/4)l_b)}$ and $Y_{i+1}'Y_i' = Y_i'Y_{i+1}'$ projects onto (a ``right'' subspace of the image of $Y_i'$)  $\mathcal R_{[(i+1/4)l_b, (i+3/4)l_b)}$. Note that there are $l_b/2$ many blocks between the images of $Y_{i-1}'Y_i'$ and $Y_{i+1}'Y_i'$. 
\begin{figure}[htp] 
    \centering
    \includegraphics[width=16.5cm]{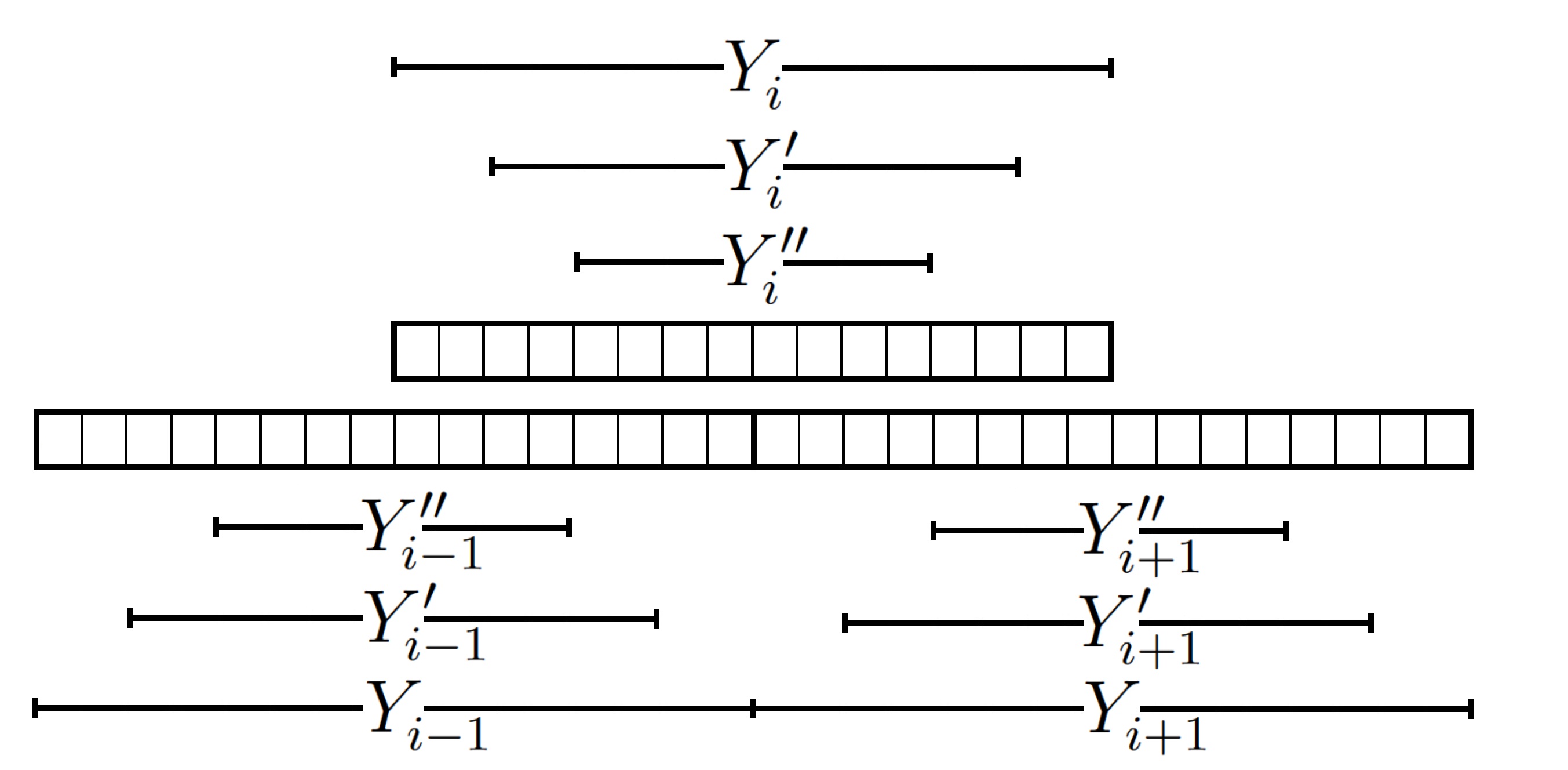}
    \caption{\label{YsubspacesImage}Typical consecutive projections $Y_i, Y_i', Y_i''$ are illustrated above with $l_b=8$.}
\end{figure}

Note that there are two distinguishable subspaces  $\mathcal Y_i$ that we address now, namely $\mathcal Y_1$ and $\mathcal Y_{n_{b}}$. We define these spaces and their respective projections $Y_i, Y_i', Y_i''$ to address the issue that $\mathcal Y_1$ and $\mathcal Y_{n_b}$ only intersect the other intervals on one side. The first subspace $\mathcal Y_1$ is a full length interval that only intersects other $\mathcal Y_i$ on its right side. And the last subspace $\mathcal Y_{n_b}$ has not yet been specified, but we do that now. The idea is that in order to apply the Lieb-Robinson estimates, we want all the intervals to have length $2l_b$ (and and the last interval has length at least $2l_b$ and less than $4l_b$), so we make $\mathcal Y_{n_b} = \mathcal R_{[(n_b-1)l_b, n_{win}]}$. For these two subspaces we now define the respective projections $Y_i, Y_i', Y_i''$ as we did above. 

For all $i$, $Y_i$ projects onto $\mathcal Y_i$. $Y_1'$ projects onto $\mathcal R_{[0,(1+3/4)l_b)}$. $Y_1''$ projects onto $\mathcal R_{[0,(1+1/2)l_b)}$. $Y_{n_b}'$ projects onto $\mathcal R_{[(n_b-3/4)l_b, n_{win}]}$. $Y_{n_b}''$ projects onto $\mathcal R_{[(n_b-1/2)l_b, n_{win}]}$. The point is that the $Y_i''$ form a resolution of the identity and there are $\frac{1}{4}l_b$ many blocks between the blocks projected onto by $Y_i''$ and the blocks not projected onto $Y_i'.$ To avoid multiple cases when dealing with these projections, let $Y_0$ project onto $\mathcal R_{[0,l_b)}$, $Y_0'$ project onto $\mathcal R_{[0, 3l_b/4)}$, $Y_{n_b+1}$ project onto $\mathcal R_{[n_bl_b, n_{win})}$, and $Y_{n_{b}+1}'$ project onto $\mathcal R_{[(n_b+1/4)l_b, n_{win})}$, so the ``left'' and ``right'' subspaces are well-defined for all $i = 1,\dots, n_b$.

\section{Construction and Properties of the spaces $\mathcal N_i$}\label{N_i}
We construct spaces $\mathcal N_i$ that essentially encompass the small eigenvalue eigenvectors of $\rho$, are subspaces of $\mathcal Y_i'$, and have various properties that we explore in this section. The major modifications in this section of the proofs from \cite{Hastings} were suggested in \cite{PC}.

\begin{defn}
Let $\rho_i = Y_i' \rho Y_i'$ be $\rho$ projected onto the range of $Y_i'$ and $\hat B_i$ a position operator for subspaces $\{\mathcal R_j\}$ of $\mathcal Y_i$ such that $\hat B_i$ is $-I$ restricted to the range of $Y_{i-1}'Y_i'$, is $I$ restricted to the range of $Y_{i+1}'Y_i'$, and linearly interpolates as multiples of the identities on blocks $\mathcal R_j$ with $\hat B_i$ on $\mathcal R_j$ being $\left(\frac{2}{l_b/2+1}[j-(i+1/4)l_b]+1\right)I$ for $(i-1/4)l_b \leq j < (i+1/4)l_b$. 

Let $\chi \in (0,1)$ be some constant that we will pick later. The following is part of our analogue of Lemma 4 from \cite{Hastings}.
\end{defn}
\begin{lemma}\label{defining N_i}
For $l_b$ large enough, there exist subspaces $\mathcal N_i$ (projected onto by $N_i$) such that $N_i \leq Y_i'$, where
\[E_{[0, G(l_b)/l_b]}(\rho_i) \leq N_i \leq Y_i'-E_{[2G(l_b)/l_b, \infty)}(\rho_i)\]
and
\[ \|[N_i, \hat B_i]   \| \leq 1-\chi < 1.\]

If $\delta$ is the constant from Lemma \ref{application of Lin's theorem} for $\epsilon = 1-\chi$, the condition on $l_b$ of $G(l_b) > 32C_{\mathcal F^{1,1}_0}/\delta$ can be taken as ``large enough''.
\end{lemma}
\begin{proof}
The idea is to apply Lemma \ref{application of Lin's theorem} to $\hat B_i$ and some function of $\rho_i$.

Because $\chi < 1$ is independent of $l_b$ and  that we require $G(l_b)$ go to infinity as $l_b$ does, we know that eventually $G(l) > 32C_{\mathcal F^{1,1}_0}/\delta$, where $\delta$ is a chosen value for $\epsilon = 1-\chi$ from Lemma \ref{application of Lin's theorem}.

Because $\rho_i$ and $\hat B_i$ both act on $\mathcal Y_i'$, $\rho_i$ is block tridiagonal, and $\hat B_i$ is a direct sum of multiples of identity matrices, with differences between consecutive multiples at most $\frac{4}{l_b+2}$, we see that $\|[\rho_i, \hat B_i]\| \leq \frac{16}{l_b+2}$. Let $f_{l_b}(x) = 1-2{\mathcal F}^{G(l_b)/l_b, G(l_b)/l_b}_0(x)$. Then 
\[E_{[-1,-1/2]}(f_{l_b}(\rho_i)) \geq E_{[0, G(l_b)/l_b]}(\rho_i)\] and
\[E_{[1/2, 1]}(f_{l_b}(\rho_i)) \geq E_{[2G(l_b)/l_b, \infty)}(\rho_i).\]
By equation (\ref{scaling under scaling}) and by Proposition \ref{B-R proposition} we get
\[ \|[f_{l_b}(\rho_i), \hat B_i]\| \leq 2C_{{\mathcal F}^{G(l_b)/l_b, G(l_b)/l_b}_0}\|[\rho_i, \hat B_i]\| \leq \frac{32l_b}{(l_b+2)G(l_b)}C_{{\mathcal F}^{1, 1}_0} <  \frac{32C_{{\mathcal F}^{1, 1}_0}}{G(l_b)}.\] So, for $l_b$ large enough we can apply Lemma \ref{application of Lin's theorem} to obtain $N_i$ such that 
\[E_{[0, G(l_b)/l_b]}(\rho_i) \leq E_{[-1,-1/2]}(f_{l_b}(\rho_i)) \leq N_i\perp E_{[1/2, 1]}(f_{l_b}(\rho_i)) \geq E_{[2G(l_b)/l_b, \infty)}(\rho_i)\]
and $\|[N_i, \hat{B}_i]\| < 1-\chi$.
\end{proof}

\begin{remark}
Instead of requiring $l_b$ to be large enough, we could have instead required that $G(l_b)$ have a lower bound that is large enough. However, how large would be left undetermined due to the statement of Lin's theorem applied to get Lemma \ref{application of Lin's theorem}. Similar remarks could be made elsewhere for estimates of our construction.
\end{remark}

\begin{defn}
Let $N^e = \sum_{i \,  even} N_i, N^o = \sum_{i \, odd}N_i$. Note that these are projections because the first item of the next result shows that the $\mathcal N_i$ are nonconsecutively orthogonal.
\end{defn}

We are interested in controlling the orthogonality of the spaces $\mathcal N_i$, the elements of the $\mathcal N_i$ as elements of eigenspaces of $\rho$, and the representations of eigenvectors (for small eigenvalues) of $\rho$ by the spaces $\mathcal N_i$.
The following is the completion of our analogue of Lemma 4 from \cite{Hastings}.
The last inequality in Item 3 was suggested by Hastings (\cite{PC}).
\begin{lemma}\label{properties of N_i}
For $l_b$ large enough defined by and for the $\mathcal N_i$ satisfying Lemma \ref{defining N_i} we have the following properties:
\begin{enumerate}
\item The spaces $\mathcal N_i$ are nonconsecutively orthogonal.
\item For any $v \in \mathcal N_i$, $0 \leq (v,\rho v) \leq \frac{2G(l_b)}{l_b}|v|^2$.
\item For $v$ in the range of $E_{[0, 1/l_b]}(\rho)$, there are $n_i \in \mathcal N_i$ such that 
\[|v - \sum_i n_i| \leq T(l_b)\sqrt{n_{b}}|v|,\]
\[\sum_i|n_i|^2 \leq  |v|^2,\]
and
\[\left|\sum_{i \, even}n_i\right|^2 \leq  |N^e v||v| + T(l_b)\sqrt{n_b}|v|^2.\]
\item $\|Y_{i+1}'N_{i}Y_{i-1}'\| \leq 1/2 - \chi/2$.
\end{enumerate}
\end{lemma}
\begin{proof}
Item 1 follows because $\mathcal N_i \subset \mathcal Y_i$ and the $\mathcal Y_i$ are nonconsecutively orthogonal by construction.

Item 2 is trivial once we remember that $\rho_i = Y_i'\rho Y_i'$, $N_i \leq E_{[0, 2G(l_b)/l_b]}(\rho_i)$, and $N_i\leq Y_i'$ so $(v, \rho v) = (v, \rho_i v)$.

We prove Item 4. Because  $ Y_{i-1}'Y_i' = E_{\{-1\}}(\hat B_i), Y_{i+1}'Y_i' = E_{\{1\}}(\hat B_i)$, and $N_i \leq Y_i'$, by Proposition \ref{comm-proj} we obtain
\begin{align*}\|Y_{i+1}'N_iY_{i-1}'\| &= \|(Y_{i+1}'Y_i')N_i(Y_i'Y_{i-1}')\| = \|E_{\{1\}}(\hat B_i)N_iE_{\{-1\}}(\hat B_i)\| \\
&\leq \frac{\|[N_i, \hat B_i]\|}{2} \leq \frac{1-\chi}{2} .
\end{align*}

We now prove Item 3. Recall that $G \geq 2.$ Let $D ={\mathcal F}_0^{G(l_b)/2l_b, G(l_b)/2l_b}(\rho),$ \\ $D_i' = {\mathcal F}_0^{G(l_b)/2l_b, G(l_b)/2l_b}(\rho_i)$. Note that $0 \leq D \leq E_{[0,G(l_b)/l_b]}(\rho)$ and $0 \leq D_i' \leq E_{[0,G(l_b)/l_b]}(\rho_i)\leq N_i$. We will apply a Lieb-Robinson estimate to $D, D'$ to obtain the $n_i$ with the desired properties.

Because $G \geq 2$ so that ${\mathcal F}_0^{G(l_b)/2l_b, G(l_b)/2l_b}(t) = 1$ for $ |t| \leq G(l_b)/2l_b \geq 1/l_b$, we have that for $v$ in the range of $E_{[0,1/l_b]}(\rho)$, $Dv = v$. If $v_i = Y_i''v$, then $v = \sum_i v_i$ is an orthonormal decomposition. 

Let $S'$ be the set of eigenvalues of $\hat{B}_i$ on $Y_i'$ and $S''$ the set of eigenvalues of $\hat{B}_i$ on $Y_i''$.
Now, $\rho_i = Y_i'\rho Y_i'$ is tridiagonal with respect to the eigenspaces of $\hat{B}_i$, so $\rho_i$ is finite range with distance less than $\Delta = 5/(l_b + 2)$. 
Because $\operatorname{dist}(S'', \R\setminus S') =  \frac{l_b}{4}\left(\frac{4}{l_b + 2}\right)$, we see that $\operatorname{dist}(S'', \R\setminus S') /\Delta \geq l_b/5$. Hence using Corollary \ref{operator Lieb-Robinson modified} along with (\ref{invariant under scaling}), (\ref{defining G(l)}) and (\ref{defining T(l)}), we have
\[\|(D-D_i')Y_i''\| \leq T(l_b).\]

So, set $n_i = D_i'v_i$. Then because $Dv =v$,
\[|v - \sum_i n_i| = |\sum_i (Dv_i - D_i'v_i)|\leq T(l_b)\sum_i|v_i|\leq T(l_b)\sqrt{n_{b}}|v|\]
and
\[\sum_i|n_i|^2 = \sum_i|D_i'v_i|^2 \leq \sum_i|v_i|^2 = |v|^2.\] 

When restricting to just the even indices we get a different bound using the orthogonality of $n_i = D_i'v_i \in \mathcal N_i$ for even $i$. Recall that $D$ is self-adjoint and the $v_i$ are orthogonal so:
\begin{align*}
\sum_{i\, even}&|n_i|^2 \leq \sum_{i\, even}|v_i|^2 = \left(v, \sum_{i\, even}v_i\right) = \left(Dv, \sum_{i\, even}v_i\right) = \left(v, \sum_{i\, even}Dv_i\right) \\
&\leq \left|\left(v, \sum_{i \, even} D_i'v_i \right)\right| 
+ \left|\left(v, \sum_{i \, even}(Dv_i - D_i'v_i) \right)\right| \\
&\leq \left|\left(v, \sum_{i \, even} n_i \right)\right|  + \sum_{i \, even}|(v, Dv_i - D_i'v_i)|
\leq \left|\left(N^ev, \sum_{i \, even} n_i \right)\right|  + T(l_b)|v|\sum_{i \, even}|v_i| \\
&\leq  |N^e v|\left|\sum_{i \, even}n_i\right| + T(l_b)\sqrt{n_b}|v|^2 = |N^e v|\left(\sum_{i \, even}|n_i|^2\right)^{1/2} + T(l_b)\sqrt{n_b}|v|^2\\
&\leq|N^e v|\sqrt{\sum_{i \, even}|v_i|^2} + T(l_b)\sqrt{n_b}|v|^2
\leq |N^e v||v| + T(l_b)\sqrt{n_b}|v|^2.
\end{align*}

\end{proof}

\begin{remark}
The last inequality in Item 3 gives us control of the norm of $\sum_{i\, even} n_i$ if we have some nontrivial bound for the norm of $N^ev$.

Although a similar result is true for $\sum_{i \, odd}|n_i|^2$, we will not use it.
\end{remark}
\begin{remark}
Note that because $Y_{i+1}'\perp Y_{i-1}'$ the inequality $\|Y_{i+1}N_iY_{i-1}\| \leq 1/2 - \chi/2$ is a slight improvement of the general property of that if $P, Q$ are any projections then $\|(1-P)QP\| \leq 1/2$.

This property follows by applying Jordan's lemma to reduce to the $2$-dimensional case where one has $P$ projecting onto the first basis vector $e_1$, $1-P$ projecting onto the second basis vector $e_2$, and $Q$ projecting onto a vector $v = \cos\theta e_1 + \sin\theta e_2$. A calculation then shows that $\|(1-P)QP\| = |\cos\theta\sin\theta| = \frac{1}{2}|\sin2\theta|$. 
\end{remark}

We now proceed to define the key subspaces $\mathcal N_i'$ as in the discussion before Lemma 5 of \cite{Hastings}.
\begin{defn}Fix some $\eta \in (0,\chi/4)$. Let $i$ be odd. Apply Jordan's lemma to $N_i$ and $N^e$ to get an orthonormal basis $\{n_s\}$ of $\mathcal N_i$ such that $(N^en_s, n_t) = 0$ if $s \neq t$.
Then let $\mathcal N_i'$ be the subspace of $\mathcal N_i$ generated by the $n_s$ such that $|N^en_s|^2 \leq 1/2 + \eta$. Then let $N_i'$ project onto $\mathcal N_i'$ and $N' = \sum_{i\, odd}{N_i}'$.
\end{defn}

\begin{remark}\label{N_i' property}
Note that $\mathcal N_i'$ is not defined for $i$ even, so $N'$ unambiguously projects onto the orthogonal sum of $N_i'$ for $i$ odd.
We also have that $N_j N^o N_i = N_j N^e N_i = N_j N' N_i = 0$ for $|i-j|\geq 3$.
By our definitions,  if $n_i' \in \mathcal N_i'$ is a unit vector then we can express it as $n_i' = \sum_s c_s n_s$, where the $n_s$ are as above. Then by definition, \begin{align}\label{N^en_i'}|N^en_i'|^2 = |\sum_{s} c_sN^en_s|^2 = \sum_{s} |c_sN^en_s|^2 \leq (1/2+\eta)\sum_{s} |c_s|^2 = 1/2+\eta\end{align}
and consequently
\begin{align}\label{(1-N^e)n_i'}|(1-N^e)n_i' |^2\geq 1/2 - \eta.\end{align}
Note that if the unit vector $n_i \in \mathcal N_i \ominus \mathcal N_i'$ then $n_i$ is expressed as linear combination of the $n_s$ with $|N^e n_s|^2 > 1/2 + \eta$ and hence the opposite inequalities hold: \begin{align}\label{N^en_i}|N^en_i|  > \sqrt{1/2+\eta}\end{align}
and 
\begin{align}\label{(1-N^e)n_i}|(1-N^e)n_i |< \sqrt{1/2 - \eta}.\end{align}
\end{remark}

See that the definition of $N_i'$ is intended to remove a part of $\mathcal N_i$ that has large projection onto $\mathcal N_{i-1}\oplus \mathcal N_{i+1}$. 
This does not quite seem to be enough control of the orthogonality to give a result like Equation (\ref{degree of orthogonality}), but, along with $\|Y_{i+1}'N_iY_{i-1}'\|< 1/2 - \chi/2$, it will suffice for our purposes.

\begin{defn}
Let $\mathcal U \subset \mathcal R$ be defined as the orthogonal complement of the span of $\bigcup_{i \, even}\mathcal N_i \cup \bigcup_{i\, odd}\mathcal N_i'$. Then let $\mathcal W = A\mathcal U \subset \mathcal X$ and $P$ project onto $\mathcal W$.
\end{defn}

\begin{remark}\label{linindep}
$\mathcal U^\perp$ is expressed as the \emph{direct} sum $\bigoplus_{i\, odd} \mathcal N_{i+1} \oplus \mathcal N_{i}'$, because these subspaces are linearly independent. This is actually a consequence of the following lemma, because otherwise the lemma could not possibly provide any bound. So, despite linear independence being a consequence of the following lemma, we first provide a simple argument for this because it makes us more comfortable with dealing with linear combinations of vectors from these subspaces and because our argument is illustrative of the argument in the lemma.

We show that this is a direct consequence of the fact that $\|N_{i+2}'N_{i+1}N_i'\| < \frac12$ and how when constructing $\mathcal N_i'$ we ``cut out'' a part of $\mathcal N_i$ that is more parallel to the image of $N^e$ than its orthogonal complement. 

We show the 
linear independence as follows.
Suppose that $n_i \in \mathcal N_i$ for $i$ even and $n_i \in \mathcal N_i'$ for $i$ odd such that
\begin{align}\label{linind}\sum_{i\, even} n_i= -\sum_{i\, odd}n_i.\end{align}
By nonconsecutive orthogonality, Equation (\ref{linind}) shows that
\begin{align}\label{pythagequality}\sum_{i\, even} |n_i|^2=\left|\sum_{i\, even} n_i\right|^2= \left|\sum_{i \, odd}n_i\right|^2 =\sum_{i\, odd}|n_i|^2 .\end{align}
However, applying $1-N^e$ to Equation (\ref{linind}) also shows that
\[0=\left|\sum_{i\, even} (1-N^e)n_i\right|^2= \left|\sum_{i\, odd}(1-N^e)n_i\right|^2.
\]

Following the calculation below in Equation (\ref{Inner Product}), we see that for $i$ odd
\begin{align*}
|((1-N^e)n_i&,(1-N^e)n_{i+2})|=|(Y_{i+2}'N_{i+1}Y_i'(1-N^e)n_i,(1-N^e)n_{i+2})|\\
&\leq \|Y_{i+2}'N_{i+1}Y_i'\||(1-N^e)n_i||(1-N^e)n_{i+2}|.
\end{align*}
So, by Equation (\ref{semi-pythagorean}),
\[0=\left|\sum_{i\, odd}(1-N^e)n_i\right|^2\geq (1-2\|Y_{i+2}'N_{i+1}Y_i'\|)\sum_{i\, odd}|(1-N^e)n_i|^2.\]
So, $\|Y_{i+2}'N_{i+1}Y_i'\|<\frac12$ implies that $(1-N^e)n_i=0$ for $i$ odd. So, by Equation (\ref{(1-N^e)n_i'}), $n_i=0$ for $i$ odd and hence $n_i = 0$ for all $i$ by Equation (\ref{pythagequality}).
\end{remark}

The above remark reflects some of the ideas behind the following lemma: reformulate the result in terms of some linear combination, single out the contributing terms from the $\mathcal N_i'$ for $i$ odd using $1-N^e$, show that the norms of these terms satisfy the desired decaying condition, then extend this decay to the terms from the $\mathcal N_i$ for $i$ even. 

The following lemma is our adaptation of Lemma 5 of \cite{Hastings} with improvements to the proof suggested by Hastings (\cite{PC}).
\begin{lemma}\label{exponential decay} There are constants $C_1 = C_1(\chi,\eta) \geq 1$ and $\alpha = \alpha(\chi,\eta)\in (0,1)$ such that for any $y_i \in \mathcal Y_i$, there are (unique) $n^i_j \in \mathcal N_j$ for $j$ even and $n^i_j \in \mathcal N_j'$ for $j$ odd such that $U^{\perp}y_i = \sum_{j} n^i_j$ and $|n_j^i| \leq C_1\alpha^{|i-j|}|y_i|$.

Consequently, for $c_\alpha = (1 + \alpha + \alpha^{-1})$ and $C_2 = c_\alpha C_1, \|Y_j U Y_i\| \leq C_2\alpha^{|i-j|}.$
\end{lemma}
\begin{proof}
Note that for the proof of the first part, it will be important to consider different values of $|n_j^i|$ for the same unit vector $n_j^i/|n_j^i|$. So, for the proof we will write $n_j^i = a_j^im_j^i$, where $m_j^i$ is always a unit vector (even if $|a_j^i| = |n_j^i| = 0$).

We first prove the second result. Note that because $C_2 \geq 1/\alpha$, the result is trivial for $|i-j| \leq 1$. If $|i-j|\geq 2$ then $\mathcal Y_i$ and $\mathcal Y_j$ are orthogonal, so $Y_j U Y_i = -Y_jU^\perp Y_i$. Then with the first result, we know that for a unit vector $y_i \in \mathcal Y_i$ we can write $U^\perp y_i = \sum_k n_k^i$ where $n^i_k \in \mathcal N_k$ for $k$ even, $n^i_k \in \mathcal N_k'$ for $k$ odd, and $|n_k^i| \leq C_1\alpha^{|i-k|}$. Then 
\[|Y_jU^\perp y_i| = |\sum_k Y_jn_k^i| \leq |n_{j-1}^i|+ |n_j^i| + |n_{j+1}^i| \leq c_\alpha C_1\alpha^{|i-j|}.\]

So, we now prove the first part. Let $y_i \in \mathcal Y_i$. 
By the definition of $U^\perp$, we can write $U^\perp y_i$ as a linear combination of unit vectors $m^i_j \in \mathcal N_j$ for $j$ even and $m^i_j \in \mathcal N_j'$ for $j$ odd. 
Removing some elements from $\N \cap [1, n_b]$ we obtain a set $S$ such that $\{m_s^i\}_{s\in S}$ is linearly independent and still having $y_i$ in their span. 
Note that the statement that the vectors $m_j^i$ are unit vectors implicitly assumes that we are excluding indices $j$ where $\mathcal N_j=0$ for $j$ even or $\mathcal N_j'=0$ for $j$ odd, so some restriction of indices is necessary. 

What we want to do is to take advantage of the properties $\mathcal N_s'$. Because $N^e \leq U^\perp$, we will isolate the $a_s^i$ for $s$ odd by applying $(1-N^e)$ to our representation of $U^\perp y_i$. This gives
$(1-N^e)U^\perp y_i = \sum_{s \, odd}a_s^i(1-N^e)m_s^i$. 
We want to find relationships between the inner products of the terms $(1-N^e)m_s^i$ for $s$ odd and use them to obtain control of the $a^i_s$ for $s$ odd. We then extend this control to $s$ even.

We focus on $s, t$ odd. The first statement of Lemma \ref{properties of N_i} implies that $N_aN_b = 0$ if $|a-b| \geq 2$. Consequently $(1-N^e)m_s^i = (1-N_{s-1}-N_{s+1})m_{s}^i = m_{s}^i-N_{s-1}m_{s}^i-N_{s+1}m_{s}^i \in \mathcal N_{[s-1,s+1]}$. So, if $((1-N^e)m_s^i, (1-N^e)m_t^i) \neq 0$ then $\operatorname{dist}([s-1,s+1],[t-1,t+1]) \leq 1$ and since $s, t$ are odd, we have $|s-t|\leq 2$. In the case $s = t$, we have $|(1-N^e)m_s^i|^2 \geq 1/2 - \eta$ by Equation (\ref{(1-N^e)n_i'}).

For $|s-t| = 2$, we can assume that $t = s+2$. Then, roughly speaking, the only contribution to the inner product comes from the orthogonality of the ranges of $N_{s+1}N_s'$ and  $N_{s+1}N_{s+2}'$ which will give some control by the last statement of Lemma \ref{properties of N_i}.
In more detail, we have the orthogonal decomposition, $m_s^i = Y_{s-1}'m_s^i + Y_{s+1}'m_s^i + (1-Y_{s-1}'-Y_{s+1}')m_s^i$ from which we can set $c_s = |Y_{s-1}'m_s^i|$ and $d_s = |Y_{s+1}'m_s^i|$ so that $c_s^2 + d_s^2 \leq 1$. Recall that $t = s+2$ so $N_{s+1}=N_{t-1}$ and using the fact that the $Y_a'$ commute we see that
\begin{align}
|((1-N^e)m_s^i,&(1-N^e)m_{s+2}^i)| \nonumber=|(1-N^e)m_s^i,  m_{s+2}^i)| \\
&= ((1-N_{s-1}-N_{s+1})m_s^i,  m_{s+2}^i)|
= |(N_{s+1}m_s^i,  m_{s+2}^i)|\nonumber\\
&=|(Y_{s+2}'Y_{s+1}'N_{s+1}Y_{s+1}'Y'_{s}m_s^i,  m_{s+2}^i)|
\nonumber\\
&=|((Y_{s+2}'N_{s+1}Y'_{s})Y_{s+1}'m_s^i,  Y_{s+1}'m_{s+2}^i)|\nonumber\\
&\leq d_sc_{s+2}\|Y_{s+2}'N_{s+1}Y_s'\| \leq (1/2 - \chi/2)d_sc_{s+2}.\label{Inner Product}
\end{align}

Then because  $(1-N^e)U^\perp y_i = \sum_{s \, odd}a_s^i(1-N^e)m_s^i$, for any odd $t$ we have
\[\left((1-N^e)m^i_t, (1-N^e)U^\perp y_i - \sum_{s \, odd}  a^i_s(1-N^e)m^i_s\right) = 0\] so
\begin{align}\label{MatrixEq}
((1-N^e)m^i_t, (1-N^e)U^\perp y_i) = \sum_{s \, odd} a^i_s((1-N^e)m^i_t, (1-N^e)m^i_s).
\end{align}

We express this as a matrix equation, but first list the odd elements $s$ in $S$ in increasing order $s_1, s_2, \dots$. We let $x$ be some small positive number, which we will choose later. We reformulate a scaled version of (\ref{MatrixEq}) by setting \[\vec v_k = \frac{2(1+x)}{1-2\eta}((1-N^e)m^i_{s_k}, (1-N^e)U^\perp y_i),\;\;
\vec a_l = a^i_{s_l},\] and
\[M_{k, l} = \frac{2(1+x)}{1-2\eta}((1-N^e)m^i_{s_k}, (1-N^e)m^i_{s_l}))\] to obtain $M \vec a = \vec v$. 
Once we show that $M$ is invertible, we will obtain a representation of $\vec a$ as $M^{-1}\vec v$ and then will use the properties of $M$ to get control of the $a_s^i$. 

Note that the scaling factor is so that the bounds on the diagonal entries of $M$ are much clearer. In particular, by construction $M$ is a tridiagonal self-adjoint matrix and we have the inequalities: $M_{k,k} \geq \frac{2(1+x)}{1-2\eta}(1/2-\eta) = 1+x\geq(c_{s_k}^2+d_{s_k}^2)+x$ and $|M_{k,k+1}| \leq \frac{(1+x)(1-\chi)}{(1-2\eta)}d_{s_k}c_{s_{k+1}}$ if $s_{k+1}=s_k+2$ and $M_{k,k+1}=0$ otherwise. 

Now, because $\eta<\chi/4$, we choose $x = \frac{\chi}{2-2\chi} > 0$ so $M_{k,k} \geq 1+x$ and $|M_{k,k+1}| \leq d_{s_k}c_{s_{k+1}}$.
By Lemma \ref{positive}, $M -xI$ is positive, so all the eigenvalues of $M$ are at least $x$. Proposition \ref{exponential decay of inverse} shows that there are constants $\tilde C>0, \tilde \alpha \in (0,1)$ such that $|(M^{-1})_{l,k}|\leq \tilde C \tilde{\alpha}^{|l-k|}$ and these constants only depend on $x$ and an upper bound on the spectrum of $M$, which can be taken to be $6(1+x)/(1-2\eta)$. 

Now, $N^e$ and $U^\perp$ commute and $U^\perp m^i_{s_k}=m^i_{s_k}$ so we see that \[U^\perp(1-N^e)m^i_{s_k} = (1-N^e)U^\perp m^i_{s_k}=(1-N^e) m^i_{s_k}\in \mathcal N_{[s_k-1,s_k+1]}\subset \mathcal Y_{[s_k-1,s_k+1]}.\]
Because $\vec v_k = \frac{2(1+x)}{1-2\eta}((1-N^e)m^i_{s_k}, U^\perp y_i)$ and $y_i \in \mathcal Y_i$, we see that $\vec v_k = 0$ for $|s_k - i| \geq 3$.  The increasing sequence of odd integers $s_1, s_2, \dots$ might have gaps. These gaps cause $M$ to be a block diagonal matrix, each block being a tridiagonal matrix. Then $M^{-1}$ also has this same block structure with exponential decay of the entries away from the diagonal. 
We illustrate this with the example that $s_1 = 5, s_2 = 7, s_3 = 11, s_4 = 13, s_5=15$ and so $M$ has the following block structure:\[\begin{pmatrix} M_{1,1} & M_{1,2} &&&\\ M_{2,1} & M_{2,2} &&&\\ &&M_{3,3}&M_{3,4}&\\ &&M_{4,3}&M_{4,4}&M_{4,5}\\&&&M_{5,4}&M_{5,5}\end{pmatrix}.\]

When bounding $a^i_{s_l}$ we can restrict to the block of $M$ in which $s_l$ lies. Let this block be $s_{a}, \dots, s_b$. Then for $a\leq k, k'\leq b$, $|s_{k}-s_{k'}|=2|k-k'|$. We also have the inequality
\[||s_l-i|-|s_l-s_k||\leq |s_k-i|.\]
So,
\begin{align*}
|a^i_{s_l}|&\leq \sum_{a\leq k \leq b}|(M^{-1})_{l,k}||\vec v_k|\leq \frac{2(1+x)}{1-2\eta}|y_i|\sum_{\substack{a\leq k \leq b\\|s_k-i|\leq 2}}|(M^{-1})_{l,k}|\leq\\
& \frac{2(1+x)}{1-2\eta}\tilde C|y_i|\sum_{|s_k-i|\leq 2} \tilde{\alpha}^{|s_l-s_k|/2}\leq \frac{2(1+x)}{1-2\eta}3\tilde C|y_i|\alpha^{-1} \tilde{\alpha}^{|s_l-i|/2}=:\tilde{\tilde C}\alpha^{|s_l-i|}|y_i|.
\end{align*}

Now, we can use the smallness of the odd coefficients to deduce the smallness of the even coefficients. If $j$ is even and $|i-j| \geq 2$, then because $Y_j y_i = 0$ and the $m_k^i$ are unit vectors, we see that
\begin{align*}
|a_j^i| &=\left|Y_j(\sum_{s\, even}a_{s}^im_{s}^i)\right| =\left|Y_j(y_i - \sum_{s\, odd}a_{s}^im_{s}^i)\right| = |Y_j( a_{j-1}^im_{j-1}^i + a_{j+1}^im_{j+1}^i)| \\
&\leq \tilde{\tilde{C}}(\alpha^{|i-j|+1} + \alpha^{|i-j|-1})|y_i| = \tilde{\tilde{C}}(\alpha + \alpha^{-1})\alpha^{|i-j|}|y_i|.
\end{align*}
Setting $C_1 = \tilde{\tilde{C}}(\alpha + \alpha^{-1})$, we obtain the result.
\end{proof}

\section{Properties of $\mathcal U$ and $\mathcal W$}\label{U and W}
We prove some properties of $\mathcal U$ and $\mathcal W$ here that will be used to complete the proof of Lemma \ref{main lemma}. We begin with the first property that was motivated by Section \ref{Many subspaces}.

This result is stated (albeit with a different exponent of $l_b$) in \cite{Hastings} and the proof is based on that in \cite{Hastings}, supplemented by suggested adjustments from \cite{PC}.
\begin{lemma}\label{(u, rho u)}
There is a(n explicit) constant $C_3$ depending only on $\eta$ such that for $l_b$ large and $u \in \mathcal U$, \[|Au| \geq \sqrt{\frac{1}{C_3l_b}}|u|.\]
It suffices that $l_b$ is large enough that $n_bT(l_b) <( (1-\sqrt{1-\eta})/6)^2$ and also large enough specified by Lemma \ref{defining N_i}.
\end{lemma}
\begin{proof}
Recall that $|Au|= \sqrt{(u, \rho u)}$. It suffices to show that there is a constant $C < 1$ such that for $w$ in the range of $E_{[0,1/l_b]}(\rho)$, $|Uw| \leq C|w|$. This is because if this is so then using the orthonormal eigendecomposition of $\rho$ expressed in bra-ket notation $\rho = \sum_{\lambda}\lambda |v_{\lambda}\rangle\langle v_{\lambda}|$ gives
\begin{align}
(u, \rho u) &\geq \nonumber  \frac{1}{l_b}\sum_{\lambda > 1/l_b}|(v_\lambda,u)|^2 = \frac{1}{l_b}|E_{(1/l_b,\infty)}(\rho)u|^2 = \frac{1}{l_b}(|u|^2 -|E_{[0,1/l_b]}(\rho)u|^2 ) \\
&= \frac{1}{l_b}(|u|^2 - \nonumber \max_{|v|=1}|(v,E_{[0,1/l_b]}(\rho)u)|^2 ) = \frac{1}{l_b}\left(|u|^2 - \max_{E_{[0,1/l_b]}(\rho)v = v, |v|=1}|(v,u)|^2 \right) \\
&\geq \frac{1}{l_b}(|u|^2 - C^2|u|^2 ) = \frac{1-C^2}{l_b}|u|^2. \label{1230}
\end{align}

If $w$ is in the range of $E_{[0,1/l_b]}(\rho)$, we have $n_i \in \mathcal N_i$, satisfying Item 3 of Lemma \ref{properties of N_i}. Because $U^\perp \geq N^e$ and $U^\perp \geq N_i'$, we see that $U n_i = 0$ for $i$ even and $Un_i = U(N_i-N_i')n_i$ for $i$ odd. Now, like in the proof of Lemma \ref{exponential decay}, we will obtain control of the terms $n_i$ for $i$ odd and then extend this control to the rest of terms. The difficulty is that there is no clear way to obtain an approximate decomposition as below of $w \sim w^{even} + w^{odd}$ that are approximately (for $l_b$ or $n_b$ large) orthogonal. So, we work around this.

Let $w^{even} = \sum_{i \, even}n_i$ and $w^{odd} = \sum_{i\, odd}n_i$, so $w$ is approximately equal to $w^{even}+w^{odd}$, as $|w - w^{even}-w^{odd}|\leq T(l_b)\sqrt{n_b}|w|$. Also because $Uw^{even} = 0$,
\[|Uw - Uw^{odd}| = |U(w -w^{even}-w^{odd})| \leq T(l_b)\sqrt{n_b}|w|.\] For $s = 1, 3$, let $[s]$ be the set of all natural numbers equivalent to $s$ modulo four. Define $w^1$ and $w^3$ by $w^s = \sum_{i \in [s]}(N_i-N_i')n_i$. So, $w^s$ is expressed as a series of orthogonal vectors and $U(w^1 + w^3) = U(\sum_{i\, odd}(N_i-N_i')n_i)=U(\sum_{i\, odd}n_i) = Uw^{odd}$.

Now because $(w^1,w^3) = 0$,
\begin{align*}
|U^\perp&(w^1 + w^3)|^2 = |U^\perp w^1|^2 + |U^\perp w^3|^2 + 2\Re(U^\perp w^1, w^3) \\
&\geq |U^\perp w^1|^2 + |U^\perp w^3|^2 - 2|(U w^1, Uw^3)| \geq |U^\perp w^1|^2 + |U^\perp w^3|^2 - |U w^1|^2 -| Uw^3|^2.
\end{align*}
By Equation (\ref{N^en_i}), for a unit vector $n \in \mathcal N_i \ominus \mathcal N_i'$ we know that $|U^\perp n|^2 \geq |N^e n|^2> 1/2 + \eta$. We can apply this to $n = (N_i-N_i')n_i$. If $i, j \in [s]$ are not equal then $|i-j|\geq 4$. So, if $n_i \in \mathcal N_i$ then $U^\perp (N_i-N_i')n_i = (N_{i-1}+N_{i+1})(N_i-N_i')n_i\in\mathcal N_{[i-1,i+1]}$ and hence $U^\perp(N_i-N_i') n_i \perp U^\perp(N_i-N_i') n_j$. This implies that 
\begin{align*}
|U^\perp w^s|^2 &=|\sum_{i \in [s]} U^\perp (N_i-N_i')n_i| ^2 =\sum_{i \in [s]} |U^\perp (N_i-N_i')n_i|^2\\
&\geq (1/2+\eta)\sum_{i \in [s]} |(N_i-N_i')n_i|^2=(1/2+\eta)|w^s|^2. \end{align*}

So because, $|U^\perp w^s|^2-|w^2|^2 = -|U w^s|^2$,
\begin{align*}
0 &\leq 2(|U^\perp w^s|^2-(1/2+\eta)|w^s|^2)=2|U^\perp w^s|^2-|w^s|^2-2\eta|w^s|^2\\
&=|U^\perp w^s|^2-|Uw^s|^2- 2\eta|w^s|^2.\end{align*} 
Hence we obtain $|U^\perp (w^1 + w^3)|^2 \geq 2\eta(|w^1|^2 + |w^3|^2) = 2\eta |w^1 + w^3|^2$, so 
\begin{align*}
|Uw| &\leq |U(w^1+w^3)|+T(l_b)\sqrt{n_b}|w| \leq \sqrt{1-2\eta}|w^1 + w^3|+T(l_b)\sqrt{n_b}|w| \\
&= \sqrt{1-2\eta}|w^{odd}|+T(l_b)\sqrt{n_b}|w|.
\end{align*}

Now, we will give an upper bound for $|Uw|$ in terms of $|w|$ based on the cases when $|N^ew|$ is small and when it is not really that small. When $|N^ew|/|w|$ is not too small, we can use 
\[|Uw| = |U(1-N^e)w| \leq |(1-N^e)w| = \sqrt{|w|^2 - |N^ew|^2}.\]
When $|N^ew|/|w|$ is small, we remember the third equation of Item 3 of Lemma \ref{properties of N_i} which gives $|w^{even}| \leq \sqrt{|N^ew||w| + T(l_b)\sqrt{n_b}|w|^2}$. Note that by the assumptions of the lemma, $T(l_b) < 1$. Hence, by what we have done:
\begin{align*}
|Uw| &\leq \sqrt{1-2\eta}|w^{odd}|+T(l_b)\sqrt{n_b}|w|\\
&\leq \sqrt{1-2\eta}|w-w^{even}|+(1+\sqrt{1-2\eta})T(l_b)\sqrt{n_b}|w|\\
&\leq \sqrt{1-2\eta}|w| +\sqrt{1-2\eta} \sqrt{|N^ew||w| + T(l_b)\sqrt{n_b}|w|^2} +2T(l_b)\sqrt{n_b}|w| \\
&\leq \sqrt{1-2\eta}(1 + \sqrt{|N^ew|/|w|})|w| +(\sqrt{1-2\eta}\sqrt[4]{n_b} + 2\sqrt{n_b})\sqrt{T(l_b)} |w|.
\end{align*}
Recall that $\eta < \chi/4 < 1/4$ and consider $p = \left(\sqrt{\frac{1-\eta}{1-2\eta}}-1\right)^2$. Then $p > 0$ and because $\eta < 3/7$, we have $p < 1$. Then for $|N^ew| \leq p|w|$, the latter estimate gives $|Uw| \leq \sqrt{1-\eta}|w| +  3\sqrt{n_bT(l_b)} |w|$. If $|N^ew| > p|w|,$ then the first estimate gives
$|Uw| \leq \sqrt{1-p^2}|w|$.

So, picking $l_b$ large enough that $3\sqrt{n_bT(l_b)} < (1-\sqrt{1-\eta})/2$, we obtain $C_3$ and $C = \max\left(\frac{1+\sqrt{1-\eta}}{2}, \sqrt{1-p^2}\right)$ using (\ref{1230}).
\end{proof}

Recall that $P = P_{\mathcal W}$.
Note that, as discussed above, this result implies that $A$ is injective on $\mathcal U$. Other consequences are the following inequalities for vectors in $\mathcal W = A\mathcal U$ from Section 5 of \cite{Hastings}.
\begin{lemma}\label{estimates for W}
For $l_b$ large (determined by Lemma \ref{defining N_i} and Lemma \ref{(u, rho u)}), we have:
\begin{enumerate}
\item For $w \in \mathcal W$, there are $x_i \in \mathcal X_i$ such that $w =  \sum_i x_i$ and 
\[|w| \geq \sqrt{\frac{1}{C_3l_b}}\sqrt{\sum_i |x_i|^2}.\]
\item Let $x \in \mathcal X$ and $y \in \mathcal Y$ such that $x = Ay$. Then
\[|(1-P)x| \leq C_4\sqrt{\frac{G(l_b)}{l_b}}|y|,\]
where $C_4 = C_1\sqrt{2c_\alpha}\left(\displaystyle\frac{1+\alpha}{1-\alpha}\right).$
\item Let $x \in \mathcal X$ and $x_i \in \mathcal X_i$ such that $x = \sum_i x_i$, then
\[|(1-P)x| \leq C_4 \sqrt{\frac{G(l_b)}{l_b}}\sqrt{\sum_i |x_i|^2}\]
\end{enumerate}
\end{lemma} 
\begin{proof}
For the proof of the first item, note that for $w \in \mathcal W$, there is a (unique) $u \in \mathcal U$ such that $w = Au$ so
\[|w| \geq \sqrt{\frac{1}{C_3l_b}}|u|.\] Now, $\mathcal U \subset \mathcal R$, so there are orthogonal $r_i \in \mathcal R_i$ (because the spaces themselves are orthogonal) such that $\sum_i r_i = u$, hence $|u|^2 = \sum_i |r_i|^2$. Now, $w = Au = \sum_i Ar_i$, so if we set $x_i = Ar_i \in \mathcal X_i$, we obtain $w = \sum_i x_i$. Because $A$ is an isometry when restricted to each $\mathcal R_i$, $|x_i| = |r_i|$ hence we obtain the estimate
\[|w|  \geq \sqrt{\frac{1}{C_3l_b}}|u| = \sqrt{\frac{1}{C_3l_b}}\sqrt{\sum_i |r_i|^2} = \sqrt{\frac{1}{C_3l_b}}\sqrt{\sum_i |x_i|^2}.\]

We now prove the second statement. Because $x \in \mathcal X$, there is a $y \in \mathcal Y (= \mathcal R)$ such that $x = Ay$. Because we chose $n_b$ odd, we can express $y$ as the orthogonal set $\sum_{i \, odd} y_i$ with $y_i = Y_iy$. We find $n_j^i \in \mathcal N_j$ for $j$ even (resp. $\mathcal N_j'$ for $j$ odd) such that $U^\perp y_i = \sum_j n_j^i$ and $|n_j^i| \leq C_1\alpha^{|i-j|}|y_i|$. The idea is that because $\alpha^{|i-j|}$ has exponential decay, the following estimate is like that of an approximation of the identity. 

Recall that $\rho
\geq0$ is tridiagonal (with respect to the spaces $\mathcal R_i$) so $(n^{i_1}_{j_1}, \rho n^{i_2}_{j_2}) = 0 $ if $|j_1 - j_2| \geq 2$. We now finish the proof of the second statement of the lemma with the following estimate using $(1-P)AU=0$ and Lemma \ref{properties of N_i}.2:
\begin{align*}
|(1-P)&x|^2 =|(1-P)AU^\perp y|^2 \leq |A\sum_i U^\perp y_i|^2 \\
&= |A\sum_i\sum_j n_j^i|^2 = (\sum_{i_1}\sum_{j_1} n_{j_1}^{i_1}, \rho\sum_{i_2}\sum_{j_2} n_{j_2}^{i_2}) = \sum_{i_1, j_1, i_2, |j_2 - j_1| \leq 1}( n_{j_1}^{i_1}, \rho n_{j_2}^{i_2}) \\
&\leq \sum_{i_1, j_1, i_2, |j_2 - j_1| \leq 1}( n_{j_1}^{i_1}, \rho n_{j_1}^{i_1})^{1/2}( n_{j_2}^{i_2}, \rho n_{j_2}^{i_2})^{1/2}
= \frac{2G(l_b)}{l_b}\sum_{i_1, j_1, i_2, |j_2 - j_1| \leq 1}|n_{j_1}^{i_1}| |n_{j_2}^{i_2}|\\
&\leq \frac{2C_1^2G(l_b)}{l_b}\sum_{i_1, j_1, i_2, |j_2 - j_1| \leq 1}\alpha^{|i_1-j_1|}|y_{i_1}|\alpha^{|i_2-j_2|}|y_{i_2}|\\
&
\leq \frac{2c_\alpha C_1^2 G(l_b)}{l_b}\sum_{i_1, i_2,j_1}\alpha^{|i_1-j_1|}|y_{i_1}|\alpha^{|i_2-j_1|}|y_{i_2}|\\
&=\frac{2c_\alpha C_1^2 G(l_b)}{l_b}\sum_j\left(\sum_{i}\alpha^{|i-j|}|y_{i}|\right)^2 = \frac{2c_\alpha C_1^2 G(l_b)}{l_b}\|\alpha^{|i|}\ast |y_i|\|_{L^2(\Z)}^2\\
&\leq \frac{2c_\alpha C_1^2 G(l_b)}{l_b}\|\alpha^{|i|}\|_{L^1(\Z)}^2 \||y_i|\|_{L^2(\Z)}^2 = \frac{2c_\alpha C_1^2 G(l_b)}{l_b}\left(\frac{1+\alpha}{1-\alpha}\right)^2 |y|^2,
\end{align*}
where we have used Minkowski's inequality (Theorem 1.2.10 in \cite{Grafakos}) and 
\[\|\alpha^{|i|}\|_{L^1(\Z)} = 1 + 2\alpha\sum_{i \geq 0}\alpha^i = 1 + \frac{2\alpha}{1-\alpha} = \frac{1+\alpha}{1-\alpha}.\]

The last statement of the lemma follows from the second, because we can pick $r_i \in \mathcal R_i$ such that $Ar_i = x_i$ so then $|r_i|^2 = |x_i|^2$. Setting $y = \sum_i r_i$ gives $Ay = x$ and $|y|^2 = \sum_i |r_i|^2$ so that
\[|(1-P)x| \leq C_4\sqrt{\frac{G(l_b)}{l_b}}|y| = C_4\sqrt{\frac{G(l_b)}{l_b}}\sqrt{\sum_i |r_i|^2} = C_4\sqrt{\frac{G(l_b)}{l_b}}\sqrt{\sum_i |x_i|^2}.\]
\end{proof}

\section{Verifying Items 1, 2, and 3.}\label{Verifying}

This section completes the proof of Lemma \ref{main lemma}. Note the similarities between the calculations for the spaces $\tilde{\mathcal X}_i$ in Section \ref{General Approach}.

We first verify Item 1: $\|P^\perp P_{\mathcal V_1}\| \leq \epsilon_1$. Let $v \in \mathcal V_1$. We will find a particular $x \in \mathcal X$ such that $|v-x|$ is small and hence apply the third item of Lemma \ref{estimates for W} to \[|(1-P)v| \leq |(1-P)(v-x)|+|(1-P)x|.\]  

Because $\sum_i \tau_i = \sum_i {\mathcal F}^{0,\kappa}_{\omega(i)}(J)S_1 = S_1$, there is a $\textbf{x} \in \C^d$ such that $v = S_1\textbf{x} = \sum_i \tau_i \textbf{x}$. Because $S_1$ is an isometry, $|v| = |\textbf{x}|$. Then $x_i =\tau_i(1-Z_i) \textbf{x} \in \mathcal X_i$ and 
\[|v - \sum_{i=0}^{n_{win}} x_i|^2 = |\sum_{i=0}^{n_{win}} \tau_iZ_i\textbf{x}|^2 \leq 2\sum_{i=0}^{n_{win}} |\tau_iZ_i\textbf{x}|^2 \leq 2\lambda_{min}(n_{win}+1)|\textbf{x}|^2 = 2\lambda_{min}(n_{win}+1)|v|^2,\]
where the constant $2$ comes from the calculation in equation (\ref{constant 2}), because the ranges of the $\tau_i$ are nonconsecutively orthogonal.

Because $x_i = \tau_i(1-Z_i)\textbf{x}$, we have that 
\[\sum_i|x_i|^2 \leq \sum_i|\tau_i\textbf{x}|^2 = \sum_{i\, odd}|\tau_i\textbf{x}|^2 + \sum_{i\, even}|\tau_i\textbf{x}|^2 = 
|\sum_{i\, odd}{\mathcal F}^{0,\kappa}_{\omega(i)}(J)v|^2 + |\sum_{i\, even}{\mathcal F}^{0,\kappa}_{\omega(i)}(J)v|^2 \leq  2|v|^2,\]
because the functions ${\mathcal F}^{0,\kappa}_{\omega(i)}(t)$ have nonconsecutively disjoint supports so  $\|\sum_{i\, even}{\mathcal F}^{0,\kappa}_{\omega(i)}(J)\|$, $\|\sum_{i\, odd}{\mathcal F}^{0,\kappa}_{\omega(i)}(J)\| \leq 1$.
Then by the third item of Lemma \ref{estimates for W}, we have that 
\begin{align}\nonumber
|(1-P)v|& \leq |(1-P)\sum_i x_i| + |v-\sum_i x_i| \leq C_4 \sqrt{\frac{G(l_b)}{l_b}}\sqrt{\sum_i |x_i|^2} + \sqrt{2\lambda_{min}(n_{win}+1)}|v|\\
&\leq \left(C_4 \sqrt{\frac{2G(l_b)}{l_b}}+\sqrt{2\lambda_{min}(n_{win}+1)}\right)|v|.
\end{align}
We later express the bounds for each item in terms of $L$.

We verify Item 2: $\|(1-P)JP\|\leq \epsilon_1$.
Let $w \in \mathcal W$. There is a $u \in \mathcal U \subset \mathcal Y$ such that $w = Au = AUu$. Write $u$ as the orthogonal sum $\sum_{i\, odd} y_i$, where $y_i = Y_i u$.
Then  \[(1-P)Jw = \sum_i (1-P)JAU y_i.\] 

Now, we want to use the following two facts. First, each $\mathcal X_k$ is approximately an eigenspace for $J$ with eigenvalue $\omega(k)$ so each $\mathcal Y_i = \mathcal X_{[(i-1)l_b, (i+1)l_b)}$ (except $i = n_b$, where it is $\mathcal X_{[(i-1)l_b, n_{win}]}$ having length less than $4l_b$) is also an approximate eigenspace with eigenvalue $\omega(il_b)$. Second, there is exponential decay for $Y_j U Y_i$. Then we note that each $AUy_i \in \mathcal W$ so $(1-P)J(AUy_i) = (1-P)[J - \omega(il_b)](AUy_i)$. This is the place that we use any special property of the projection $P$, because we then remove it by using $\|1-P\| \leq 1$. 
This has a similar feel to the proof of the second item of Lemma \ref{estimates for W}.

In more detail, by Lemma \ref{exponential decay} let $y_j^i = Y_j Uy_i \in \mathcal Y_j$ for $j$ odd so $Uy_i = \sum_{j\, odd}y_j^i$ and $|y_j^i| \leq C_2\alpha^{|i-j|}|y_i|$. Then
\begin{align*}
|(1-&P)Jw|^2 =\left|\sum_i(1-P)[J - \omega(il_b)]AU y_i\right|^2=\left|(1-P)\sum_i \sum_{j \, odd}[J - \omega(il_b)]A y_j^i\right|^2 \\
&\leq \left|\sum_i\sum_{j \, odd} [J - \omega(il_b)]A y_j^i\right|^2 = \sum_{i_1,i_2}\,\sum_{j_1, j_2 \, odd} \left([J - \omega(i_1l_b)]Ay_{j_1}^{i_1}, [J - \omega(i_2l_b)]Ay_{j_2}^{i_2}\right).
\end{align*}
Note $Ay_j^i \in \mathcal X_{[(j-1)l_b,(j+1)l_b)}$ for $j < n_b$. The following work all applies when $j = n_b$ except that the upper limits for intervals and sums are both $n_{win}$. Because $\mathcal X_k$ is a subspace of the range of ${\mathcal F}^{0,\kappa}_{\omega(k)}(J)$, we see that both the $\mathcal X_k$ and the $J(\mathcal X_{k})$ are nonconsecutively orthogonal. Consequently, we continue our calculation as
\begin{align}\label{inner product}
|(1-P)&Jw|^2 \leq \sum_{\substack{i_1,i_2 \\j_1 \, odd}}\sum_{\substack{j_2 \, odd \\ |j_2-j_1| \leq 2}} \left([J - \omega(i_1l_b)]Ay_{j_1}^{i_1}, [J - \omega(i_2l_b)]Ay_{j_2}^{i_2}\right).
\end{align}
We now estimate $|[J - \omega(il_b)]Ay_{j}^{i}|^2$. Because $y_{j}^{i}\in \mathcal Y_j= \mathcal R_{[(j-1)l_b,(j+1)l_b)}$, we write it as an orthogonal sum $\sum_{k=(j-1)l_b}^{(j+1)l_b-1} r_k$ for $r_k \in \mathcal R_k$. We then obtain
\[|[J - \omega(il_b)]Ay_{j}^{i}|^2 = \left|[J - \omega(il_b)]\sum_{k} Ar_k\right|^2= \sum_{k_1, k_2}\left([J - \omega(il_b)]A r_{k_1},[J - \omega(il_b)]A r_{k_2}\right).\]

Note that $J-\omega(k)$ restricted to $\mathcal X_k$ has norm bounded by $\kappa$. Now,  $\mathcal R_k\subset \mathcal Y_j$ for some $j<n_b$ so  
\begin{align*}
|[J - \omega(il_b)]Ar_k|
&\leq |[J - \omega(k)]Ar_k| + |\omega(k) - \omega(jl_b)||Ar_k| + |\omega(il_b)-\omega(jl_b)||Ar_k|  
\\
&\leq 2\kappa(1+(1+|i-j|)l_b)|Ar_k| = 2\kappa(1+(1+|i-j|)l_b)|r_k|.
\end{align*}
Let $C_{i,j}:=2\kappa(1+(2+|i-j|)l_b)$ to account for the cases $j < n_b$ and $j = n_b$ because $\mathcal Y_j$ has length at most $4l_b$.
So, keeping in mind that the range of $k$, (we can define $r_k = 0$ outside this range) we see that

\begin{align*}
|[J -& \omega(il_b)]Ay_{j}^{i}|^2 = \sum_{k_1}\sum_{|k_1-k_2|\leq 1}\left([J - \omega(il_b)]A r_{k_1},[J - \omega(il_b)]A r_{k_2}\right)\\
&\leq C_{i,j}^2\sum_{k_1}\sum_{|k_1-k_2|\leq 1}|r_{k_1}||r_{k_2}| =  C_{i,j}^2\sum_{\sigma =-1}^1\sum_{k_1}|r_{k}||r_{k+\sigma}|\\
&\leq C_{i,j}^2\sum_{\sigma =-1}^1\left(\sum_{k}|r_{k}|^2\right)^{1/2}\left(\sum_{k}|r_{k+\sigma}|^2\right)^{1/2} \\
&\leq 3C_{i,j}^2\sum_{k}|r_{k}|^2 = 3C_{i,j}^2|y_j^i|^2.
\end{align*}
Using $2\sqrt{3}\kappa(1+(2+|i-j|)l_b)) \leq 2\sqrt{3}(\kappa l_b)(3+|i-j|) =: K(3+|i-j|)$ and $C_\alpha = \max_x (|x|+3)\alpha^{|x|/2}$, we insert our above calculations into Equation (\ref{inner product}) to get
\begin{align*}\label{inner product}
|(1&-P)Jw|^2 \leq K^2\sum_{\substack{i_1,i_2 \\j_1 \, odd}}\sum_{\substack{j_2 \, odd \\ |j_2-j_1| \leq 2}} (3+|i_1-j_1|)(3+|i_2-j_2|)|y_{j_1}^{i_1}||y_{j_2}^{i_2}| \\
&\leq (C_2K)^2\sum_{\substack{i_1,i_2 \\j_1 \, odd}}\sum_{\substack{j_2 \, odd \\ |j_2-j_1| \leq 2}} (3+|i_1-j_1|)(3+|i_2-j_2|)|y_{i_1}|\alpha^{|i_1-j_1|}|y_{i_2}|\alpha^{|i_2-j_2|}\\
&\leq c_{\alpha}(C_\alpha C_2K)^2\sum_{j \, odd}\,\sum_{i_1,i_2}  |y_{i_1}|\alpha^{|i_1-j|/2}|y_{i_2}|\alpha^{|i_2-j|/2} \\
&= c_{\alpha}(C_\alpha C_2K)^2\sum_{j \, odd}\,\left(\sum_{i}  |y_{i}|\alpha^{|i-j|/2}\right)^2 \leq c_{\alpha}(C_\alpha C_2K)^2\left(\frac{2+\alpha}{2-\alpha}\right)^2\sum_k|y_k|^2 \\
&= c_{\alpha}K^2C_\alpha^2C_2^2\left(\frac{2+\alpha}{2-\alpha}\right)^2|u|^2 =: (K_\alpha\kappa l_b |u|)^2,
\end{align*}
where $K_\alpha$ is a constant.
Because $u \in \mathcal U$ with $w = Au$,  Lemma \ref{(u, rho u)} shows that
\[|(1-P)Jw| \leq \kappa l_b  K_\alpha\sqrt{C_3l_b} |w| = Const.\kappa l_b^{3/2} |w|.\]

Now, we address the third item. For $w \in \mathcal W$, we will bound $P_{\mathcal V_L}(w)$. Now, by Lemma \ref{estimates for W}, we can write $w = \sum x_i$, $x_i \in \mathcal X_i$ such that 
\begin{align}\label{x-wIneq}
\sqrt{\sum_i |x_i|^2} \leq \sqrt{C_3l_b}|w|.\end{align}
We will bound each $P_{\mathcal V_L}(x_i)$ using the Lieb-Robinson estimates. Let $\hat B$ be a position operator on $\mathcal B$ being $jI$ on $\mathcal V_j$. Then $J$ is tridiagonal with respect to these blocks so it satisfies the conditions of Corollary \ref{operator Lieb-Robinson} with $\Delta = 2$. So,
\[
\|P_{\mathcal V_L}{\mathcal F}^{0,\kappa}_{\omega(i)}(J)P_{\mathcal V_1}\|\leq \int_{|k| \geq \frac{L-1}{2e^2}} |\hat{\mathcal F}_0^{0, \kappa}(k)|dk  + \|\hat{\mathcal F}_0^{0, \kappa}\|_{L^1(\R)}e^{-\frac{L-1}{2}}.
\]
By equations (\ref{invariant under scaling}) and (\ref{defining F(L)}) and the definition of $\kappa$,
\[\|P_{\mathcal V_L}{\mathcal F}^{0,\kappa}_{\omega(i)}(J)P_{\mathcal V_1}\| \leq \int_{|k| \geq \frac{L-1}{e^2n_{win}}} |\hat{\mathcal F}_0^{0, 1}(k)|dk + \|\hat{\mathcal F}_0^{0, 1}\|_{L^1(\R)}e^{-\frac{L-1}{2}} = S(L),
\]
as defined in Section \ref{smooth partitions}. Now, since $x_i \in \mathcal X_i$, we have that there are $\textbf{x}_i \in \C^d$ such that $x_i =  {\mathcal F}^{0,\kappa}_{\omega(i)}(J)S_1(1-Z_i)\textbf{x}_i$. We can pick $\textbf{x}_i$ in the kernel of $Z_i$ so that $x_i = {\mathcal F}^{0,\kappa}_{\omega(i)}(J)S_1\textbf{x}_i$, with $S_1\textbf{x}_i \in \mathcal V_1$ and
$|x_i|^2 = (\tau_i^\ast \tau_i\textbf{x}_i, \textbf{x}_i) \geq \lambda_{min}|\textbf{x}_i|^2 = \lambda_{min}|S_1\textbf{x}_i|^2$. So,
\[|P_{\mathcal V_L}x_i| = |P_{\mathcal V_L}{\mathcal F}^{0,\kappa}_{\omega(i)}(J)S_1\textbf{x}_i| \leq S(L)|S_1\textbf{x}_i| \leq \frac{S(L)}{\lambda_{min}^{1/2}}|x_i|.\]

By Equation (\ref{x-wIneq}), we have
\[|P_{\mathcal V_L}(w)| \leq \sqrt{n_{win}+1}\left(\sum_{i=0}^{n_{win} }|P_{\mathcal V_{L}}(x_i)|^2\right)^{1/2} \leq S(L)\sqrt{\frac{C_3(n_{win}+1)l_b}{\lambda_{min}}}|w|.\]

So that we have gotten our estimates, we write $n_b \sim L^{\beta_0}$, $n_{win} \sim L^{\beta_1}/F(L)$, $\lambda_{min}\sim 1/(L^{\beta_2}(n_{win}+1))\sim F(L)/L^{\beta_1+\beta_2}$, for $\beta_0, \beta_1, \beta_2 > 0,$ and in the definition of $n_{win}$ we had $\beta_1 \leq 1$. Recall that $l_b \sim n_{win}/n_b\sim L^{\beta_1-\beta_0}/F(L)$ and $\kappa \sim 2F(L)/L^{\beta_1}$.

We then get 
\[\epsilon_3 \lesssim \sqrt{G\left(\frac{L^{\beta_1-\beta_0}}{F(L)}\right)F(L)}\;L^{\beta_0/2 - \beta_1/2} + L^{-\beta_2/2},\]
\[\epsilon_4 \lesssim \frac{F(L)}{L^{\beta_1}}\left(\frac{L^{\beta_1-\beta_0}}{F(L)}\right)^{3/2}= \frac{L^{\beta_1/2 - 3\beta_0/2}}{F(L)^{1/2}},\]
and
\[\epsilon_5 \lesssim S(L)\sqrt{\frac{(n_{win}+1)l_b}{\lambda_{min}}}\sim \frac{S(L)}{F(L)^{3/2}}L^{3\beta_1/2+\beta_2/2-\beta_0/2}.\]
So, if $\epsilon_3$ and $\epsilon_4$ have similar rates, because we are assuming that $G(l)$ increases slower than any power of $l$ and $F(L)$ grow slower than any power of $L$, we get $\beta_0 - \beta_1 = \beta_1 - 3\beta_0 = -\beta_2$, hence $\beta_0 = \beta_1/2$. This gives a rate of $L^{-\beta_1/4}$. However, since $\beta_1 \leq 1$, the best that this gives us is $L^{-1/4}$, which is why we pick $\beta_1=1$. We pick $\beta_2 = 1/2$ and $\beta_0 = 1/2$. We get $\gamma_2 = 1/4$.

This ends the proof of Lemma \ref{main lemma}.

\section{Almost commuting Hermitian and Normal matrices}\label{A and N}

Section 1 of
\cite{HastingsLoring} formulates various almost commuting - nearly commuting problems under certain ``geometric'' restrictions.
These examples include the two that  are directly related to the primary reformulations of Lin's theorem: 
\begin{enumerate}
\item The geometry of a square, which is $\{(x_1,x_2)\in \R^2:|x_1|, |x_2| \leq 1\}$, giving rise to almost commuting Hermitian $H_1, H_2$ with $\|H_1\|, \|H_2\| \leq 1$, and 
\item The geometry of the disk, which is $\{z \in \C : |z|\leq 1\}$, giving rise to an almost normal $N$ with $\|N\| \leq 1$.
\end{enumerate}

Some other examples (which we will discuss below) are that of almost commuting Hermitian  and unitary matrices (which is related to the geometry of a cylinder/anulus) and two almost commuting unitaries (the geometry of the torus). 
The latter does not always have nearby commuting matrices, but \cite{Tobias} showed that if both unitaries have a spectral gap then one obtains nearby commuting unitaries. The proof given there involves using a matrix logarithm to reduce the unitary matrix with a spectral gap into a Hermitian matrix.

Both of the geometries discussed above (and more general geometries) are addressed by Enders and Shulman in \cite{Dimension} providing a dimensional and cohomology criterion for liftings in the spirit of Voiculescu's comment in \cite{Voiculescu}.
In this section we discuss an approach to a few of these types of results following the ideas presented in previous sections. 

We make some remarks on the continuity of $f:\R\to \C$ or $f:\C\to\C$ as a function on normal matrices. A continuous function $f$ of a complex variable is operator continuous (see Proposition II.2.3.3 of \cite{Blackadar}). It is important to note that the operator modulus of continuity of $f$ is bounded below by the modulus of continuity of $f$ as a function of a real/complex variable. The proof in \cite{Blackadar} shows that the operator modulus of continuity can be chosen independent of the $C^\ast$-algebra. See \cite{OC} for a more detailed analysis of the operator modulus of continuity. As seen in a remark before Lemma 5.2 of \cite{OC}, which references proofs given in Section 10 of \cite{OHZ}, if $f$ is a continuous function, then its operator modulus of continuity is equivalent to its ``commutator modulus of continuity''. That is,
\[\Omega_{f,\mathfrak F}(\delta) =\sup\{\|f(A)-f(B)\|: A,B \mbox{ normal}; \sigma(A),\sigma(B)\subset \mathfrak F; \|A-B\|\leq \delta\}\]
is equivalent to
\[\Omega_{f,\mathfrak F}^\flat(\delta) =\{\|[f(A),R]\|: A \mbox{ normal}, \sigma(A)\subset \mathfrak F, R \mbox{ self-adjoint}, \|R\|\leq 1, \|[A,R]\|\leq \delta\}.\]

For a summary of Operator Lipschitz functions and the norm $\|-\|_{\operatorname{OL}(\R)}$ see \cite{OL}. 
Note that  $\|f\|_{\operatorname{OL}(\R)} \leq C_f$, which can be seen by using the Fourier representation of $f$ and $\|e^{ikx}\|_{\operatorname{OL}(\R)} = |k|$ from \cite{OL}.
Although $f$ being a Lipschitz function of a real/complex variable does not ensure that $f$ is Operator Lipschitz, its operator Lipschitz constant is at least its Lipschitz constant.

We present a na\"ive non-proof that given two commuting self-adjoint contractions $A_1, A_2$ that almost commute with the self-adjoint contraction $B$, there exist commuting self-adjoint $A_1', A_2', B'$ that are close to $A_1, A_2, B$, respectively. 
The idea is that because $[A_1,A_2] = 0$, we can simultaneously diagonalize these matrices and hence there exists a self-adjoint matrix $A$ whose eigenspaces are subspaces of the eigenspaces of $A_1$ and the eigenspaces of $A_2$. Therefore, we can write $A_1 = f_1(A), A_2 = f_2(A)$ for some continuous functions $f_1, f_2 : \R \to \R$.

We hope that $[A,B]$ is small so that there are commuting $A', B'$  near $A$ and $B$, respectively and that if we define $f_1(A') = A_1', f_2(A') = A_2'$ then $A_1'$ is near $A_1$ and $A_2'$ is near $A_2$. It is a given that $A_1', A_2', B'$ all commute, but our other ``hopes'' are not guaranteed to be realizable. The first issue is the commutator of $[A,B]$ might not be small. The second issue is that we might not have control of how close $A_i' = f_i(A')$ is to $A_i$ if the functions $f_i$ depend on the original matrix $A$.

Exploring an example more in line with what we will be doing, suppose that the spectrum of a normal matrix $N$ is a subset of the graph $\{(x,f(x)): x\in[-1,1]\}$, where we identify $\R^2$ with $\C$. If $f$ is an operator Lipschitz function, then set $A_1 = (N + N^\ast)/2$, $A_2 = (N-N^\ast)/2i$. There are commuting Hermitian $A_1', B'$. Because $A_2=f(A_1)$, setting $A_2' = f(A_1')$ makes $N' = A_1' + iA_2'$ normal and close to $N$ because
$\|A_2-A_2'\| = \|f(A_1) - f(A_1')\|$.  We then have that $N, B$ are nearly commuting. 

What we mean by ``a  nice function'' is that either we have one fixed operator Lipschitz function $f$ or we have a family of  functions $\mathfrak F$ with bounded operator Lipschitz norm from which we pick $f$.

To illustrate how the above naive construction can fail, consider $N$ to have spectrum in $[0,1]^2$.  We can perturb $N$ (very) slightly so the spectrum is in $\{(\frac{k+\epsilon j}{n},\frac{j}{n}): -n \leq k,j \leq n\}$ for some small $\epsilon > 0$. 
Then an interpolating function $f$ will display behavior like a saw-tooth function with slope at least $\epsilon$ and hence has growing operator modulus of continuity. 
The example of Davidson in \cite{Davidson} is of a weighted shift operator in $M_n(\C)$ (where $n-1$ is a square) that is almost normal and almost commutes with a diagonal Hermitian matrix.
Following his use of \cite{Berg}, one sees that we get a normal matrix that is unitarily equivalent to $\bigoplus_{k=0}^{M} c_k S_{d_k}$, where $S_m$ is a cyclic permutation matrix of length $m$ and $c_k=\frac{M-k}{M}$. This normal matrix has spectrum in the union of concentric circles of radii $c_k$ centered at the origin. The limit points of its spectrum as $n \to \infty$ is the entire unit disk. \cite{Davidson} proved that this normal matrix and a Hermitian matrix do not have nearby commuting matrices. 

A way that the above construction works is if the spectrum of $N$ lies within a simple curve $\Gamma$ with distinct endpoints  such that there is a closed (and bounded) interval $I = [-a,a]$ and operator Lipschitz functions $f:I\to\Gamma\subset \C$ and $g:\C\to I$ such that $f$ is a bijection and $g\circ f = \operatorname{id}_I$. Then if $N$ almost commutes with $B$, $A=g(N)$ almost commutes with $B$ as discussed above. Then nearby $A,B$ there are commuting Hermitian $A', B'$. Scaling $A'$ by a number slightly less than $1$ if necessary as in \cite{HastingsLoring}, we can assume that  the spectrum of $A$ is in $I$. Setting $N' = f(A')$, we see that $N'$ commutes with $B'$ and
\[\|N' -N\| = \|f(A')-f(A)\|\]
is small when  $\|[N,B]\|$ is small. If $\Gamma$ is the disjoint union of such curves, then one can proceed in a similar way by forming block matrices gotten by collecting the spectrum lying in different curves.

We will explore how a type of one-dimension spectrum of a normal matrix $N$ can guarantee that if $N$ and a Hermitian matrix $A$ almost commute, then there are nearby commuting  normal and Hermitian matrices.
We first approach the situation of an almost commuting Hermitian contraction $A$ and unitary $U$. This provides an alternative approach to the proof in \cite{HastingsLoring} regarding an ``almost representation of the disk''.
\begin{prop}
Let $A$ be a Hermitian contraction and $U$ unitary. Then if $A, U$ are almost commuting, they are nearly commuting with the same asymptotic rates as in Theorem \ref{maintheorem}.
\end{prop}
\begin{proof}
Let $\delta = \|[A,U]\|$. Picking  $\Delta = \delta^{\gamma_0}$, we apply Corollary \ref{finite range normal} to obtain a finite range (of distance $\Delta$ with respect to the normal $U$) $H$ such that
$\|A-H\|\leq Const.\delta^{1-\gamma_0}$ and $\|[H,N]\|\leq Const.\delta^{\gamma_0}$.

Recall that in Section \ref{reformulations}, where we had $B$ self-adjoint, we partitioned its spectrum into intervals much larger than $\Delta$. We used the tridiagonal nature of $H$ by grouping together eigenvalues of $B$ so as to apply Lemma \ref{modified lemma} to a ``pinching'' of $H$ by certain projections to obtain the subspaces $\mathcal W$.

We do the same here: we partition the spectrum of $U$ into arcs $I_i$ by forming $n_{win}\sim \Delta^{-\gamma_1}/2\pi$ many intervals, where $0 < \gamma_1 < 1$. 
This gives the subspaces $\mathcal B_i$ on which $H$ restricted is block tridiagonal with many blocks and on which $U$ takes on potentially many eigenvalues, but all in the interval $I_i$.
For $H$ projected onto $\mathcal B_i$, we then obtain the almost invariant subspaces $\mathcal W_i$ satisfying the properties specified in Lemma \ref{modified lemma} including $\tilde{\mathcal B}_i = \mathcal W_i^\perp \oplus \mathcal W_{i+1} \subset \mathcal B_i \oplus \mathcal B_{i+1}$. We then have that $\tilde{\mathcal B}_i$ is almost invariant for $H$. 
This construction is essentially identical, with the key distinction that although there are ``first'' and ``last'' subspaces, we index the spaces $\mathcal B_i$ cyclically and must (at least for our construction to work) define $\tilde{\mathcal B}_{n_{cut}} = \mathcal W_{n_{cut}}^\perp \oplus \mathcal W_{n_{cut}+1} = \mathcal W_{n_{cut}}^\perp \oplus \mathcal W_{1} = \tilde{\mathcal B}_{0}.$

Let $U'$ be a multiple of the identity on $\tilde{\mathcal B}_{i}$, with the multiple being the center of the arc $I_i$. We project $H$ onto the $\tilde{\mathcal B}_{i}$ to obtain $H'$ and we conclude.
\end{proof}

Note that we then obtain the following form of Osborne's result in \cite{Tobias} that only requires one matrix to have a spectral gap. We only need to transform one unitary into a Hermitian matrix, but we use a modification of a fractional linear transformation instead for simplicity.
\begin{prop}\label{Gap}
Fix $\theta \in (0,\pi)$. Then for unitaries $U, V$ where $V$ has a spectral gap of angular radius $\theta$, there are nearby commuting unitaries $U',V'$. If $\epsilon = \epsilon(\delta)$ derived from Lin's theorem for a unitary and self-adjoint pair, then \[\|U-U'\|\leq \epsilon\left(\frac{\|[U,V]\|}{\sin\theta}\right), \;\;\|V-V'\|\leq \frac{2\sin\theta}{1-\cos\theta}\epsilon\left(\frac{\|[U,V]\|}{\sin\theta}\right).\] 
\end{prop}
\begin{proof}
Just as in \cite{Tobias}, we can multiply $V$ by a phase so that the spectral gap is centered at $1$. Let \[c_\theta = \frac{1-\cos\theta}{\sin\theta}.\] Set $f(x) = \frac{x-c_\theta i}{x+c_\theta i}$ and $g(z) = ic_\theta\frac{1+z}{1-z}$. A simple calculation (or a simple comparison to fractional linear transformations in \cite{Stein}) shows that $f$ maps $\R$ to the unit circle and $g$ maps the unit circle to $\R$ as its inverse. For $z = e^{i\varphi}$, one can see that
\[g(e^{i\varphi}) = -c_\theta\frac{\sin\varphi}{1-\cos\varphi}\]
using $(1+e^{i\varphi})(1-e^{-i\varphi})=2i\sin\varphi$ and $|1-e^{i\varphi}|^2=2(1-\cos\varphi)$.

We define the self-adjoint matrix $W = g(V)$.
We included the scaling factor $c_\theta$ in the definitions of $f$ and $g$ so that $\|W\|\leq1$.  Now, using $\frac{z+a}{z-a} = 1 + \frac{2a}{z-a}$ we see that
\begin{align*}
\|[U, W]\| &= 2c_\theta\|[U, (V-1)^{-1}]\|= 2c_\theta\|(V-1)^{-1}(\,(V-1)U-U(V-1)\,)(V-1)^{-1} \| \\
&\leq2c_\theta\|(V-1)^{-1}\|^2\|[U,V-1]\|  \leq \frac{c_\theta}{1-\cos\theta}\|[U,V]\|.
\end{align*}
We then can find nearby commuting $U', W'$ where $U'$ is unitary and  $W'$ is self-adjoint. Then $V' = f(W')$ is unitary and we have
\begin{align*}
\|V-V'\|&\leq \|(W-ic_\theta)(W+ic_\theta)^{-1} - (W'-ic_\theta)(W'+ic_\theta)^{-1}\| \\
&= 2c_\theta\|(W+ic_\theta)^{-1}-(W'+ic_\theta)^{-1}\|\leq \frac2{c_\theta}\|W-W'\|
\end{align*}
by the resolvent identity calculation in Example 1 of \cite{OL} or just by arguments similar to what we did above. The result then follows.
\end{proof}

Now we can attempt to do the same sort of construction for a normal matrix $N$ whose spectrum belongs to some set $\Gamma$.  We define the following condition on a $\Gamma \subset \C$, but note that it is not necessarily optimal.

\vspace{0.1in}

\noindent ($\ast$): 
$\Gamma$ will be the union of $k_m$ many simple curves $\Gamma^k$ which include their two endpoints. Let $\{z_1, \dots, z_n\}$ be the endpoints of the $\Gamma^k$. The curves $\Gamma^k$ are not permitted to intersect except at the endpoints. A  closed closed curve can be included in this framework by picking two points on the curve and break the curve into two simple curves sharing two distinct endpoints.

Let $\Delta_0, C_s>0, \gamma_1\in(0,1)$ and $0<\Delta< \ell_s/2$, where $\ell_s=\ell_s(\Delta)\sim \Delta^{\gamma_s}$. The constants $C^k,\ell^k$ have the same restrictions as $C_s, \ell_s$, respectively. 
Let $\Gamma^{k,\Delta} = \Gamma^k\setminus\bigcup_s B_{\ell_s}(z_s),\tilde\Gamma^{k,\Delta,s} = (B_{l_s}(z_s)\cap \Gamma^k) \setminus B_{l_s/2}(z_s)$.

We require the following properties to hold for $\Delta \leq \Delta_0$:
\begin{enumerate}
\item If $z_s$ is not an endpoint of $\Gamma^k$, then $\dist(z_s,\Gamma^k)> \ell_s$.
\item If $z_s$ is an endpoint of $\Gamma^k$, then $\Gamma$ only intersects $\partial B_{\ell_s}(z_s)$ in a single point, which is a point in $\Gamma^k$.
\item $\dist( \tilde\Gamma^{k,\Delta}, \tilde\Gamma^{k',\Delta}\cup \Gamma^{k',\Delta})\geq (1-\delta_{k,k'})\Delta$.
\item Suppose that $\Gamma^k$ has endpoints $z_{s_1}$ and $z_{s_2}$. There is a partition of $\Gamma^k$ into disjoint consecutive arcs $A^{k,\Delta}_1, \dots, A^{k,\Delta}_{r^{k,\Delta}}\subset \Gamma^{k,\Delta}$ satisfying the following properties:
\begin{enumerate}
\item $\frac{1}{C^k}\ell^k\leq \diam( A^{k,\Delta}_i)\leq C^k\ell^k$.
\item $\dist( A^{k,\Delta}_i, A^{k,\Delta}_{i'})\geq (1-\delta_{i,i'})\Delta$.
\item Each $A^{k,\Delta}_i$ can be broken up into consecutive subarcs $A^{k,\Delta}_{i,j}$ each satisfying\\
$\diam(A^{k,\Delta}_{i,j})\leq C^k\Delta$ and $\dist( A^{k,\Delta}_{i,j}, A^{k,\Delta}_{i,j'})\geq (1-\delta_{j,j'})\Delta.$
\item If $z_s$ is an endpoint of $\Gamma^k$ then there is a unique value $i_s^k$ of $i$ and $j_s^k$ of $j$ so that $A^{k,\Delta,s}_{i,j}$ is within a distance of $\Delta$ from $ B_{\ell_s}(z_s)$. One of the endpoints of this subarc will then be on $\partial B_{\ell_s}(z_s)$. \end{enumerate}
\end{enumerate}

How small that $\Delta$ needs to be will depend on how close the $\Gamma^k$ are near and away from $z_j$. It is necessarily the case that if there are many $\Gamma^k$ with endpoint $z_s$, then $l/\Delta$ must be large. This is illustrated by Figure 4 below, in particular how visually near the intersection point one can compare the size of the black ``dot'' formed by the intersection of line segments depicted with positive width.
\begin{figure}[htp]
    \centering
    \includegraphics[width=7cm]{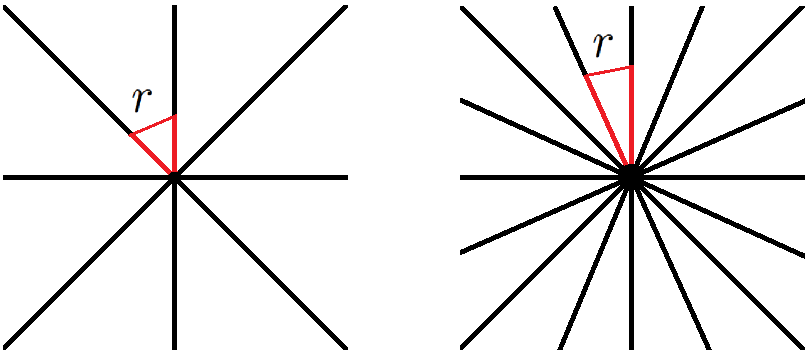}
    \caption{For a given $r > 0$, compare how far one has to move from the point of intersection so that each line is at least a distance of $r$ from the other lines.}
\end{figure}

\begin{figure}[htp]
    \centering
    \includegraphics[width=7cm]{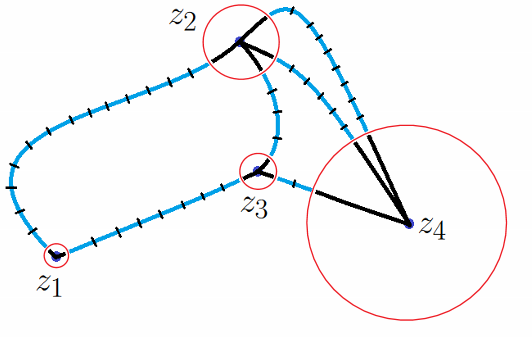}
    \caption{  $\Gamma^{k,\Delta}$ is depicted in blue, $\partial B_{l_j}(z_j)$ in red, and $\Gamma^{k,\Delta}\cap B_{l_s}(z_s)$ in black. The tick marks on the blue curves separate the arcs $A_i^{k,\Delta}$ outside of the disks.
    Although not all end points are shared by multiple curves $\Gamma^k$ in general, it is this case in the illustration above.}
\end{figure}

We then have the following:
\begin{prop}
Let $N$ be normal with spectrum in $\Gamma$ satisfying $(\ast)$ and $A$ self-adjoint. Then if $A, N$ are almost commuting, they are nearly commuting.
\end{prop}
\begin{proof}
Let $\delta = \|[A,N]\|$. We apply the argument earlier to obtain $H$ finite range with respect to $N$ of distance $\Delta \sim \delta^{\gamma_0}$ such that $\|A-H\|\leq Const.\delta^{1-\gamma_0}, \|[H,N]\| \leq Const. \delta$. Recall $(\ast)$.

Let $\mathcal B^{k,\Delta}_i$, $\mathcal V^{k,\Delta}_{i,j}$ be the range of the spectral projection of $N$ onto $A^{k,\Delta}_i$, $A^{k,\Delta}_{i,j}$, respectively. We then obtain $\mathcal W^{k,\Delta}_i$ which are subspaces of the $\mathcal B^{k,\Delta}_i$ that are almost invariant under $P_{\mathcal B^{k,\Delta}_i}HP_{\mathcal B^{k,\Delta}_i}$ such that $\mathcal V^{k,\Delta,s}_{i,1}\subset \mathcal W^{k,\Delta}_i\perp {\mathcal V}^{k,\Delta,s}_{i,r^{k,\Delta,s}_i}$. Let $\mathcal W^{k,\Delta,\perp}_i=\mathcal B^{k,\Delta}_i\ominus\mathcal W^{k,\Delta}_i$.

We now need to deal with the subspace $\mathcal B^{\Delta,s}$ projected onto by $E_{B_{l_s}(z_s)}(N)$. The issue is that because $z_s$ can be the endpoint of multiple curves $\Gamma^k$, we have to use a slightly different approach to form almost invariant subspaces of $\mathcal B^{\Delta,s}$ so that everything works out.
This is the only part of the proof where one uses a ``two dimensional'' grouping of eigenvalues of $N$ (using ``polar coordinates'').

Fix an endpoint $z_s$. Let $S_1=B_{\ell_s/2}(z_s)$, $S_i^k=\Gamma^k\cap B_{(i-1)\Delta+\ell_s/2}(z_s)\setminus B_{(i-2)\Delta+\ell_s/2}(z_s)$, $S_{r^{\Delta,s}}^k=\Gamma^k\cap     B_{l_s}(z_s)\setminus B_{(r^{\Delta,s}-1)\Delta+\ell_s/2}(z_s)$, where $\ell_s/2-1< r^{\Delta,s}\Delta\leq \ell_s/2$. 
What we do is break $\mathcal B^{\Delta,s}$ into orthogonal subspaces $\mathcal V^{\Delta,s}_1,\mathcal V^{k,\Delta,s}_2, \dots, \mathcal V^{k\Delta,s}_{r^{\Delta,s}}$ by letting $E_{S_1}(N)$ 
project onto $\mathcal V^{\Delta,s}_1$ and
$E_{S_i^k}(N)$ project onto $\mathcal V^{k,\Delta,s}_i$ for $1< i\leq r^{\Delta,s}$. Consider  $H^{\Delta,s}=P_{\mathcal B^{\Delta,s}}H P_{\mathcal B^{\Delta,s}}$ as a block matrix with respect to the subspaces $\mathcal V^{\Delta,s}_1,\mathcal V^{k,\Delta,s}_2, \dots, \mathcal V^{k\Delta,s}_{r^{\Delta,s}}$.

For each $k$, we restrict $H$ to $\mathcal B^{k,\Delta,s}=\mathcal V^{\Delta,s}_1\oplus\bigoplus_{i=2}^{r^{\Delta,s}}\mathcal V^{k,\Delta,s}_i$ 
to obtain $H^{k,\Delta,s}$. $H^{k,\Delta,s}$ is tridiagonal with respect to the subspaces $\mathcal V^{\Delta,s}_1,\mathcal V^{k,\Delta,s}_2, \dots, \mathcal V^{k\Delta,s}_{r^{\Delta,s}}$ which correspond to a set covering $\Gamma^k\cup S_1$ in $B_{\ell_s}(z_s)$. Then there is a subspace $\mathcal W^{\Delta,s}_k\subset \mathcal B^{k,\Delta,s}$ that is almost invariant under $H^{k,\Delta,s}$ such that
$\mathcal V^{k,\Delta,s}_{r^{\Delta,s}}\subset \mathcal W^{\Delta,s}_k\perp \mathcal V^{\Delta,s}_1$.

\begin{figure}[htp] 
    \centering
    \includegraphics[width=7cm]{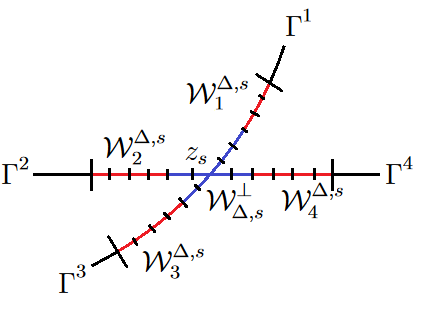}
    \caption{The point $z_s$ is an endpoint for $\Gamma^1, \Gamma^2, \Gamma^3,$ and $\Gamma^4$. The small tick marks are the intersections of circles centered at $z_s$ of radii $i\Delta$. The almost invariant subspaces
    $\mathcal W^{\Delta,s}_k$ are illustrated in red and $\mathcal W^\perp_{\Delta,s}$ in blue. Although these subspaces do not necessarily lie entirely within the span of the illustrated eigenspaces of $N$, $\mathcal W^\perp_{\Delta,s}$ contains all eigenspaces of $N$ within distance $l_s/2$ of $z_s$ and  $\mathcal W^{\Delta,s}_k$ contains the eigenspaces depicted by part of the (red) spectrum near the large tick marks.}
\end{figure}
Now, $P_{\mathcal W^{\Delta,s}_k}\leq E_{\tilde{\Gamma}^{k,\Delta,s}}(N)$ so not only are the subspaces $\mathcal W^{\Delta,s}_k$ orthogonal for different values of $k$, but $H^{\Delta,s}$ maps $\mathcal W^{\Delta,s}_k$ into $\mathcal W^{\Delta,s}_k\oplus\mathcal V^{\Delta,s}_1$ due to the assumption on the distance between the $\tilde{\Gamma}^{k,\Delta,s}$ for different values of $k$. Hence, $H^{\Delta,s}$ is equal to $H^{k,\Delta,s}$ on $\mathcal W^{\Delta,s}_k$ and $H^{\Delta,s}(\mathcal W^{\Delta,s}_k)$ is perpendicular to $\mathcal W^{\Delta,s}_{k'}$ for $k' \neq k$ since $ \mathcal W^{\Delta,s}_{k'}\perp\mathcal V^{\Delta,s}_1$ by construction.
Then we find that $\mathcal W^{\Delta, s}_k$ are almost invariant under $H^{\Delta,s}$ because 
\begin{align*}
(1-P_{\mathcal W^{\Delta,s}_k})H^{\Delta,s}P_{\mathcal W^{\Delta,s}_k}=P_{\mathcal B^{k,\Delta,s} \ominus \mathcal W^{\Delta, s}_k}H^{k,\Delta,s}P_{\mathcal W^{\Delta,s}_k}
\end{align*}
has small norm.

Let $\mathcal W^\perp_{\Delta,s}=\mathcal B^{\Delta,s}\ominus\bigoplus_k\mathcal W^{\Delta,s}_k$. Because $H^{s,\Delta}$ is Hermitian, we have  broken $\mathcal B^{\Delta,s}$ into almost invariant subspaces $\mathcal W^{\Delta,s}_k,\mathcal W^\perp_{\Delta,s}$.

We now form the new basis of subspaces $\tilde{\mathcal B}$. This is a simple process, but its explanation is complicated by the fact that because of parity issues, we cannot just say ``join $\mathcal W$'s with $\mathcal W^\perp$'s'' as before. 
We let the subspaces $\tilde{\mathcal B}$ include three types of subspaces. The first subspaces that we include are $\mathcal W^\perp_{\Delta,s}$. The other ``$\mathcal W$'' subspaces contain a ``${\mathcal V}$'' subspace which has to be matched with its neighboring ``${\mathcal V}$'' subspace. One of the subspaces $\mathcal W^{k,\Delta}_i, \mathcal W^{k,\Delta,\perp}_i$ contains that subspace $\mathcal V^{k,\Delta}_{i_s^k,j_s^k}$, call this subspace $\mathcal W_{i_s^k, j_s^k}$. We then include the subspace $\mathcal W_{i_s^k, j_s^k}\oplus \mathcal W^{\Delta, s}_k$.  We also include the direct sums of the remaining consecutive subspaces $\mathcal W^{k,\Delta,\perp}_{i}\oplus\mathcal W^{k,\Delta}_{i+1}$.
    
Having defined the subspaces $\tilde{\mathcal B}$, we proceed to defining the nearby commuting matrices. We define $N'$ to be the normal matrix that is block identity picking an eigenvalue from the arcs (or union of arcs) that we used to make the spaces $\tilde{\mathcal B}$. We then have that $N'$ and $H'=\sum_{\tilde{\mathcal B}}P_{\tilde{\mathcal B}}HP_{\tilde{\mathcal B}}$ commute and these matrices are within a distance of a power of $\Delta$ from $N$ and $H$, respectively.
\end{proof}

\newpage

\textbf{ACKNOWLEDGEMENTS}. The author would like to thank Eric A. Carlen for introducing the issue of \cite{Hastings} to the author and for helpful discussions about it, Matthew Hastings for some corrections of and clarifications concerning \cite{Hastings}, and Ilya Kachkovskiy for providing feedback concerning the  characterization of the constructiveness of \cite{KS} in the first version of this paper.

This research was partially supported by NSF grants DMS-2055282 and DMS-1764254.


\begin{thebibliography}{10}

\bibitem{OC} A.B. Aleksandrova, V.V. Peller. {\it Estimates of operator moduli of continuity.} Journal of Functional Analysis 261 (2011) 2741–2796.
 
\bibitem{OHZ} A. B. Aleksandrov and V. V. Peller. {\it Operator H\"older-Zygmund Functions.} \href{https://arxiv.org/abs/0907.3049}{	arXiv:0907.3049}

\bibitem{OL} A. B. Aleksandrov and V. Peller. {\it Operator Lipschitz functions.} 2016 Russ. Math. Surv.71 605.

\bibitem{Subnormal} J. Bastian and K. Harrison. {\it Subnormal Weighted Shifts and Asymptotic Properties of Normal Operators} Proc. Amer. Math. Soc. 42 (1974), 475-479. \url{https://doi.org/10.1090/S0002-9939-1974-0380491-X }

\bibitem{Benzi} M. Benzi and G. Golub. {\it Bounds for the Entries of Matrix Functions with Applications to Preconditioning.} BIT Numerical Mathematics (1999) 39: 417. \url{https://doi.org/10.1023/A:1022362401426}

\bibitem{berg-olsen} I. Berg and C. Olsen. {\it A note on Almost Commuting Matrices.} Proc. R. Ir. Acad. Vol. 81A(1), 43-47 (1981)

\bibitem{Berg} I. Berg. {\it On Approximation of Normal Operators by Weighted Shifts.} Michigan Math. J. 21(4): 377-383 (July 1975).

\bibitem{Bhatia Perturbation} R. Bhatia, C. Davis, and A. McIntosh. {\it Perturbation of Spectral Subspaces and Solution of Linear Operator Equations}. Linear Algebra and its Applications. 52/53:45-67 (1983) \url{https://doi.org/10.1016/0024-3795(83)80007-X}

\bibitem{Extremal} R. Bhatia, C.Davis and P. Koosis. {\it An extremal problem in Fourier analysis with applications to operator theory.} Journal of Functional Analysis 82, 138-150 (1989)\\ \url{https://doi.org/10.1016/0022-1236(89)90095-5}

\bibitem{Bhatia} R. Bhatia and P. Rosethal. {\it How and Why to Solve the Operator Equation $AX - XB = Y.$} Bull. London Math. Soc. 29 (1997) 1-21 \url{ https://doi.org/10.1112/S0024609396001828} 

\bibitem{Blackadar} B. Blackadar. {\it Operator Algebras - Theory of $C^\ast$-Algebras and von Neumann Algebras.} Springer-Verlag. Berlin Heidelberg. 2006.

\bibitem{B&R} O. Bratteli and D. Robinson. {\it Operator Algebras and Quantum Statistical Mechanics.} Volume 1. 1979.

\bibitem{Choi}
M.-D. Choi. {\it Almost Commuting Matrices Need not be Nearly Commuting.} Proceedings of American Mathematical Society. Vol 102. No 3. March 1988.

\bibitem{Davidson} K. Davidson. {\it Almost Commuting Hermitian Matrices.} Math. Scand. 56 (1985), 222-240. 

\bibitem{ByExample} K. Davidson. {\it $C^\ast$-Algebras by Example.} AMS. 1996.

\bibitem{Handbook} K. Davidson and S. Szarek ``Local Operator Theory, random matrices and Banach spaces'' {\it Handbook of the Geometry of Banach Spaces. Vol I.} Elsevier Science. 2001.

\bibitem{three projections} C. Davis. {\it Generators of the Ring of Bounded Operators.} Proc. Amer. Math. Soc. Vol. 6, No. 6 (Dec., 1955), pp.970-972.

\bibitem{Davis} C. Davis. {\it The Rotation of Eigenvectors by a Perturbation.} Journal of Mathematical Analysis and Applications 6, 159-173 (1963) \url{https://doi.org/10.1016/0022-247X(63)90001-5}

\bibitem{D&K} C. Davis and W. Kahan. {\it The Rotation of Eigenvectors by a Perturbation. III.} SIAM J. Numer. Anal., 7(1), 1–46. 1970. \url{https://doi.org/10.1137/0707001}

\bibitem{Exp} S. Demko, W. Moss, and P. Smith. {\it Decay Rates for Inverses of Band Matrices}. Mathematics of Computation, 43(168), 491-499. (1984) doi:10.2307/2008290

\bibitem{Dixmier} J. Dixmier. {\it Position relative de deux vari\'et\'es lin\'eaires ferm\'ees dans un espace de Hilbert.}

\bibitem{Dimension} D. Enders and T. Shulman. {\it Almost Commuting Matrices, Cohomology, and Dimension.} \href{https://arxiv.org/abs/1902.10451}{	arXiv:1902.10451}

\bibitem{F-R} P. Friis and M. R{\o}rdam. {\it Almost commuting self-adjoint matrices - a short proof of Huaxin Lin's theorem.} J. reine angew. Math 479 (1996). 121-131.

\bibitem{Glebsky} L. Glebsky. {\it Almost commuting matrices with respect to normalized Hilbert-Schmidt norm.} \href{https://arxiv.org/abs/1002.3082}	{arXiv:1002.3082}

\bibitem{Grafakos} L. Grafakos. {\it Classical Fourier Analysis}. 2nd Edition. 2008.

\bibitem{Halmos} P. Halmos. {\it Two Subspaces.} Trans. Amer. Math. Soc. 144 (1969), 381–389.

\bibitem{Halmosproblems}
P. Halmos. {\it Some unsolved problems of  unknown depth  about operators  on Hilbert space}. Proceedings  of  the Royal  Society  of Edinburgh, 76A, 67-76, 1976. \url{https://doi.org/10.1017/S0308210500019491}

\bibitem{Hastings1sted} M. Hastings. {\it Making Almost Commuting Matrices Commute.} Commun. Math. Phys. 291, 321–345 (2009)

\bibitem{Hastings} M. Hastings. {\it Making Almost Commuting Matrices Commute.} Version 4. \href{https://arxiv.org/abs/0808.2474}{arXiv:0808.2474}

\bibitem{PC} M. Hastings. Personal Communication.

\bibitem{HastingsLoring} M. Hastings and T. Loring. {\it Almost commuting matrices, localized Wannier functions, and the quantum Hall effect.} J. Math. Phys. 51, 015214 (2010). \url{https://doi.org/10.1063/1.3274817}

\bibitem{OgataC} D. Herrera. {\it Constructing Nearby Commuting Matrices
for Reducible Representations of $su(2)$
with an Application to Ogata’s Theorem.} Preprint (to appear on arXiv.org).

\bibitem{Jordan} C. Jordan. {\it Essai sur la g\'eom\'eacutrie \`a $n$ dimensions.}

\bibitem{KS} I. Kachkovskiy and Y. Safarov. {\it Distance to Normal Elements in $C^\ast$-Algebras of Real Rank Zero.} Journal of the American Mathematical Society, 2016-01, Vol.29 (1), p.61-80 

\bibitem{Lin} H. Lin.  {\it Almost commuting self-adjoint matrices and applications}. Fields. Inst. Commun. 13, 193 (1995).

\bibitem{LS} T. Loring and A. S{\o}rensen. {\it Almost-Commuting Self-Adjoint Matrices - The Real and Self-Dual Cases.} Reviews in Mathematical Physics, 28(07): 1650017, 2016.

\bibitem{Lux} W. Luxemburg and R. Taylor. {\it Almost Commuting Matrices are near Commuting Matrices}. Indagationes Mathematicae (Proceedings).
Volume 73, 1970, pp 96-98

\bibitem{QMA} D. Nagaj,  P. Wocjan,  Y. Zhang. {\it Fast Amplification of QMA.} \url{https://arxiv.org/abs/0904.1549}

\bibitem{Ogata} Y. Ogata. {\it Approximating macroscopic observables in quantum spin systems with commuting matrices.} Journal of Functional Analysis. Vol. 264, Issue 9, 1 May 2013, pp 2005-2033.

\bibitem{Tobias} T. Osborne. {\it Almost Commuting Unitaries with Spectral Gap are Near Nearly Commuting Matrices.}  Proc. Amer. Math. Soc. Vol 137, No 12, December 2009, pp4043-4048

\bibitem{P-S} C. Pearcy and A. Shields. {\it Almost Commuting Matrices.} Journal of Functional Analysis. 33, 332-338 (1979).

\bibitem{Pironio} S. Pironio. et al. {\it Device-independent quantum key distribution secure against collective attacks.} 2009 New J. Phys.11 045021

\bibitem{Quantum Algorithms} A. Prakash. {\it Quantum Algorithms for Linear Algebra and Machine Learning.} Dissertation. \url{https://www2.eecs.berkeley.edu/Pubs/TechRpts/2014/EECS-2014-211.html}

\bibitem{Rosenthal} P. Rosenthal. {\it Are Almost Commuting Matrices Near Commuting Matrices?} The American Mathematical Monthly, Vol. 76, No. 8 (Oct., 1969), pp. 925-926.

\bibitem{Said} M. Said. {\it Almost Commuting Elements in Non-Commutative Symmetric Operator Spaces}. Dissertation.

\bibitem{Stein} E. Stein and R. Shakarchi. {\it Complex Analysis.} Princeton University Press. 2003.

\bibitem{Szarek} S. Szarek. {\it On Almost Commuting Hermitian Operators.} Rocky Mountain Journal of Mathematics. Vol. 20, No. 2, Spring 1990.

\bibitem{Voiculescu} D. Voiculescu. {\it Asymptotically commuting finite rank unitary operators without commuting approximants.} Acta Sci. Math. (Szeged) 45:1-4(1983), 429-431 1983.


\end{thebibliography}
\end{document}